\pdfoutput=1
\RequirePackage{ifpdf}
\ifpdf % We are running pdfTeX in pdf mode
\documentclass[pdftex]{sigma}
\else
\documentclass{sigma}
\fi

\usepackage{euscript}
\usepackage[cmtip,all]{xy}

\newlength{\whitecirclewidth}
\newlength{\whitecircleheight}
\newlength{\blackcirclewidth}
\newlength{\blackcircleheight}
\newlength{\corssandcircleheight}
\newlength{\corssandcirclewidth}
\newlength{\nodesize}
\newcommand{\whitecirclesymbol}{\makebox[\whitecirclewidth]{\Large$\circ$}}

\newcommand{\blackcirclesymbol}{\makebox[\blackcirclewidth]{\Large$\bullet$}}

\newcommand{\crossandcirclesymbol}{\makebox[\corssandcirclewidth]{\large$\otimes$}}

\newcommand{\abovewnode}[1]%
{
*-<\nodesize>{\raisebox{0pt}[\whitecircleheight][0pt]{$\overset{\rlap{$\displaystyle
#1$}}{\whitecirclesymbol}$}}}

\newcommand{\belowwnode}[1]%
{
*-<\nodesize>{\raisebox{0pt}[\whitecircleheight][0pt]{$\underset{\rlap{$\displaystyle
#1$}}{\whitecirclesymbol}$}}}

\newcommand{\abovebnode}[1]%
{
*-<\nodesize>{\raisebox{0pt}[\whitecircleheight][0pt]{$\overset{\rlap{$\displaystyle
#1$}}{\blackcirclesymbol}$}}}

\newcommand{\belowbnode}[1]%
{
*-<\nodesize>{\raisebox{0pt}[\whitecircleheight][0pt]{$\underset{\rlap{$\displaystyle
#1$}}{\blackcirclesymbol}$}}}

\newcommand{\abovecnode}[1]%
{
*-<\nodesize>{\raisebox{0pt}[\whitecircleheight][0pt]{$\overset{\rlap{$\displaystyle
#1$}}{\crossandcirclesymbol}$}}}

\newcommand{\belowcnode}[1]%
{
*-<\nodesize>{\raisebox{0pt}[\whitecircleheight][0pt]{$\underset{\rlap{$\displaystyle
#1$}}{\crossandcirclesymbol}$}}}

\renewcommand{\ge}{\epsilon}
\numberwithin{equation}{section}
\numberwithin{theorem}{section}
\numberwithin{proposition}{section}
\numberwithin{lemma}{section}
\numberwithin{corollary}{section}
\numberwithin{definition}{section}
\numberwithin{example}{section}
\numberwithin{remark}{section}
\numberwithin{note}{section}
\newtheorem{Thm}{Theorem}[section]
\newtheorem{Cor}[Thm]{Corollary}
\newtheorem{Lem}[Thm]{Lemma}
\newtheorem{Prop}[Thm]{Proposition}
\newtheorem{Cons}[Thm]{Consequence}
{\theoremstyle{definition}
\newtheorem{Def}[Thm]{Definition}

\newtheorem{Rem}[Thm]{Remark}
}

\newcommand{\IP}[2]{\langle#1 , #2\rangle}     %inner product.

\DeclareMathOperator{\pr}{pr}
\DeclareMathOperator{\Id}{Id}
\DeclareMathOperator{\ad}{ad}
\DeclareMathOperator{\Hom}{Hom}
\DeclareMathOperator{\Sym}{Sym}
\DeclareMathOperator{\Span}{span}
\DeclareMathOperator{\diag}{diag}
\DeclareMathOperator{\Tr}{Tr}

\begin{document}

\allowdisplaybreaks

\renewcommand{\PaperNumber}{008}

\FirstPageHeading

\ShortArticleName{Systems of Dif\/ferential Operators and Generalized Verma Modules}

\ArticleName{Systems of Dif\/ferential Operators\\
and Generalized Verma Modules}

\Author{Toshihisa KUBO} 
\AuthorNameForHeading{T.~Kubo}

\Address{Graduate School of Mathematical Sciences, The University of Tokyo,\\
3-8-1 Komaba, Meguro-ku, Tokyo 153-8914, Japan}

\Email{\href{mailto:toskubo@ms.u-tokyo.ac.jp}{toskubo@ms.u-tokyo.ac.jp}}
\URLaddress{\url{http://www.ms.u-tokyo.ac.jp/~toskubo/}}

\ArticleDates{Received April 08, 2013, in f\/inal form January 17, 2014; Published online January 24, 2014}

\Abstract{In this paper we close the cases that were left open in our earlier works on the study of
conformally invariant systems of second-order dif\/ferential operators for degenerate principal series.
More precisely, for these cases, we f\/ind the special values of the systems of dif\/ferential operators,
and determine the standardness of the homomorphisms between the generalized Verma modules, that come from
the conformally invariant systems.}


\Classification{22E46; 17B10; 22E47}

\section{Introduction}
%\label{chap:intro}

Conformally invariant systems are systems of dif\/ferential operators that are equivariant under an action
of a~Lie algebra.
More precisely, let $\EuScript{V} \to M$ be a~vector bundle over a~smooth mani\-fold~$M$ and ${\mathfrak
g}_0$ a~Lie algebra of f\/irst-order dif\/ferential operators acting on smooth sections on $\EuScript{V}$.
A~linearly independent list $D_1,\ldots, D_m$ of dif\/ferential operators on~$\EuScript{V}$ is then said to
be \emph{conformally invariant} if, for each $X \in {\mathfrak g}_0$, the condition
\begin{gather*}
[X,D_i]\in\sum^n_{j=1}C^\infty(M)D_j
\end{gather*}
holds for all $1\leq i \leq n$, where $[X, D_i] = XD_i - D_i X$.
(For the precise def\/inition see, for example, Section~\ref{SS:Prelim} of~\cite{BKZ09}.) In this article we shall
take $M=\bar{N}_0Q_0/Q_0$, an open dense submanifold of a~certain class of f\/lag manifolds $G_0/Q_0$, and
take vector bundle $\EuScript{V} \to M$ to be the restriction of a~homogeneous line bundle $\mathcal{L}_s
\to G_0/Q_0$ to $\bar{N}_0Q_0/Q_0$.

Many examples of conformally invariant systems implicitly or explicitly exist in the literature, especially
in the case of $m=1$.
The Laplacian $\Delta$ on $\mathbb{R}^n$ and wave operator $\square$ on the Minkowski space
$\mathbb{R}^{3,1}$ are, for instance, two outstanding examples.
The Yamabe operator on $S^n$ is another example in conformal geometry.
(See, for instance, the series of works~\cite{KO03a, KO03b}, and~\cite{KO03c} of Kobayashi--{\O}rsted.) In
the case of $m \geq 2$, an example may f\/ind in a~work of Davidson--Enright--Stanke~\cite{DES91}.
It would be also interesting to point out that such systems of operators are implicitly presented in the
work~\cite{Wallach79} of Wallach and also related to a~project of Dobrev on constructing intertwining
dif\/ferential operators of arbitrary order
(see, for example,~\cite{Dobrev-a12,Dobrev-a13}, and~\cite{Dobrev88}).

The theory of conformally invariant systems consisting of one dif\/ferential operator started with a~work
of Kostant~\cite{Kostant75}.
He called such dif\/ferential operators \emph{quasi-invariant}.
The notion of conformally invariant systems generalizes that of quasi-invariant dif\/ferential operators,
and it agrees with the notion of conformal invariance introduced by Ehrenpreis~\cite{Ehrenpreis88}.
As a~generalization of quasi-invariant dif\/ferential operators conformally invariant systems are also
related with a~work of Huang~\cite{Huang93} for intertwining dif\/ferential operators.

As described above the theory of conformally invariant systems can be viewed as a~geometric-analytic theory.
Further, it is also closely related to the explicit construction of homomorphisms between generalized Verma
modules.
Homomorphisms between generalized Verma modules (equivalently, intertwining dif\/ferential operators
between degenerate principal series representations) have received a~lot of attentions from many points of
view.
See, for instance,~\cite{BE89, Boe85, BC86, Jakobsen85,KKR03, Juhl09,KOSS-a13, KP-a13, KR00, Lepowsky77,
Matumoto06,Matumoto-a12, Somberg-a13}, and the references therein.
We wish to note that the list of the articles is far from complete.
We hope that the reader will understand that it is impossible to attempt to give an exhaustive list, as
there are numerous relevant contributions in the literature.

In~\cite{Kostant75}, Kostant showed that quasi-invariant dif\/ferential operators explicitly yield
homomorphisms between appropriate generalized Verma modules.
Conformally invariant systems also yield concrete homomorphisms between appropriate generalized Verma
modules~\cite{BKZ09}.
We shall describe our work for this direction more carefully later in this introduction.

We now turn to a~project on conformally invariant systems.
A project for conformally invariant systems started with the work of~\cite{BKZ08} and~\cite{BKZ09}.
Other progress on the project are reported in~\cite{Kable11, Kable12A, Kable12C, Kable12B, Kable13}
and~\cite{ Kubo11, KuboThesis2, KuboThesis1}.
To see a~recent development of the theory of conformally invariant systems the introduction of the
work~\cite{Kable13} of Kable is very helpful.
Besides the recent development one can also f\/ind relationships between quasi-invariant dif\/ferential
operators and the classic work of Maxwell on harmonic polynomials in the introduction.
Descriptions on Heisenberg Laplacian and Heisenberg ultrahyperbolic equation may also deserve one's
attention.
(In the particular direction the introduction of~\cite{Kable12C} has more details.)

The present work is also part of the project.
The aim of this paper is to close the cases that were left open in~\cite{KuboThesis1}
and~\cite{KuboThesis2}.
To describe our work of this paper more precisely, we now brief\/ly review the works in the papers.
Let $G$ be a~complex, simple, connected, simply-connected Lie group with Lie algebra ${\mathfrak g}$.
Give a~$\mathbb{Z}$-grading ${\mathfrak g} = \bigoplus_{j=-r}^r {\mathfrak g}(j)$ on ${\mathfrak g}$ so
that ${\mathfrak q} = {\mathfrak g}(0) \oplus \bigoplus_{j>0} {\mathfrak g}(j)$
$={\mathfrak l}\oplus{\mathfrak n}$ is a~maximal parabolic subalgebra.
Let $Q = N_G({\mathfrak q}) = LN$.
For a~real form ${\mathfrak g}_0$ of ${\mathfrak g}$, def\/ine $G_0$ to be an analytic subgroup of $G$ with
Lie algebra ${\mathfrak g}_0$.
Set $Q_0 = N_{G_0}({\mathfrak q})$.
We consider a~line bundle $\mathcal{L}_{s} \to G_0/Q_0$ for $s\in \mathbb{C}$.
As the homogeneous space $G_0/Q_0$ admits an open dense submanifold $\bar{N}_0Q_0/Q_0$, we restrict our
bundle to this submanifold.
By slight abuse of notation we refer to the restricted bundle as $\mathcal{L}_s$.
The systems that we shall construct act on smooth sections of the restricted bundle $\mathcal{L}_s \to
\bar{N}_0$.

Our systems of operators are constructed from $L$-irreducible constituents $W$ of ${\mathfrak g}(-r+k)
\otimes {\mathfrak g}(r)$ for $1\leq k \leq 2r$.
We call the systems of operators \emph{$\Omega_k$ systems}.
(We shall describe the construction more precisely in Section~\ref{SS:Omega}.)
An $\Omega_k$ system is a~system of $k$th-order dif\/ferential operators.
There is no reason to expect that $\Omega_k$ systems are conformally invariant on $\mathcal{L}_{s}$ for
arbitrary $s \in \mathbb{C}$; the conformal invariance of $\Omega_k$ systems depends on the complex
parameter $s$ for the line bundle~$\mathcal{L}_{s}$.
We then say that an $\Omega_k$ system has \emph{special value $s_k$} if the system is conformally invariant
on the line bundle $\mathcal{L}_{s_k}$.

In~\cite{KuboThesis1}, the parabolic subalgebra ${\mathfrak q}={\mathfrak l} \oplus {\mathfrak n}$ was
taken to be a~maximal parabolic subalgebra of \emph{quasi-Heisenberg type}, that is, a~maximal parabolic
subalgebra with nilpotent radical ${\mathfrak n}$ with conditions that $[{\mathfrak n}, [{\mathfrak n},
{\mathfrak n}]] = 0$ and $\dim([{\mathfrak n},{\mathfrak n}]) > 1$.
In this setting the $\Omega_k$ systems for $k\geq 5$ are zero.
We then sought the special values for the $\Omega_1$ system and $\Omega_2$ systems.
While the special value $s_1$ for the $\Omega_1$ system was determined for each parabolic subalgebra
${\mathfrak q}$ as $s_1 = 0$, we left three cases open for $\Omega_2$ systems.
To describe the open cases, let us now start explaining brief\/ly the classif\/ication for the irreducible
constituents $W$ that contribute to $\Omega_2$ systems.
In~\cite{KuboThesis1}, we f\/irst observed that if irreducible constituents $W$ contribute to $\Omega_2$
systems then their highest weights are of the form $\mu+\ge$, where $\mu$ is the highest weight for
${\mathfrak g}(1)$ and $\ge$ is some weight for~${\mathfrak g}(1)$.
We called such irreducible constituents \emph{special} and classif\/ied as Type~1a, Type~1b, Type~2, and
Type~3 with respect to certain technical conditions for the highest weight $\mu+\ge$.
(The precise conditions will be given in Def\/inition~\ref{Def:Type}.) Table~\ref{table-1} summarizes the
types of special constituents.
Here, for example, ``$B_n(i)$'' indicates the maximal standard parabolic subalgebra of ${\mathfrak g}$ of
type~$B_n$, which is determined by the $i$th simple root $\alpha_i$.
In the table $V(\mu+\epsilon_\gamma)$, $V(\mu+\epsilon_{n\gamma})$, and $V(\mu+\epsilon_{n\gamma}^\pm)$
denote the special constituents with highest weights $\mu+\epsilon_\gamma$, $\mu+\epsilon_{n\gamma}$, and
$\mu+\epsilon_{n\gamma}^\pm$, respectively.
(We shall precisely describe these highest weights in Section~\ref{SS:Omega_2} and
Section~\ref{SS:D(n-2)}).
As illustrated in Table~\ref{table-1}, there are one, two, or three special constituents.
A dash in the column for $V(\mu + \epsilon_{n\gamma})$ indicates that there is no special constituent
$V(\mu+\epsilon_{n\gamma})$ in the case.

In~\cite{KuboThesis1}, under the assumption that ${\mathfrak q}$ is not of type~$D_n(n-2)$, we found the
special values $s_2$ for $\Omega_2$ systems for the Type~1a and Type~2 constituents.
The technique that we used allowed us to handle each case uniformly.
However, since the technique relied on some technical conditions on the highest weights, we could not apply
it to the systems coming from Type~1b and Type~3 constituents.

\begin{table}[h]
\caption{Types of special constituents.}\label{table-1}\vspace{1mm}
\centering
\begin{tabular}{ccccc}
\hline
Parabolic subalgebra &$V(\mu+\epsilon_\gamma)$ &$V(\mu+\epsilon_{n\gamma})$
\\
\hline
$B_n(i),\quad 3\leq i \leq n-2$ &Type 1a&Type 1a
\\
$B_n(n-1)$ &Type 1a&\fbox{Type 1b}
\\
$B_n(n)$ &Type 2&$-$
\\
$C_n(i), \quad 2\leq i \leq n-1$ &\fbox{Type 3}&Type 2
\\
$D_n(i), \quad 3\leq i \leq n-3$ &Type 1a&Type 1a
\\
$E_6(3)$ &Type 1a&Type 1a
\\
$E_6(5)$ &Type 1a&Type 1a
\\
$E_7(2)$ &Type 1a&$-$
\\
$E_7(6)$ &Type 1a&Type 1a
\\
$E_8(1)$ &Type 1a&$-$
\\
$F_4(4)$ &Type 2&$-$
\\
\hline
\end{tabular}
\end{table}

\begin{table}[h]
\centering
\begin{tabular}{c|c|c|c}
\hline
Parabolic subalgebra &$V(\mu+\epsilon_\gamma)$ &$V(\mu+\ge^+_{n\gamma})$ &$V(\mu+\ge^-_{n\gamma})$
\\
\hline
\fbox{$D_n(n-2)$} & Type 1a & Type 1a & Type 1a\tsep{4pt}\bsep{4pt}
\\
\hline
\end{tabular}
%\label{table-2}
\end{table}

If $V(\mu+\ge):=V(\mu+\epsilon_\gamma)$, $V(\mu+\epsilon_{n\gamma})$, or $V(\mu+\epsilon_{n\gamma}^\pm)$
then the missing cases may be summarized as follows:
\begin{enumerate}\itemsep=0pt
\item[1)] the maximal parabolic subalgebra ${\mathfrak q}$ is of type~$D_n(n-2)$,
\item[2)] the special constituent
$V(\mu+\ge)$ is of Type~1b, and
\item[3)] the special constituent $V(\mu+\ge)$ is of Type~3.
\end{enumerate}
These are the cases boxed in Table~\ref{table-1}.
Our goal in this article is to f\/ind the special values~$s_2$ of~$\Omega_2$ systems in these three cases.
In contrast to~\cite{KuboThesis1}, in which each case was treated as uniformly as possible, we handle the
three cases individually in this paper.
For the case that ${\mathfrak q}$ is of type~$D_n(n-2)$, we f\/irst observe that each special constituent
is of Type~1a.
We then directly apply the technique used in~\cite{KuboThesis1}.
For Type~1b and Type~3 cases, we use an explicit realization of a~Lie algebra ${\mathfrak g}$.
In this way certain computations can be carried out easily.

The special values $s_2$ for the type~$D_n(n-2)$, Type~1b, and Type~3 cases are $s_2=1$
(Corollary~\ref{Cor:Special}), $s_2 = 1$ (Theorem~\ref{Thm:SpValType1b}), and $s_2=n-i+1$
(Theorem~\ref{Thm:SpValType3}), respectively.
Now, with these results together with ones in~\cite{KuboThesis1}, if $\Delta$ and $\Delta({\mathfrak
g}(1))$ denote a~(f\/ixed) root system and the set of roots contributing to ${\mathfrak g}(1)$,
respectively, then we obtain the following beautiful consequence.

\begin{Cons}
Let ${\mathfrak q}$ be a~maximal parabolic subalgebra of quasi-Heisenberg type.
The special value $s_2$ of the $\Omega_2$ system associated to the special constituent $V(\mu+\ge)$ is
\begin{gather*}
s_2=
\begin{cases}
\displaystyle{\frac{|\Delta_{\mu+\ge}({\mathfrak g}(1))|}{2}-1} & \text{if}\  V(\mu+\ge)\quad\text{is of Type~1},
\\
-1 & \text{if}\quad V(\mu+\ge)\ \text{is of Type~2},
\\
n-i+1 & \text{if}\quad V(\mu+\ge)\ \text{is of Type~3},
\end{cases}
\end{gather*}
where $|\Delta_{\mu+\ge}({\mathfrak g}(1))|$ is the number of elements of
$\Delta_{\mu+\ge}({\mathfrak g}(1)):= \{\alpha \in \Delta({\mathfrak g}(1)) \, | \, (\mu+\ge) -\alpha \in
\Delta\}$.
\end{Cons}

Here we combine Type~1b with Type~1a.
This is because it turned out that the special value for Type~1b case can be given by the same formula
as for  Type~1a case
(see Remark~\ref{Rem:SpValType1b}). If $\Omega_2|_{V(\mu+\ge)^*}$ denotes the $\Omega_2$ system coming from
the special constituent $V(\mu+\ge)$ then Table~\ref{table-3} exhibits the special values for all the
$\Omega_2$ systems under consideration.

\begin{table}[h] \caption{Line bundles with special values.}\label{table-3}\vspace{1mm}
\centering
\begin{tabular}{c|c|c}
\hline
Parabolic subalgebra &$\Omega_2|_{V(\mu+\epsilon_\gamma)^*}$ & $\Omega_2|_{V(\mu+\epsilon_{n\gamma})^*}$
\\
\hline
$B_n(i),\quad 3\leq i \leq n-2$ & $\mathcal{L}\big( (n- i - \frac{1}{2})\lambda_i\big)$ &
$\mathcal{L}(\lambda_i)$
\\
$B_n(n-1)$ & $\mathcal{L}\big( \frac{1}{2}\lambda_{n-1} \big)$ & $\mathcal{L}(\lambda_{n-1})$
\\
$B_n(n)$ & $\mathcal{L}(-\lambda_n)$ & $-$
\\
$C_n(i),\quad 2 \leq i \leq n-1$ &$\mathcal{L}\big((n - i + 1)\lambda_i\big)$ & $\mathcal{L}(-\lambda_i)$
\\
$D_n(i),\quad 3 \leq i \leq n-3$ & $\mathcal{L}\big((n - i - 1)\lambda_i\big)$& $\mathcal{L}(\lambda_i)$
\\
$E_6(3)$ & $\mathcal{L}(\lambda_3)$ & $\mathcal{L}(2\lambda_3)$
\\
$E_6(5)$ & $\mathcal{L}(\lambda_5)$ & $\mathcal{L}(2\lambda_5)$
\\
$E_7(2)$ & $\mathcal{L}(2\lambda_2)$ & $-$
\\
$E_7(6)$ & $\mathcal{L}(\lambda_6)$ & $\mathcal{L}(3\lambda_6)$
\\
$E_8(1)$ & $\mathcal{L}(3\lambda_1)$ & $-$
\\
$F_4(4)$ & $\mathcal{L}(-\lambda_4)$ & $-$
\\
\hline
\end{tabular}
\end{table}
\begin{table}[h]
\centering
\begin{tabular}
{c|c|c|c}
\hline
Parabolic subalgebra &$\Omega_2|_{V(\mu+\ge_{\gamma})^*}$ &
$\Omega_2|_{V(\mu+\ge^+_{n\gamma})^*}$&$\Omega_2|_{V(\mu+\ge^-_{n\gamma})^*}$
\\
\hline
$D_n(n-2)$ & $\mathcal{L}(\lambda_{n-2})$ & $\mathcal{L}(\lambda_{n-2})$ & $\mathcal{L}(\lambda_{n-2})$
\\
\hline
\end{tabular}
%\label{table-4}
\end{table}

Here $\lambda_i$ denotes the fundamental weight for the simple root $\alpha_i$ and
$\mathcal{L}(s\lambda_i):=\mathcal{L}_s$, the restricted line bundle over $\bar{N}_0$.

Now we turn to~\cite{KuboThesis2}.
To describe the work of the paper, we f\/irst recall that Kostant showed in~\cite{Kostant75} that
a~quasi-invariant dif\/ferential operator gives a~homomorphism between suitable scalar generalized Verma
modules with explicit image of a~highest weight vector.
As a~full generalization of quasi-invariant dif\/ferential operators, it is then shown in~\cite{BKZ09} that
a~conformally invariant system also explicitly yields a~homomorphism between appropriate generalized Verma
modules.
Here the generalized Verma modules are not necessarily of scalar-type.

A homomorphism between generalized Verma modules is called \emph{standard} if it comes from a~homomorphism
between corresponding (full) Verma modules, and called \emph{non-standard} other\-wise~\cite{Lepowsky77}.
While standard homomorphisms are well-understood~\cite{Boe85,Lepowsky77}, the classif\/ication of
non-standard homomorphisms is still an open problem.

In~\cite{KuboThesis2}, we classif\/ied the standardness of the homomorphisms $\varphi_{\Omega_k}$ between
generalized Verma modules arising from the conformally invariant $\Omega_1$ and $\Omega_2$ systems.
While the map $\varphi_{\Omega_1}$ was shown to be standard for each parabolic subalgebra~${\mathfrak q}$,
the classif\/ication for the map~$\varphi_{\Omega_2}$ was incomplete.
This is because of the lack of the special values of the three cases mentioned above (the boxed cases in
Table~\ref{table-1}).
Thus, in this paper, we also determine whether or not the maps~$\varphi_{\Omega_2}$ coming from the
conformally invariant $\Omega_2$ systems in the three cases are standard.
The classif\/ication results are given in Theorem~\ref{Thm:Map_D(n-2)} (type~$D_n(n-2)$ case),
Theorem~\ref{Thm:Map_Type1b} (Type~1b case), and Theorem~\ref{Thm:Map_Type3} (Type~3 case).
It turned out that in each case the map $\varphi_{\Omega_2}$ is non-standard.
With the results from~\cite{KuboThesis2}, Table~\ref{table-5} summarizes the classif\/ication of the
standardness of the maps~$\varphi_{\Omega_2}$.
\begin{table}[h] \caption{The classif\/ication of $\varphi_{\Omega_2}$.}\label{table-5}\vspace{1mm}
\centering
\begin{tabular}{c|c|c}
\hline
Parabolic subalgebra &$\Omega_2|_{V(\mu+\ge_{\gamma})^*}$ & $\Omega_2|_{V(\mu+\ge_{n\gamma})^*}$
\\
\hline
$B_n(i),\quad 3\leq i \leq n-2$ & standard & non-standard
\\
$B_n(n-1)$ & standard & non-standard
\\
$B_n(n)$ & standard & $-$
\\
$C_n(i),\quad 2 \leq i \leq n-1$ &non-standard & standard
\\
$D_n(i),\quad 3 \leq i \leq n-3$ & non-standard & non-standard
\\
$E_6(3)$ & non-standard & non-standard
\\
$E_6(5)$ & non-standard & non-standard
\\
$E_7(2)$ & non-standard & $-$
\\
$E_7(6)$ & non-standard & non-standard
\\
$E_8(1)$ & non-standard & $-$
\\
$F_4(4)$ & standard & $-$
\\
\hline
\end{tabular}
\end{table}
\begin{table}[h]
\centering
\begin{tabular}{c|c|c|c}
\hline
Parabolic subalgebra &$\Omega_2|_{V(\mu+\ge_{\gamma})^*}$ &
$\Omega_2|_{V(\mu+\ge^+_{n\gamma})^*}$&$\Omega_2|_{V(\mu+\ge^-_{n\gamma})^*}$
\\
\hline
$D_n(n-2)$ & non-standard & non-standard & non-standard
\\
\hline
\end{tabular}
%\label{table-6}
\end{table}

Recall that, for each parabolic subalgebra ${\mathfrak q}$, the special value $s_1$ for the $\Omega_1$
system is $s_1 = 0$ and that the map $\varphi_{\Omega_1}$ is standard.
Now, with the results and ones in Tables~\ref{table-1} and~\ref{table-5}, we obtain the following interesting
consequence.
\begin{Cons}
%\label{Cons:13}
Let ${\mathfrak q}$ be a~maximal parabolic subalgebra of quasi-Heisenberg type.
The map~$\varphi_{\Omega_k}$ for $k=1,2$ is non-standard if and only if the special value~$s_k$ of the
$\Omega_k$ system is a~positive integer.
\end{Cons}

Here we wish to note that there is another interesting relationship between the special values of
$\Omega_k$ systems and the associated scalar generalized Verma modules.
In~\cite{BKZ08}, when the nilpotent radical ${\mathfrak n}$ of parabolic subalgebra ${\mathfrak
q}={\mathfrak l} \oplus {\mathfrak n}$ is a~Heisenberg algebra, it was shown that all the f\/irst
reducibility points of the scalar generalized Verma modules associated with~$\Omega_k$ systems were
accounted by the special values of the systems.
In this case, as of quasi-Heisenberg case, the~$\Omega_k$ systems for $k \geq 5$ are zero.
To obtain the account the special values of~$\Omega_3$ systems and~$\Omega_4$ systems as well as these of
the~$\Omega_1$ system and~$\Omega_2$ systems were used.
In the quasi-Heisenberg case the results in~\cite{KuboThesis1} and in this paper for the~$\Omega_1$ system
and~$\Omega_2$ systems do not give an account for the f\/irst reducibility points.
The special values of $\Omega_3$ systems and~$\Omega_4$ systems thus seem to require also in the
quasi-Heisenberg case.
We will report the spacial values of the systems elsewhere.

Before closing this introduction let us make two remarks on this paper, although we understand that this
introduction is already long enough.
The f\/irst is on the technique to determine the special values of~$\Omega_2$ systems.
In~\cite{KuboThesis1}, to determine the special values, we used some reduction techniques.
Although the techniques signif\/icantly reduced the amount of computations, as several technical formulas
on dif\/ferential operators were used, the computations were still somewhat not straightforward.
Now it is known that special values can be obtained by computations in generalized Verma modules
(see, for instance,~\cite{Kable12C}). Then, in this paper, by combining the idea used in~\cite{KuboThesis1}
with that for generalized Verma modules, we are successful to simplify computations further.
The technique is given in Proposition~\ref{Prop:HL}.

The second remark is on the def\/inition of special constituents.
As stated in the introduction in~\cite{KuboThesis2}, there were certain discrepancy on the terminology
``special constituents'' between~\cite{KuboThesis1} and~\cite{KuboThesis2}.
In~\cite{KuboThesis1} (and in an earlier paragraph of this introduction) we def\/ined the special
constituents for $\Omega_2$ systems as ones whose highest weights satisfy certain conditions.
On the other hand, in~\cite{KuboThesis2}, we redef\/ined such constituents for any $\Omega_k$ systems as
irreducible constituents that contribute to $\Omega_k$ systems.
The reason why we redef\/ined special constituents as in~\cite{KuboThesis2} is that the later def\/inition
works not only for $\Omega_2$ systems but also for any $\Omega_k$ systems, and also that the irreducible
constituents with highest weights satisfying the technical conditions are highly expected to contribute to
$\Omega_2$ systems.
Indeed, in~\cite{KuboThesis1}, the implication was verif\/ied except the three missing cases.
In this paper we showed that the implication does hold also in the three cases.
The results are described in Section~\ref{SS:D(n-2)} for the type~$D_n(n-2)$ case, and stated in
Propositions~\ref{Prop:Tau2Type1b} and~\ref{Prop:Ttau} for the Type~1b and Type~3 cases, respectively.
Now the two notions of special constituents do agree for $\Omega_2$ systems.

We now outline the rest of this paper.
This paper consists of f\/ive sections (with this introduction) and one appendix.
In Section~\ref{SS:Prelim} we review the works~\cite{KuboThesis1} and~\cite{KuboThesis2}.
In particular we give the precise construction of $\Omega_k$ systems.
We also review about maximal parabolic subalgebras ${\mathfrak q}$ of quasi-Heiseberg type in this section.
In Section~\ref{SS:D(n-2)}, when ${\mathfrak q}$ is of type~$D_n(n-2)$, we f\/ind the special values of the
$\Omega_2$ systems and determine the standardness of the map $\varphi_{\Omega_2}$.
The special values are given in Corollary~\ref{Cor:Special} together with Theorem~\ref{Thm:DnSpecial}.
The standardness of $\varphi_{\Omega_2}$ is determined in Theorem~\ref{Thm:Map_D(n-2)}.
Sections~\ref{SS:Type1b} and~\ref{SS:Type3}
are devoted to the $\Omega_2$ systems arising from  Type~1b special constituent
and Type~3 special constituent, respectively.
The special values and the standardness of $\varphi_{\Omega_2}$ are determined in
Theorems~\ref{Thm:SpValType1b} and~\ref{Thm:Map_Type1b} for Type~1b case and Theorems~\ref{Thm:SpValType3}
and~\ref{Thm:Map_Type3} for Type~3 case, respectively.
Finally, in Appendix~\ref{SS:Data}, we collect the miscellaneous data that will be helpful for the work of
this paper.

\section{Preliminaries}
\label{SS:Prelim}

The purpose of this section is to summarize the framework established in~\cite{KuboThesis1}
and~\cite{KuboThesis2}.
The notation and conventions remain in force in the rest of this paper.

\subsection[A specialization of a~vector bundle $\EuScript{V} \to M$]{A specialization of a~vector bundle
$\boldsymbol{\EuScript{V} \to M}$}
%\label{SS:Setup}

We start with recalling from~\cite{KuboThesis1} the smooth manifold $M$ and the vector bundle $\EuScript{V}
\to M$ to work with.
This is nothing but the non-compact picture of a~degenerate principal series representation.

Let $G$ be a~complex, simple, connected, simply-connected Lie group with Lie algebra ${\mathfrak g}$.
Fix a~maximal connected solvable subgroup $B$, and write ${\mathfrak b} ={\mathfrak h} \oplus {\mathfrak
u}$ for its Lie algebra with ${\mathfrak h}$ the Cartan subalgebra and~${\mathfrak u}$ the nilpotent
subalgebra.
Let ${\mathfrak q} \supset {\mathfrak b}$ be a~standard parabolic subalgebra of~${\mathfrak g}$.
If $Q = N_G({\mathfrak q})$ then write $Q = LN$ for the Levi decomposition of~$Q$.
Let ${\mathfrak g}_0$ be a~real form of~${\mathfrak g}$ in which the complex parabolic subalgebra
${\mathfrak q}$ has a~real form~${\mathfrak q}_0$.
Let $G_0$ be the analytic subgroup of $G$ with Lie algebra~${\mathfrak g}_0$.
Def\/ine $Q_0 = N_{G_0}({\mathfrak q}) \subset Q$, and write $Q_0=L_0N_0$.
We will work with $G_0/Q_0$ for a~class of maximal parabolic subgroup $Q_0$ whose complexif\/ied Lie
algebra~${\mathfrak q}$ is a~maximal parabolic subalgebra of ${\mathfrak g}$ of quasi-Heisenberg type
(see Section~\ref{SS:Omega_2}).

Let $\Delta = \Delta({\mathfrak g},{\mathfrak h})$ denote the set of roots of ${\mathfrak g}$ with respect
to ${\mathfrak h}$.
Write $\Delta^+$ for the positive system attached to ${\mathfrak b}$ and denote by $\Pi$ the set of simple
roots.
We write ${\mathfrak g}_\alpha$ for the root space for $\alpha \in \Delta$.
For each subset $S \subset \Pi$, let ${\mathfrak q}_S$ be the corresponding standard parabolic subalgebra.
Write ${\mathfrak q}_S = {\mathfrak l}_S \oplus {\mathfrak n}_S$ with Levi factor ${\mathfrak l}_S =
{\mathfrak h} \oplus \bigoplus_{\alpha \in \Delta_S}{\mathfrak g}_\alpha$ and nilpotent radical ${\mathfrak
n}_S = \bigoplus_{\alpha \in \Delta^+ \backslash \Delta_S} {\mathfrak g}_\alpha$, where $\Delta_S = \{
\alpha \in \Delta \, | \, \alpha \in \Span(\Pi \backslash S) \}$.
If ${\mathfrak q}$ is a~maximal parabolic subalgebra of ${\mathfrak g}$ then there exists a~unique simple
root $\alpha_{\mathfrak q} \in \Pi$ so that ${\mathfrak q} = {\mathfrak q}_{\{\alpha_{\mathfrak q}\}}$.
Let $\lambda_{\mathfrak q}$ be the fundamental weight of~$\alpha_{\mathfrak q}$.
The weight $\lambda_{\mathfrak q}$ is orthogonal to any roots~$\alpha$ with ${\mathfrak g}_\alpha \subset
[{\mathfrak l}, {\mathfrak l}]$.
Hence it exponentiates to a~character~$\chi_{\mathfrak q}$ of~$L$.
For $s \in \mathbb{C}$, we set $\chi^{s}:=|\chi_{\mathfrak q}|^{s}$, an analytic character for~$L_0$.
As $\chi_{\mathfrak q}$ takes real values on $L_0$, we have $d\chi^s = s\lambda_{\mathfrak q}$.
Let $\mathbb{C}_{\chi^{s}}$ be the one-dimensional representation of~$L_0$ with character~$\chi^{s}$.
The representation $\chi^{s}$ is extended to a~representation of~$Q_0$ by making it trivial on~$N_0$.
It then deduces a~line bundle $\mathcal{L}_{s}$ on $ G_0/Q_0$ with f\/iber $\mathbb{C}_{\chi^{s}}$.

The group $G_0$ acts on the space
\begin{gather*}
C^\infty_{\chi}(G_0/Q_0,\mathbb{C}_{\chi^{s}})=\{F\in C^\infty(G_0,\mathbb{C}_{\chi^{s}}
)\,|\, F(gq)=\chi^{s}(q^{-1})F(g)~\text{for all}~q\in Q_0~\text{and}~g\in G_0\}
\end{gather*}
by left translation.
The action of ${\mathfrak g}_0$ on $C^\infty_\chi(G_0/Q_0, \mathbb{C}_{\chi^{s}})$ arising from this action
is given by
\begin{gather}
\label{Eqn:Action}
(Y{\raisebox{1pt}{$\scriptstyle\bullet$}}F)(g)=\frac{d}{dt}F(\exp(-tY)g)\bigg|_{t=0}
\end{gather}
for $Y \in {\mathfrak g}_0$, where the dot $ {\raisebox{1pt} {$\scriptstyle \bullet$} } $ denotes
the action of the Lie algebra as dif\/ferential operators.
This action is extended $\mathbb{C}$-linearly to ${\mathfrak g}$ and then naturally to the universal
enveloping algebra $\mathcal{U}({\mathfrak g})$.

Let $\bar{N}_0$ be the unipotent subgroup opposite to $N_0$.
The natural inf\/initesimal action of ${\mathfrak g}$ on the image of the restriction map
$C^\infty_{\chi}(G_0/Q_0, \mathbb{C}_{\chi^{s}}) \to C^\infty(\bar N_0, \mathbb{C}_{\chi^{s}})$ induced
by~\eqref{Eqn:Action} gives an action of ${\mathfrak g}$ on the whole space $C^\infty(\bar N_0,
\mathbb{C}_{\chi^{s}})$
(see, for instance, Section~\ref{SS:Prelim} of~\cite{KuboThesis1}).
The line bundle $\mathcal{L}_{s} \to G_0/Q_0$ restricted to $\bar{N}_0$ is the trivial bundle $\bar{N}_0
\times \mathbb{C}_{\chi^{s}} \to \bar{N}_0$.
We shall work with the trivial bundle on~$\bar{N}_0$.
By slight abuse of notation, we refer to the trivial bundle over~$\bar{N}_0$ as~$\mathcal{L}_{s}$.

\subsection[The $\Omega_k$ systems]{The $\boldsymbol{\Omega_k}$ systems}
\label{SS:Omega}

A systematic construction of systems of dif\/ferential operators for the line bundle $\mathcal{L}_s \to
\bar{N}_0$ was established in~\cite{KuboThesis1}.
In this subsection we summarize the construction of the systems of dif\/ferential operators.
For a~subspace $W$ of ${\mathfrak g}$, we write $\Delta(W) = \{\alpha \in \Delta \, | \, {\mathfrak
g}_\alpha \subset W\}$ and $\Pi(W) = \Delta(W) \cap \Pi$.
We keep the notation from the previous subsection, unless otherwise specif\/ied.

Let ${\mathfrak g} = \bigoplus_{j=-r}^r {\mathfrak g}(j)$ be a~$\mathbb{Z}$-grading on ${\mathfrak g}$ with
${\mathfrak g}(1) \neq 0$.
Take $L$ to be the analytic subgroup of $G$ with Lie algebra ${\mathfrak g}(0)$.
Observe that, as $[{\mathfrak g}(0),{\mathfrak g}(j)] \subset {\mathfrak g}(j)$, each graded subspace
${\mathfrak g}(j)$ is an $L$-module and so is ${\mathfrak g}(-r+k) \otimes {\mathfrak g}(r)$.
Write $R$ for the inf\/initesimal right translation of ${\mathfrak g}_0$.
As usual, we extend it $\mathbb{C}$-linearly to ${\mathfrak g}$ and then naturally to
$\mathcal{U}({\mathfrak g})$.

We build systems of dif\/ferential operators in three steps.
\begin{description}\itemsep=0pt
\item[{\emph{Step 1:}}]
First, for $1\leq k \leq 2r$, consider the $L$-equivariant polynomial map
\begin{gather*}
\tau_k : \ {\mathfrak g}(1)\to{\mathfrak g}(-r+k)\otimes{\mathfrak g}(r),\qquad
%\\
%\phantom{\tau_k:}
X\mapsto\big(\ad(X)^k\otimes\Id\big)\omega,
\end{gather*}
where $\omega$ is the element in ${\mathfrak g}(-r) \otimes {\mathfrak g}(r)$ def\/ined by
\begin{gather*}%\label{Eqn:Omega}
\omega:=\sum_{\gamma_j\in\Delta({\mathfrak g}(r))}X^*_{\gamma_j}\otimes X_{\gamma_j}.
\end{gather*}
Here $X_{\gamma_j}$ are root vectors for $\gamma_j$, and $X^*_{\gamma_j}$ are the vectors dual to
$X_{\gamma_j}$ with respect to the Killing form $\kappa$, namely, $X^*_{\gamma_j}(X_{\gamma_t}):=
\kappa(X^*_{\gamma_j}, X_{\gamma_t}) =\delta_{j,t}$ with $\delta_{i,t}$ the Kronecker delta.
\item[{\emph{Step 2:}}]
Next, for an $L$-irreducible constituent $W$ of ${\mathfrak g}(-r+k)
\otimes {\mathfrak g}(r)$, consider the associated $L$-intertwining operator $\tilde{\tau}_k|_{W^*} \in
\Hom_L(W^*, \mathcal{P}^k({\mathfrak g}))$ def\/ined by
\begin{gather*}
\tilde{\tau}_k|_{W^*}(Y^*)(X):=Y^*(\tau_k(X)),
\end{gather*}
where $W^*$ is the dual space for $W$ with respect to the Killing form $\kappa$.
Take an irreducible constituent $W$ of ${\mathfrak g}(-r+k) \otimes {\mathfrak g}(r)$ so that
$\tilde{\tau}_k|_{W^*} \not \equiv 0$.
\item[{\emph{Step 3:}}]
Last, to the space $W^*$ dual to the space $W$ taken in Step 2, apply the
following algebraic procedure:
\begin{gather}
\label{Eqn:Comp}
W^*\stackrel{\tilde{\tau}_k|_{W^*}}{\to}\mathcal{P}^k({\mathfrak g}(1))\cong\Sym^k({\mathfrak g}
(-1))\stackrel{\sigma}{\hookrightarrow}\mathcal{U}(\bar{\mathfrak n})\stackrel{R}{\to}\mathbb{D}(\mathcal{L}
_{s})^{\bar{\mathfrak n}}.
\end{gather}
Here, $\mathcal{P}^k({\mathfrak g}(1))$ is the space of homogeneous polynomials on ${\mathfrak g}(1)$ of
degree $k$, the map $\sigma: \Sym^k({\mathfrak g}(-1)) \to \mathcal{U}(\bar {\mathfrak n})$ is the
symmetrization operator, and $\mathbb{D}(\mathcal{L}_{s})^{\bar {\mathfrak n}}$ is the space of $\bar
{\mathfrak n}$-invariant dif\/ferential operators for $\mathcal{L}_{s}$.
\end{description}

Let $\Omega_k|_{W^*}: W^* \to \mathbb{D}(\mathcal{L}_{s})^{\bar {\mathfrak n}}$ be the composition of the
linear maps described in~\eqref{Eqn:Comp}, namely, $\Omega_k|_{W^*} = R \circ \sigma \circ
\tilde{\tau}_k|_{W^*}$.
For simplicity we write $\Omega_k(Y^*)=\Omega_k|_{W^*}(Y^*)$ for the dif\/ferential operator arising from
$Y^* \in W^*$.

Now, given basis $\{Y^*_1, \ldots, Y^*_m\}$ for $W^*$, we have a~system of dif\/ferential operators
\begin{gather*}
\Omega_k\big(Y^*_{1}\big),\ldots,\Omega_k(Y^*_{m}).
\end{gather*}
We call the system of operators the \textit{$\Omega_k|_{W^*}$ system}.
When the irreducible constituent $W^*$ is not important, we simply refer to each $\Omega_k|_{W^*}$ system
as an $\Omega_k$ system.
We may want to note that~$\Omega_k$ systems are independent of a~choice for a~basis for~$W^*$ up to some
natural equivalence
(see Def\/inition 3.5 of~\cite{KuboThesis1}).

The $\Omega_k$ systems are not conformally invariant for arbitrary complex parameter $s \in \mathbb{C}$ for
the line bundle $\mathcal{L}_s$.
We say that an $\Omega_k$ system has \textit{special value $s_k$} if the system is conformally invariant on
the line bundle $\mathcal{L}_{s_k}$.

\subsection[The $\Omega_k$ systems and generalized Verma modules]{The $\boldsymbol{\Omega_k}$ systems and
generalized Verma modules}
\label{SS:GVM-I}

Conformally invariant systems yield non-zero $\mathcal{U}({\mathfrak g})$-homomorphisms between appropriate
gene\-ra\-lized Verma modules.
Since the theory simplif\/ies computations to f\/ind the special values of~$\Omega_k$ systems, in this
subsection, we review how conformally invariant~$\Omega_k$ systems induce such homomorphisms.

We start with a~well-known fact on the duality between degenerate principal series and generalized Verma
modules.
A generalized Verma module $M_{\mathfrak q}[E]:=\mathcal{U}({\mathfrak g})\otimes_{\mathcal{U}({\mathfrak
q})}E$ is a~$\mathcal{U}({\mathfrak g})$-module that is induced from a~f\/inite dimensional simple
${\mathfrak q}$-module $E$.
It is well-known that if $\mathcal{E} \to G_0/Q_0$ is the homogenous bundle with f\/iber $E$ then there is
a~natural pairing between the space $\Gamma(\mathcal{E})$ of smooth sections and generalized Verma module
$M_{\mathfrak q}[E^*]$, where $E^*$ is the dual space for~$E$
(see, for example,~\cite{CS90} and~\cite{HJ82}). Via this natural pairing, associated to the line bundle~$\mathcal{L}_s$ is the generalized Verma module $M_{\mathfrak q}[\mathbb{C}_{-s}]$, where $\mathbb{C}_{-s}
= \mathbb{C}_{-s\lambda_{\mathfrak q}}$ is the ${\mathfrak q}$-module derived from the $Q_0$-representation
$(\chi^{-s},\mathbb{C})$ with $d\chi = \lambda_{\mathfrak q}$.

Now, given irreducible constituent $W$ of ${\mathfrak g}(-r+k) \otimes {\mathfrak g}(r)$, we def\/ine an
$L$-intertwining operator $\omega_k|_{W^*}: W^* \to \mathcal{U}(\bar {\mathfrak n})$ by $\omega_k|_{W^*} =
\sigma \circ \tilde{\tau}|_{W^*}$, so that $\Omega_k|_{W^*} = R \circ \omega_k|_{W^*}$.
If writing $\omega_k(W^*) = \omega_k|_{W^*}(W^*)$ then we obtain the following diagram:
\begin{gather}
\begin{split}
& \xymatrix @R=.1pc{
&&W^* \ar[ddd]^{\omega_k|_{W^*}} & \\
&&&&\\
&&&&\\
&M_{\mathfrak q}[\mathbb{C}_{-s}] & \ar[l]_{\quad \; \cdot \otimes\mathbb{C}_{-s}}
\mathcal{U}(\bar{\mathfrak n})\ar[r]^R&\mathbb{D}(\mathcal{L}_s)^{\bar{\mathfrak n}} \\
&\omega_k(W^*) \otimes\mathbb{C}_{-s} & \ar@{|->}[l] \omega_k(W^*)  \ar@{|~>}[r] &
\{\Omega_k(Y^*_1), \ldots,\Omega_k(Y^*_m)\},
}
\end{split}
\label{Eqn:Comp2}
\end{gather}
where $\{Y^*_1, \ldots, Y^*_m\}$ is a~given basis for $W^*$.
Here we indicate by the squiggly arrow $\rightsquigarrow$ that the right dif\/ferentiation $R$ is only
applied for the basis elements $\omega_k(Y^*_1), \ldots, \omega_k(Y^*_m)$ in $\omega_k(W^*)$.

It can be seen from~\eqref{Eqn:Comp2} that constructing the $\Omega_k|_{W^*}$ system is equivalent to
constructing the $L$-submodule $\omega_k(W^*)\otimes \mathbb{C}_{-s}$.
Proposition~\ref{Prop:CIS-GVM} below shows that there is a~further relationship.
\begin{Prop}[\protect{\cite[Theorem 19]{BKZ09}}]
\label{Prop:CIS-GVM}
The $\Omega_k|_{W^*}$ system is conformally invariant on the line bundle $\mathcal{L}_{s_k}$ if and only if
\begin{gather*}
\omega_k\big(W^*\big)\otimes\mathbb{C}_{-s_k}\subset M_{\mathfrak q}[\mathbb{C}_{-s_k}]^{{\mathfrak n}},
\end{gather*}
where $M_{\mathfrak q}[\mathbb{C}_{-s}]^{{\mathfrak n}} = \{ v \in M_{\mathfrak q}[\mathbb{C}_{-s}] \, | \,
X \cdot v = 0$  for all $X \in {\mathfrak n}\}$.
\end{Prop}

It follows from Proposition~\ref{Prop:CIS-GVM} that if the $\Omega_k|_{W^*}$ system is conformally
invariant on $\mathcal{L}_{s_k}$ then the $L$-submodule $\omega_k(W^*)\otimes \mathbb{C}_{-s_k}$ induces
a~$\mathcal{U}({\mathfrak g})$-module $\mathcal{U}({\mathfrak g})\otimes_{\mathcal{U}({\mathfrak q})}
\big(\omega_k(W^*)\otimes \mathbb{C}_{-s_k} \big)$.
The following proposition shows that this is indeed a~generalized Verma module.
\begin{Prop}[\protect{\cite[Proposition 3.4(1)]{KuboThesis2}}]\label{Prop:L-simple}
If $W^*$ has highest weight $\nu$ then $\omega_k(W^*)\otimes \mathbb{C}_{-s}$ is the simple $L$-submodule
of $M_{\mathfrak q}[\mathbb{C}_{-s}]$ with highest weight $\nu-s\lambda_{\mathfrak q}$.
\end{Prop}

Now it follows from Propositions~\ref{Prop:CIS-GVM} and~\ref{Prop:L-simple} that a~conformally invariant
$\Omega_k$ system induces a~homomorphism from $M_{\mathfrak q}[\omega_k(W^*)\otimes \mathbb{C}_{-s_k}]$ to
$M_{\mathfrak q}[\mathbb{C}_{-s_k}]$.
Indeed, if the $\Omega_k|_{W^*}$ system is conformally invariant for $\mathcal{L}_{s_k}$ then, by
Proposition~\ref{Prop:CIS-GVM}, there exists inclusion map $\iota \in \Hom_{L} \big(
\omega_k(W^*)\otimes \mathbb{C}_{-s_k}, M_{\mathfrak q}[\mathbb{C}_{-s_k}] \big)$.
The inclusion map $\iota$ then induces a~non-zero $\mathcal{U}({\mathfrak g})$-homomorphism
\begin{gather*}
\varphi_{\Omega_k}\in\Hom_{\mathcal{U}({\mathfrak g}),L}\big(M_{\mathfrak q}
[\omega_k(W^*)\otimes\mathbb{C}_{-s_k}],M_{\mathfrak q}[\mathbb{C}_{-s_k}]\big),
\end{gather*}
that is given by
\begin{gather}
M_{\mathfrak q}[\omega_k\big(W^*\big)\otimes\mathbb{C}_{-s_k}]\stackrel{\varphi_{\Omega_k}}{\to}
M_{\mathfrak q}[\mathbb{C}_{-s_k}],
\qquad
u\otimes\big(\omega_k(Y^*)\otimes1)\mapsto u\cdot\iota\big(\omega_k(Y^*)\otimes1).
\label{Eqn823}
\end{gather}

Observe that there is a~quotient map from a~(full) Verma module to a~generalized Verma module.
A homomorphism between generalized Verma modules is called \textit{standard} if it comes from
a~homomorphism between the corresponding full Verma modules, and called \textit{non-standard} otherwise~\cite{Lepowsky77}.
In the rest of this subsection, to study the standardness of the map $\varphi_{\Omega_k}$
in~\eqref{Eqn823}, we give a~simple criterion to determine whether or not the standard homomorphism
$\varphi_{\rm std}:M_{\mathfrak q}[\omega_k(W^*)\otimes \mathbb{C}_{-s_k}] \to M_{\mathfrak
q}[\mathbb{C}_{-s_k}]$ is zero.
To do so, it is convenient to parametrize generalized Verma modules by their inf\/initesimal characters.
Therefore we write
\begin{gather*}%\label{Eqn:InfChar2}
M_{\mathfrak q}[\mathbb{C}_{-s_k\lambda_{\mathfrak q}}]=M_{\mathfrak q}(-s_k\lambda_{\mathfrak q}+\rho),
\end{gather*}
where $\rho$ is half the sum of the positive roots.
Similarly, if $W^*$ has highest weight $\nu$ then, by Proposition~\ref{Prop:L-simple}, we write
\begin{gather*}%\label{Eqn:InfChar1}
M_{\mathfrak q}[\omega_k\big(W^*\big)\otimes\mathbb{C}_{-s_k}]=M_{\mathfrak q}(\nu-s_k\lambda_{\mathfrak q}+\rho).
\end{gather*}
Now, if $v: = \omega_k(Y^*) \otimes 1$ then~\eqref{Eqn823} is expressed by
\begin{gather}
M_{\mathfrak q}(\nu-s_k\lambda_{\mathfrak q}+\rho)\stackrel{\varphi_{\Omega_k}}{\to}M_{\mathfrak q}
(-s_k\lambda_{\mathfrak q}+\rho),
\qquad
u\otimes v \mapsto u\cdot\iota(v).
\label{Eqn:GVM2}
\end{gather}

To describe the criterion ef\/f\/iciently we recall a~well-known def\/inition for a~\emph{link} of two
weights.
Let $\IP{\cdot}{\cdot}$ be the inner product on ${\mathfrak h}^*$ induced from the Killing form $\kappa$.
Write $\alpha^{\vee} = 2\alpha/\IP{\alpha}{\alpha}$.
\begin{Def}[Bernstein--Gelfand--Gelfand]%\label{Def:Link}
Let $\lambda, \delta \in {\mathfrak h}^*$ and $\beta_1, \ldots, \beta_t \in
\Delta^+$.
Set $\delta_0 = \delta$ and $\delta_i = s_{\beta_i} \dots s_{\beta_1}\delta$ for $1 \leq i \leq t$.
We say that the sequence $(\beta_1, \ldots, \beta_t)$ \textit{links} $\delta$ to $\lambda$ if
\begin{enumerate}\itemsep=0pt
\item[(1)] $\delta_t = \lambda$ and
\item[(2)] $\IP{\delta_{i-1}}{\beta_i^{\vee}} \in \mathbb{Z}_{\geq 0}$
for $1\leq i \leq t$.
\end{enumerate}
\end{Def}

Let $M(\eta)$ denote the (full) Verma module with highest weight $\eta-\rho$.
As usual, if there is a~non-zero $\mathcal{U}({\mathfrak g})$-homomorphism from $M(\eta)$ into $M(\zeta)$
then we write $M(\eta) \subset M(\zeta)$.
If $\Pi({\mathfrak l})$ denotes the set of simple roots $\alpha \in \Pi$ so that ${\mathfrak g}_\alpha
\subset {\mathfrak l}$ then the criterion is given as follows.
\begin{Prop}[\protect{\cite[Proposition 4.6]{KuboThesis2}}]\label{Prop:Std} Let $M_{\mathfrak q}(\nu-s_k\lambda_{\mathfrak q} + \rho)$ and
$M_{\mathfrak q}(-s_k\lambda_{\mathfrak q} + \rho)$ be the generalized Verma modules in~\eqref{Eqn:GVM2}.
Then the standard map from $M_{\mathfrak q}(\nu-s_k\lambda_{\mathfrak q} + \rho)$ to $M_{\mathfrak
q}(-s_k\lambda_{\mathfrak q} + \rho)$ is zero if and only if there exists $\alpha_\nu \in \Pi({\mathfrak
l})$ so that $-\alpha_\nu - s_k\lambda_{\mathfrak q} +\rho$ is linked to $\nu-s_k\lambda_{\mathfrak q}
+\rho$.
\end{Prop}

\subsection[The $\Omega_2$ systems associated to maximal parabolic subalgebras of quasi-Heisenberg type]{The
$\boldsymbol{\Omega_2}$ systems associated to maximal parabolic subalgebras\\ of quasi-Heisenberg type}
\label{SS:Omega_2}

In Sections~\ref{SS:D(n-2)},~\ref{SS:Type1b}, and~\ref{SS:Type3}, we study $\Omega_2$ systems
associated to a~certain class of maximal parabolic subalgebra ${\mathfrak q} = {\mathfrak l} \oplus
{\mathfrak n}$, which we call quasi-Heisenberg type.
More precisely we f\/ind their special values and determine the standardness of the maps
$\varphi_{\Omega_2}$.
Then, in this subsection, we recall from~\cite{KuboThesis1} some observation on the $\Omega_2$ systems
associated to such maximal parabolic subalgebras.

First, we recall from~\cite{KuboThesis1} that a~parabolic subalgebra ${\mathfrak q} = {\mathfrak l} \oplus
{\mathfrak n}$ is called \textit{quasi-Heisenberg type} if its nilpotent radical ${\mathfrak n}$
satisf\/ies the conditions that $[{\mathfrak n}, [{\mathfrak n}, {\mathfrak n}]] = 0$ and $\dim([{\mathfrak
n}, {\mathfrak n}]) >1$.
Let~$\alpha_{\mathfrak q}$ be a~simple root, so that the maximal parabolic subalgebra ${\mathfrak
q}={\mathfrak q}_{\{\alpha_{\mathfrak q}\}} = {\mathfrak l} \oplus {\mathfrak n}$ determined by~$\alpha_{\mathfrak q}$ is of quasi-Heisenberg type.
Given Dynkin type~$\mathcal{T}$ of~${\mathfrak g}$, if we write $\mathcal{T}(i)$ for the Lie algebra
together with the choice of maximal parabolic subalgebra ${\mathfrak q} = {\mathfrak q}_{\{\alpha_i\}}$
determined by $\alpha_i$ then the maximal parabolic subalgebras ${\mathfrak q}={\mathfrak l} \oplus
{\mathfrak n}$ of quasi-Heisenberg type are classif\/ied as follows:
\begin{gather*}%\label{Eqn4.0.1}
B_n(i),\quad 3\leq i\leq n,
\qquad
C_n(i),\quad 2\leq i\leq n-1,
\qquad
D_n(i),\quad 3\leq i\leq n-2,
\end{gather*}
and
\begin{gather*}%\label{Eqn4.0.2}
E_6(3),\quad E_6(5),\quad E_7(2),\quad E_7(6),\quad E_8(1),\quad F_4(4).
\end{gather*}
Here, the Bourbaki conventions~\cite{Bourbaki08} are used for the labels of the simple roots.
Note that, in type~$A_n$, any maximal parabolic subalgebra has abelian nilpotent radical, and also that, in
type~$G_2$, the nilpotent radicals of two maximal parabolic subalgebras are a~3-step nilpotent or
Heisenberg algebra.

We next observe that a~maximal parabolic subalgebra ${\mathfrak q}$ of quasi-Heisenberg type induces
a~2-grading on ${\mathfrak g}$.
As ${\mathfrak q}$ has two-step nilpotent radical, if $\lambda_{\mathfrak q}$ is the fundamental weight for
$\alpha_{\mathfrak q}$ then, for all $\beta \in \Delta$, the quotient $2\IP{\lambda_{\mathfrak
q}}{\beta}/||\alpha_{\mathfrak q}||^2$ takes the values of $0$, $\pm 1$, or~$\pm 2$
(see for example Section~\ref{SS:Type1bA} of~\cite{KuboThesis1}).
Therefore, if $H_{\mathfrak q}$ is the element
in ${\mathfrak h}$ so that $\beta(H_{\mathfrak q}) = 2\IP{\lambda_{\mathfrak q}}{\beta}/||\alpha_{\mathfrak
q}||^2$ for all $\beta \in \Delta$, and if ${\mathfrak g}(j)$ is the $j$-eigenspace of
$\ad(H_{\mathfrak q})$ on ${\mathfrak g}$ then the adjoint action of $H_{\mathfrak q}$ induces
a~2-grading ${\mathfrak g} = \bigoplus_{j=-2}^2{\mathfrak g}(j)$ on ${\mathfrak g}$ with parabolic
subalgebra ${\mathfrak q} = {\mathfrak g}(0) \oplus {\mathfrak g}(1) \oplus {\mathfrak g}(2)$, where
${\mathfrak l} = {\mathfrak g}(0)$ and ${\mathfrak n} = {\mathfrak g}(1) \oplus {\mathfrak g}(2)$.
The subalgebra $\bar {\mathfrak n}$, the nilpotent radical opposite to ${\mathfrak n}$,
is given by $\bar{\mathfrak n} = {\mathfrak g}(-1) \oplus {\mathfrak g}(-2)$.
Here we have ${\mathfrak g}(2) = {\mathfrak z}({\mathfrak n})$ and ${\mathfrak g}(-2) = {\mathfrak
z}(\bar{\mathfrak n})$, where ${\mathfrak z}({\mathfrak n})$ (resp.\
${\mathfrak z}(\bar {\mathfrak n})$) is the center of ${\mathfrak n}$ (resp.\
$\bar {\mathfrak n}$).
Thus we denote the 2-grading on ${\mathfrak g}$ by
\begin{gather}
\label{Eqn4.1.6}
{\mathfrak g}={\mathfrak z}(\bar{\mathfrak n})\oplus{\mathfrak g}(-1)\oplus{\mathfrak l}\oplus{\mathfrak g}
(1)\oplus\frak{z}(\frak{n})
\end{gather}
with parabolic subalgebra
\begin{gather*}
%\label{Eqn4.1.7}
{\mathfrak q}={\mathfrak l}\oplus{\mathfrak g}(1)\oplus{\mathfrak z}({\mathfrak n}).
\end{gather*}

Now, for $1 \leq k \leq 4$, the maps $\tau_k$ associated to the grading~\eqref{Eqn4.1.6} are given by
\begin{gather*}%\label{Eqn:TauTwoStep}
\tau_k : \  {\mathfrak g}(1)\to{\mathfrak g}(-2+k)\otimes{\mathfrak z}({\mathfrak n}),
\qquad
X\mapsto\frac{1}{k!}\big(\ad(X)^k\otimes\Id\big)\omega
\end{gather*}
with
\begin{gather}
\label{Eqn:omega}
\omega=\sum_{\gamma_j\in\Delta({\mathfrak z}({\mathfrak n}))}X^*_{\gamma_j}\otimes X_{\gamma_j}.
\end{gather}
In particular when $k=2$, we have{\samepage
\begin{gather}
\tau_2 : \  {\mathfrak g}(1)\to{\mathfrak l}\otimes{\mathfrak z}({\mathfrak n}),
\qquad
X\mapsto\frac{1}{2}\big(\ad(X)^2\otimes\Id\big)\omega.
\label{Eqn:Tau2}
\end{gather}
The $\Omega_2$ systems, the systems of dif\/ferential operators that we study, are constructed from the~$\tau_2$ map.}

By construction the $\Omega_2$ systems arise from irreducible constituents $W$ of ${\mathfrak l} \otimes
{\mathfrak z}({\mathfrak n})$ so that $\tilde{\tau}_2|_{W^*}$ is not identically zero
(see the procedures described in Subsection~\ref{SS:Omega}). Since such irreducible constituents play
a~role to determine the special values of the $\Omega_2$ systems, for the remainder of this section, we
recall from~\cite{KuboThesis1} the observation on ${\mathfrak l}$ and irreducible constituents $W$ for
which $\tilde{\tau}_2|_{W^*} \not \equiv0$.

We start with the structure of ${\mathfrak l} = {\mathfrak z}({\mathfrak l}) \oplus [{\mathfrak
l},{\mathfrak l}]$, where ${\mathfrak z}({\mathfrak l})$ is the center of ${\mathfrak l}$.
First, as ${\mathfrak q} = {\mathfrak l} \oplus {\mathfrak n}$ is the maximal parabolic subalgebra
determined by $\alpha_{{\mathfrak q}}$, we have ${\mathfrak z}({\mathfrak l}) = \mathbb{C} H_{\mathfrak q}$.
The semisimple part $[{\mathfrak l},{\mathfrak l}]$ is either simple or the direct sum of two or three
simple ideals
(see for example Appendix~A of~\cite{KuboThesis1}). To characterize the simple ideals, let $\gamma$ be the
highest root for~${\mathfrak g}$.
If ${\mathfrak g}$ is not of type~$A_n$ then there is exactly one simple root not orthogonal to $\gamma$.
It is well known that if $\alpha_\gamma$ is the unique simple root then the nilpotent radical ${\mathfrak
n}'$ of the parabolic subalgebra ${\mathfrak q}' = {\mathfrak q}_{\{\alpha_\gamma\}}$ satisf\/ies
$\dim([{\mathfrak n}', {\mathfrak n}']) = 1$.
Recall from Subsection 2.2 that we say ${\mathfrak q}={\mathfrak l} \oplus {\mathfrak n}$ is of
quasi-Heisenberg type if $\dim([{\mathfrak n}, {\mathfrak n}])>1$.
Hence, if ${\mathfrak q} = {\mathfrak q}_{\{\alpha_{\mathfrak q}\}}$ is a~parabolic subalgebra of
quasi-Heisenberg type then~$\alpha_\gamma$ is in $\Pi({\mathfrak l}) = \Pi \backslash \{\alpha_{\mathfrak
q}\}$.
In particular there is a~unique simple ideal of $[{\mathfrak l},{\mathfrak l}]$ containing the root space
${\mathfrak g}_{\alpha_\gamma}$ for~$\alpha_\gamma$.
We denote by $\frak{l}_{\gamma}$ the unique simple ideal containing ${\mathfrak g}_{\alpha_\gamma}$.
Similarly, when $[{\mathfrak l},{\mathfrak l}]$ consists of two (resp.\
three) simple ideals, we denote the other simple ideal(s) by $\frak{l}_{n\gamma}$ (resp.\
$\frak{l}_{n\gamma}^+$ and $\frak{l}_{n\gamma}^-$).
The three simple factors occur only when ${\mathfrak q}$ is of type~$D_n(n-2)$.
So, when ${\mathfrak q}$ is not of type~$D_n(n-2)$, the Levi subalgebra ${\mathfrak l}$ may decompose into
\begin{gather}
\label{Eqn4.1.8}
{\mathfrak l}=\mathbb{C}H_{\mathfrak q}\oplus\frak{l}_{\gamma}\oplus\frak{l}_{n\gamma}.
\end{gather}
Similarly, when ${\mathfrak q}$ is of type~$D_n(n-2)$, one may write
\begin{gather}
\label{Eqn:LD(n-2)}
{\mathfrak l}=\mathbb{C}H_{\mathfrak q}\oplus\frak{l}_{\gamma}\oplus\frak{l}_{n\gamma}^+\oplus\frak{l}
_{n\gamma}^-.
\end{gather}
Note that when $[{\mathfrak l}, {\mathfrak l}]$ is a~simple ideal, we have $\frak{l}_{n\gamma} = \{0\}$
($\frak{l}_{n\gamma}^\pm = \{0\}$).
It follows from the decompositions~\eqref{Eqn4.1.8} and~\eqref{Eqn:LD(n-2)} that the tensor product
${\mathfrak l} \otimes {\mathfrak z}({\mathfrak n})$ may be written as
\begin{gather*}
{\mathfrak l}\otimes{\mathfrak z}({\mathfrak n})=
\begin{cases}
\big(\mathbb{C}H_{\mathfrak q}\otimes{\mathfrak z}({\mathfrak n})\big)\oplus\big(\frak{l}_{\gamma}
\otimes{\mathfrak z}({\mathfrak n})\big)\oplus\big(\frak{l}_{n\gamma}\otimes{\mathfrak z}({\mathfrak n})\big),
\quad
{\mathfrak q}~\text{is not of type}~D_n(n-2),
\\
\big(\mathbb{C}H_{\mathfrak q}\otimes{\mathfrak z}({\mathfrak n})\big)\oplus\big(\frak{l}_{\gamma}
\otimes{\mathfrak z}({\mathfrak n})\big)\oplus\big(\frak{l}_{n\gamma}^+\otimes{\mathfrak z}({\mathfrak n}
)\big)\oplus\big(\frak{l}_{n\gamma}^-\otimes{\mathfrak z}({\mathfrak n})\big),
\\
\hspace*{75mm}
{\mathfrak q}~\text{is of type}~D_n(n-2).
\end{cases}
\end{gather*}

To build an $\Omega_2$ system, it is necessary to choose an irreducible constituent $W$ in ${\mathfrak l}
\otimes {\mathfrak z}({\mathfrak n})$ so that the $L$-intertwining map
\begin{gather*}
\tilde{\tau}_2|_{W^*}\in\Hom_{L}\big(W^*,\mathcal{P}^2({\mathfrak g}(1))\big)
\end{gather*}
is not identically zero.
Now we give a~necessary condition for irreducible constituents $W$ so that $\tilde{\tau}_2|_{W^*} \not
\equiv 0$.
To do so, for $\nu \in {\mathfrak h}^*$ with $\IP{\nu}{\alpha^\vee} \in \mathbb{Z}_{\geq 0}$ for all
$\Pi({\mathfrak l})$, we denote by $V(\nu)$ the simple ${\mathfrak l}$-module with highest weight
$\nu|_{{\mathfrak h} \cap [{\mathfrak l},{\mathfrak l}]}$.

Suppose that ${\mathfrak l}\otimes {\mathfrak z}({\mathfrak n})$ has irreducible constituent $V(\nu)$.
If the linear map $\tilde{\tau}_2\big|_{V(\nu)^*}: V(\nu)^* \to \mathcal{P}^2({\mathfrak g}(1))$ is not
identically zero then, via the isomorphism $\mathcal{P}^2({\mathfrak g}(1)) \cong \Sym^2({\mathfrak
g}(1))^*$, $V(\nu)$ should be an irreducible constituent of $\Sym^2({\mathfrak g}(1)) \subset
{\mathfrak g}(1) \otimes {\mathfrak g}(1)$.
In particular if $\mu$ is the highest weight of ${\mathfrak g}(1)$ then $\nu$ is of the form $\nu=\mu+\ge$
for some $\ge \in \Delta({\mathfrak g}(1))$.
It was shown in Lemma 4.14 of~\cite{KuboThesis1} that the highest root $\gamma$ is of this form.
However, it follows from Proposition 6.5 of~\cite{KuboThesis1} that~$V(\gamma)$ does not occur in
$\Sym^2({\mathfrak g}(1))$.
Based on this observation we give the following def\/inition.
\begin{Def}[\protect{\cite[Def\/inition 6.7]{KuboThesis1}}]
\label{Def:Special}
An irreducible constituent $V(\nu)$ of ${\mathfrak l} \otimes {\mathfrak z}({\mathfrak n})$ is called
\textit{special}\footnote{There is a~certain discrepancy on the terminology ``special constituents''.
See the comments in the introduction on this matter.} if $V(\nu)$ satisf\/ies the following two conditions:{\samepage
\begin{enumerate}\itemsep=0pt
\item[(C1)] $\nu = \mu +\ge$ for some $\ge \in \Delta({\mathfrak g}(1))$.
\item[(C2)] $\nu \neq \gamma$.
\end{enumerate}}
\end{Def}

It is shown in~\cite[Section 6]{KuboThesis1} that, for ${\mathfrak q}$ not of type~$D_n(n-2)$,
there are exactly one or two special constituents of ${\mathfrak l} \otimes {\mathfrak z}({\mathfrak n})$;
one is an irreducible constituent of $\frak{l}_{\gamma}\otimes {\mathfrak z}({\mathfrak n})$ and the other
is equal to $\frak{l}_{n\gamma} \otimes {\mathfrak z}({\mathfrak n})$.
We denote by $V(\mu+\epsilon_\gamma)$ and $V(\mu+\epsilon_{n\gamma})$ the special constituents so that
$V(\mu+\epsilon_\gamma) \subset \frak{l}_{\gamma}\otimes {\mathfrak z}({\mathfrak n})$ and
$V(\mu+\epsilon_{n\gamma}) = \frak{l}_{n\gamma}\otimes {\mathfrak z}({\mathfrak n})$.
It will be shown in Section~\ref{SS:D(n-2)} that if ${\mathfrak q}$ is of type~$D_n(n-2)$ then there
are three special constituents, namely, $V(\mu+\epsilon_\gamma) \subset \frak{l}_{\gamma}\otimes {\mathfrak
z}({\mathfrak n})$ and $V(\mu+\epsilon_{n\gamma}^\pm) = \frak{l}_{n\gamma}^\pm\otimes {\mathfrak
z}({\mathfrak n})$.

To compute the special values of $\Omega_2$ systems ef\/f\/iciently, the special constituents $V(\mu+\ge)$
are classif\/ied as Type~1a, Type~1b, Type~2, or Type~3 as follows:
\begin{Def}[\protect{\cite[Def\/inition 6.20]{KuboThesis1}}]\label{Def:Type}
Let $\mu$ be the highest weight for ${\mathfrak g}(1)$.
We say that a~special constituent $V(\mu+\ge)$ is of
\begin{enumerate}\itemsep=0pt
\item [(1)] \textit{Type~1a} if $\mu + \ge$ is not a~root with $\ge \neq \mu$ and both $\mu$ and $\ge$ are
long roots, \item [(2)] \textit{Type~1b} if $\mu + \ge$ is not a~root with $\ge \neq \mu$ and either $\mu$
or $\ge$ is a~short root, \item [(3)] \textit{Type~2} if $\mu + \ge = 2\mu$ is not a~root, or \item [(4)]
\textit{Type~3} if $\mu + \ge$ is a~root.
\end{enumerate}
\end{Def}

Table~\ref{table-1} in the introduction shows the types of special constituents.
A dash in the column for $V(\mu + \epsilon_{n\gamma})$ indicates that $\frak{l}_{n\gamma} = \{0\}$ for the
case.
(So there is no special constituent $V(\mu+\epsilon_{n\gamma})$.)

In Sections~\ref{SS:Type1b} and~\ref{SS:Type3}, we f\/ind the special values of the $\Omega_2$
systems coming from special constituents of Type~1b and Type~3, respectively.
The following proposition will play a~key role to determine the special values.
\begin{Prop}
\label{Prop:HL}
Let $V(\mu+\ge)^*$ be the dual module of an irreducible constituent $V(\mu+\ge)$ of ${\mathfrak l}\otimes
{\mathfrak z}({\mathfrak n})$ so that the operator $\Omega_2|_{V(\mu+\ge)^*}: V(\mu+\ge)^* \to
\mathbb{D}(\mathcal{L}_{s})^{\bar {\mathfrak n}}$ is non-zero.
If $X_h$ is a~highest weight vector for ${\mathfrak g}(1)$ and if $Y_l^*$ is a~lowest weight vector for
$V(\mu+\ge)^*$ then the $\Omega_2|_{V(\mu+\ge)^*}$ system is conformally invariant on $\mathcal{L}_{s}$ if
and only if in $M_{\mathfrak q}[\mathbb{C}_{-s}]$
\begin{gather}
\label{Eqn:HL}
X_\mu\cdot\big(\omega_2(Y^*_l)\otimes1_{-s}\big)=0.
\end{gather}
\end{Prop}
\begin{proof}
Observe that, by Proposition~\ref{Prop:CIS-GVM}, the $\Omega_2|_{V(\mu+\ge)^*}$ system is conformally
invariant if and only if, for all $X \in {\mathfrak n}$ and $Y^* \in V(\mu+\ge)^*$,
\begin{gather}
\label{Eqn:Prop2.22}
X\cdot\big(\omega_2(Y^*)\otimes1_{-s}\big)=X\omega_2(Y^*)\otimes1_{-s}=0.
\end{gather}
Therefore, to prove this proposition, it suf\/f\/ices to show that~\eqref{Eqn:HL}
implies~\eqref{Eqn:Prop2.22}.
Since the arguments are similar to ones for Proposition 7.13 with Lemma 3.9 and Lemma 3.12
of~\cite{KuboThesis1}, we omit the proof.
\end{proof}

\section[Parabolic subalgebra of type~$D_n(n-2)$]{Parabolic subalgebra of type~$\boldsymbol{D_n(n-2)}$}
\label{SS:D(n-2)}

In this short section we consider the $\Omega_2$ systems associated with maximal parabolic subalgebra
${\mathfrak q} = {\mathfrak l} \oplus {\mathfrak n}$ of type~$D_n(n-2)$.
In particular we f\/ind the special values for the $\Omega_2$ systems and determine the standardness for
the maps $\varphi_{\Omega_2}$.

The parabolic subalgebra ${\mathfrak q}$ of type~$D_n(n-2)$ is the maximal parabolic subalgebra determined
by the simple root $\alpha_{n-2}$; the deleted Dynkin diagram is
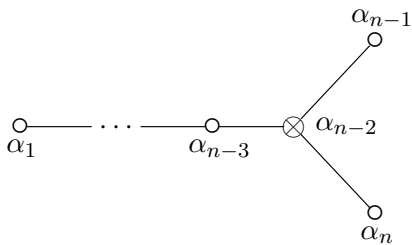
\begin{gather*}
\xymatrix{&&&&\abovewnode{\alpha_{n-1}}\\
\belowwnode{\alpha_1}\ar@{-}[r]&\cdots\ar@{-}[r]&
\belowwnode{\alpha_{n-3}}\ar@{-}[r]&
*-<\nodesize>{\crossandcirclesymbol}\save[]+<20pt,0pt>*\txt{$\alpha_{n-2}$} \restore\ar@{-}[ur]\ar@{-}[dr]& \\
&&&&\belowwnode{\alpha_n} }
\end{gather*}
with subgraphs
\begin{gather*}
\xymatrix{
\belowwnode{\alpha_1}\ar@{-}[r]&
\belowwnode{\alpha_2}\ar@{-}[r]&
\belowwnode{\alpha_3}\ar@{-}[r]&\cdots\ar@{-}[r]
&\belowwnode{\alpha_{n-3}}
&&
\belowwnode{\alpha_{n-1} }
&&
\belowwnode{\alpha_{n}}}
\end{gather*}
As the simple root $\alpha_{\mathfrak q}$ that determines the parabolic subalgebra
${\mathfrak q}$ is $\alpha_{\mathfrak q} = \alpha_{n-2}$, the fundamental weight $\lambda_{\mathfrak q}$
for $\alpha_{{\mathfrak q}}$ is $\lambda_{\mathfrak q} = \lambda_{n-2}$.
(For the def\/inition of deleted Dynkin diagrams see, for instance, Section~\ref{SS:Type1bA}
of~\cite{KuboThesis1}.)

Recall from Section~\ref{SS:Omega_2} that the Levi subalgebra ${\mathfrak l}$ may be decomposed as
\begin{gather*}%\label{Eqn:D(n-2)}
{\mathfrak l}=\mathbb{C}H_{\mathfrak q}\oplus\frak{l}_{\gamma}\oplus\frak{l}_{n\gamma}^+\oplus\frak{l}
_{n\gamma}^-.
\end{gather*}
The unique simple root $\alpha_\gamma$ that is not orthogonal to the highest root $\gamma$ is
$\alpha_\gamma = \alpha_2$.
Therefore we have $\frak{l}_{\gamma}\cong \mathfrak{sl}(n-2,\mathbb{C})$ and $\frak{l}_{n\gamma}^{\pm}
\cong \mathfrak{sl}(2,\mathbb{C})$.
For convenience we set $\frak{l}_{n\gamma}^+$ (resp.\
$\frak{l}_{n\gamma}^-$) to be the simple ideal that corresponds to the singleton for $\alpha_n$ (resp.\
$\alpha_{n-1}$).

\subsection{Special constituents and special values}
%\label{SS:D(n-2)A}

We start with f\/inding the special constituents of ${\mathfrak l} \otimes {\mathfrak z}({\mathfrak n})$.
As shown in Section~\ref{SS:Omega_2}, the tensor product ${\mathfrak l} \otimes {\mathfrak z}({\mathfrak n})$ may
be decomposed into
\begin{gather*}
{\mathfrak l}\otimes{\mathfrak z}({\mathfrak n})=\big(\mathbb{C}H_{\mathfrak q}\otimes{\mathfrak z}
({\mathfrak n})\big)\oplus\big(\frak{l}_{\gamma}\otimes{\mathfrak z}({\mathfrak n})\big)\oplus\big(\frak{l}
_{n\gamma}^+\otimes{\mathfrak z}({\mathfrak n})\big)\oplus\big(\frak{l}_{n\gamma}^-\otimes{\mathfrak z}
({\mathfrak n})\big).
\end{gather*}
With the arguments in~\cite[Section 6.1]{KuboThesis1}, one can easily check that
$\mathbb{C} H_{\mathfrak q}\otimes {\mathfrak z}({\mathfrak n})$ and $\frak{l}_{n\gamma}^\pm\otimes
{\mathfrak z}({\mathfrak n})$ are simple ${\mathfrak l}$-modules.
In fact, if we use the standard realization of roots with $\alpha_{n-1} = \varepsilon_{n-1} -
\varepsilon_{n}$ and $\alpha_n = \varepsilon_{n-1} + \varepsilon_{n}$ then
\begin{gather*}
\mathbb{C}H_{\mathfrak q}\otimes{\mathfrak z}({\mathfrak n}
)=V(\gamma)=V(\varepsilon_1+\varepsilon_2),
\\
\frak{l}_{n\gamma}^+\otimes{\mathfrak z}({\mathfrak n}
)=V(\alpha_n+\gamma)=V(\varepsilon_1+\varepsilon_2+\varepsilon_{n-1}+\varepsilon_{n}),\qquad \text{and}
\\
\frak{l}_{n\gamma}^-\otimes{\mathfrak z}({\mathfrak n})=V(\alpha_{n-1}
+\gamma)=V(\varepsilon_1+\varepsilon_2+\varepsilon_{n-1}-\varepsilon_{n}).
\end{gather*}
On the other hand, the tensor product $\frak{l}_{\gamma} \otimes {\mathfrak
z}({\mathfrak n})$ is reducible.
By using the character formula of Klimyk \cite[Corollary]{Klimyk68}, it can be shown that
\begin{gather*}
\frak{l}_{n\gamma}\otimes{\mathfrak z}({\mathfrak n})=V(2\varepsilon_1+\varepsilon_2-\varepsilon_{n-2}
)\oplus V(\varepsilon_1+\varepsilon_2)\oplus V(2\varepsilon_1).
\end{gather*}

Now one may observe that only $V(2\varepsilon_1)$ and $V(\varepsilon_1 + \varepsilon_2 + \varepsilon_{n-1}
\pm \varepsilon_{n})$ satisfy the condi\-tions~(C1) and~(C2) in Def\/inition~\ref{Def:Special}.
Thus these irreducible constituents are the special constituents of ${\mathfrak l} \otimes {\mathfrak
z}({\mathfrak n})$.
Write $\epsilon_\gamma$ and $\epsilon_{n\gamma}^\pm$ for the roots in $\Delta({\mathfrak g}(1))$ so that
$\mu+\epsilon_\gamma = 2\varepsilon_1$ and $\mu+\epsilon_{n\gamma}^\pm = \varepsilon_1 + \varepsilon_2 +
\varepsilon_{n-1} \pm \varepsilon_{n}$, where $\mu$ is the highest weight for ${\mathfrak g}(1)$.
Tables~\ref{T:Weights} and~\ref{T:Special} summarize the data for the special constituents.
\begin{table}[h] \caption{Roots $\mu$, $\epsilon_\gamma$, $\epsilon_{n\gamma}^+$, and $\epsilon_{n\gamma}^-$.}\label{T:Weights}\vspace{1mm}
\centering
\begin{tabular}{cccccccccc}
\hline
Parabolic ${\mathfrak q}$ &&$\mu$ &&$\epsilon_\gamma$ &&$\epsilon_{n\gamma}^+$&&$\epsilon_{n\gamma}^-$
\\
\hline
$D_n(n-2)$ &&$\varepsilon_1 + \varepsilon_{n-1}$ &&$\varepsilon_1 - \varepsilon_{n-1}$ &&$\varepsilon_2 +
\varepsilon_{n}$ &&$\varepsilon_2 - \varepsilon_{n}$
\\
\hline
\end{tabular}
\end{table}

\begin{table}[h] \caption{Highest weights for special constituents.}\label{T:Special}\vspace{1mm}
\centering
\begin{tabular}{ccccccccc}
\hline
Parabolic ${\mathfrak q}$ &&$V(\mu + \epsilon_\gamma)$ &&$V(\mu + \epsilon_{n\gamma}^+)$ &&$V(\mu +
\epsilon_{n\gamma}^-)$
\\
\hline
$D_n(n-2)$ &&$2\varepsilon_1$ &&$\varepsilon_1 +\varepsilon_2 + \varepsilon_{n-1} + \varepsilon_{n}$
&&$\varepsilon_1 +\varepsilon_2 + \varepsilon_{n-1} - \varepsilon_{n}$
\\
\hline
\end{tabular}
\end{table}

Observe that $\mu$, $\epsilon_\gamma$, and $\epsilon_{n\gamma}^\pm$ are all long roots and that neither
$\mu+\epsilon_\gamma$ nor $\mu+\epsilon_{n\gamma}^\pm$ is a~root.
Thus the special constituents $V(\mu+\epsilon_\gamma)$ and $V(\mu+\epsilon_{n\gamma}^\pm)$ are all of Type~1a
(see Def\/inition~\ref{Def:Type}).
As ${\mathfrak l} \otimes {\mathfrak z}({\mathfrak n})$ contains
a~special constituent of Type~1a, it follows from the argument in the proof for Proposition~7.3
of~\cite{KuboThesis1} that the $\tau_2$ map is not identically zero.
Also, the argument for Proposition~7.5 of~\cite{KuboThesis1} shows that, for $V(\mu+\ge) =
V(\mu+\epsilon_\gamma), V(\mu+\epsilon_{n\gamma}^\pm)$, the linear map $\tilde{\tau}_2|_{V(\mu+\ge)^*}$ is
not identically zero.

Now we are going to f\/ind the special values of the $\Omega_2$ systems arising from the special
constituents $V(\mu+\epsilon_\gamma)$ and $V(\mu+\epsilon_{n\gamma}^\pm)$.
If $W \subset {\mathfrak g}$ is an $\ad({\mathfrak h})$-invariant subspace then, for any weight $\nu
\in {\mathfrak h}^*$, we write
\begin{gather}
\label{Eqn:Delta}
\Delta_{\nu}(W):=\{\alpha\in\Delta(W)\,|\,\nu-\alpha\in\Delta\}.
\end{gather}
\begin{Thm}
\label{Thm:DnSpecial}
Let ${\mathfrak q}$ be the maximal parabolic subalgebra of type~$D_n(n-2)$.
If $V(\mu+\ge) = V(\mu+\epsilon_\gamma)$ or $V(\mu+\epsilon_{n\gamma}^\pm)$ then the
$\Omega_2|_{V(\mu+\ge)^*}$ system is conformally invariant on $\mathcal{L}_s$ if and only~if
\begin{gather*}
s=\frac{|\Delta_{\mu+\ge}({\mathfrak g}(1))|}{2}-1,
\end{gather*}
where $|\Delta_{\mu+\ge}({\mathfrak g}(1))|$ is the number of the elements in $\Delta_{\mu+\ge}({\mathfrak
g}(1))$.
\end{Thm}

\begin{proof}
Since the special constituents $V(\mu+\ge)$ and $V(\mu+\epsilon_\gamma^\pm)$ are of Type~1a, all the
statements in~\cite{KuboThesis1} for Type~1a special constituents hold for $V(\mu+\ge)$ and
$V(\mu+\epsilon_\gamma^\pm)$.
Now the theorem follows from the arguments in~\cite[Section 7]{KuboThesis1}.
\end{proof}
\begin{Cor}
\label{Cor:Special}
Under the same hypotheses in Theorem~{\rm \ref{Thm:DnSpecial}}, all $\Omega_2|_{V(\mu+\ge)^*}$ systems are
conformally invariant on $\mathcal{L}_1$.
\end{Cor}

\begin{proof}
By inspection we have $|\Delta_{\mu+\ge}({\mathfrak g}(1))| =4$ for each special constituent.
Now the results follow from Theorem~\ref{Thm:DnSpecial}.
\end{proof}

\subsection[The standardness of the map $\varphi_{\Omega_2}$]{The standardness of the map $\boldsymbol{\varphi_{\Omega_2}}$}
%\label{SS:D(n-2)C}

In the rest of this section we determine whether or not the maps $\varphi_{\Omega_2}$ coming from the
$\Omega_2$ systems are standard.

Observe that
\begin{gather*}
V(\mu+\epsilon_\gamma)^*=V(2\varepsilon_1)^*=V(-2\varepsilon_{n-2})
\end{gather*}
and
\begin{gather*}
V\big(\mu+\epsilon_{n\gamma}^\pm\big)^*=V(\varepsilon_1+\varepsilon_2+\varepsilon_{n-1}\pm\varepsilon_{n}
)^*=V(-\varepsilon_{n-3}-\varepsilon_{n-2}+\varepsilon_{n-1}\pm\varepsilon_{n}).
\end{gather*}
It follows from Corollary~\ref{Cor:Special} that the special value $s_2$ is $s_2 = 1$ for each
special constituent.
Therefore if we denote by $\varphi_{(\Omega_2, \mu+\epsilon_\gamma)}$ (resp.\
$\varphi_{(\Omega_2, \mu+\epsilon_{n\gamma}^\pm)}$) the homomorphism induced by the
$\Omega_2|_{V(\mu+\epsilon_\gamma)^*}$ system (reps.
the $\Omega_2|_{V(\mu+\epsilon_{n\gamma}^\pm)^*}$ system) then, by~\eqref{Eqn:GVM2}, we have
\begin{gather}
\label{Eqn:3.3GVM1}
\varphi_{(\Omega_2,\mu+\epsilon_\gamma)}:M_{\mathfrak q}(-2\varepsilon_{n-2}-\lambda_{n-2}
+\rho)\to M_{\mathfrak q}(-\lambda_{n-2}+\rho)
\end{gather}
and
\begin{gather}
\label{Eqn:3.3GVM2}
\varphi_{\big(\Omega_2,\mu+\epsilon_{n\gamma}^\pm\big)}:M_{\mathfrak q}((-\varepsilon_{n-3}
-\varepsilon_{n-2}+\varepsilon_{n-1}\pm\varepsilon_{n})-\lambda_{n-2}+\rho)\to M_{\mathfrak q}
(-\lambda_{n-2}+\rho).
\end{gather}

\begin{Thm}
\label{Thm:Map_D(n-2)}
If ${\mathfrak q}$ is the maximal parabolic subalgebra of type~$D_n(n-2)$ then the standard maps between
the generalized Verma modules in~\eqref{Eqn:3.3GVM1} and~\eqref{Eqn:3.3GVM2} are zero.
Consequently, the maps $\varphi_{(\Omega_2,\mu+\epsilon_\gamma)}$ and
$\varphi_{(\Omega_2,\mu+\epsilon_{n\gamma}^\pm)}$ are non-standard.
\end{Thm}
\begin{proof}
To prove this theorem, by Proposition~\ref{Prop:Std}, it suf\/f\/ices to show that there exits
$\alpha_{\nu} \in \Pi({\mathfrak l})$ so that $-\alpha_\nu - \lambda_{n-2} +\rho$ is linked to
$\nu-\lambda_{n-2} + \rho$ for $\nu = -2\varepsilon_{n-2}$, $-\varepsilon_{n-3} -\varepsilon_{n-2} +
\varepsilon_{n-1} \pm \varepsilon_{n}$.
Since the arguments are the same for each case, we only demonstrate a~proof for~\eqref{Eqn:3.3GVM1}.
To show for~\eqref{Eqn:3.3GVM1}, observe that
\begin{gather*}
-2\varepsilon_{n-2}=-2(\varepsilon_{n-2}-\varepsilon_{n-1})-(\varepsilon_{n-1}
-\varepsilon_n)-(\varepsilon_{n-1}+\varepsilon_n)
\end{gather*}
with $\varepsilon_{n-2}-\varepsilon_{n-1} \in \Delta({\mathfrak g}(1))$ and $\varepsilon_{n-1} \pm
\varepsilon_n \in \Pi({\mathfrak l})$
(see Appendix~\ref{SS:Data}).
We claim that $(\varepsilon_{n-1}-\varepsilon_n,
\varepsilon_{n-2}-\varepsilon_{n-1})$ links $-(\varepsilon_{n-1}+\varepsilon_n) - \lambda_{n-2} + \rho$ to
$-2\varepsilon_{n-2}-\lambda_{n-2} + \rho$, that is,
\begin{gather*}
s_{\varepsilon_{n-2}-\varepsilon_{n-1}}s_{\varepsilon_{n-1}-\varepsilon_n}(-(\varepsilon_{n-1}
+\varepsilon_n)-\lambda_{n-2}+\rho)=-2\varepsilon_{n-2}-\lambda_{n-2}+\rho
\end{gather*}
with
\begin{gather*}
\IP{-(\varepsilon_{i+1}-\varepsilon_{i+2})-\lambda_{n-2}+\rho}{(\varepsilon_{n-1}-\varepsilon_n)^\vee}
\in\mathbb{Z}_{\geq0}
\end{gather*}
and
\begin{gather*}
\IP{s_{\varepsilon_{n-1}-\varepsilon_n}(-(\varepsilon_{n-1}+\varepsilon_n)-\lambda_{n-2}+\rho)}
{(\varepsilon_{n-2}-\varepsilon_{n-1})^\vee}\in\mathbb{Z}_{\geq0}.
\end{gather*}
A direct computation shows that it indeed holds.
For~\eqref{Eqn:3.3GVM2}, one may use $(\varepsilon_{n-3} - \varepsilon_{n-2},
\varepsilon_{n-2}-\varepsilon_{n-1})$ as a~link.
Now the proposed statements follow.
\end{proof}

\section{Type 1b special constituent}
\label{SS:Type1b}

In this section we study the $\Omega_2$ system associated to the Type~1b special constituent.
It follows from Table~\ref{table-1} in the introduction that Type~1b special constituent occurs only
when the parabolic subalgebra ${\mathfrak q}$ is of type~$B_n(n-1)$.
The simple root $\alpha_{\mathfrak q}$ that determines the parabolic subalgebra~${\mathfrak q}$ is then
$\alpha_{\mathfrak q} = \alpha_{n-1}$.
We write $\lambda_{{\mathfrak q}} = \lambda_{n-1}$ for the fundamental weight $\lambda_{{\mathfrak q}}$ for
$\alpha_{{\mathfrak q}}$.
The deleted Dynkin diagram for ${\mathfrak q}$ is
\begin{gather*}
\xymatrix{\belowwnode{\alpha_1}\ar@{-}[r]&
\belowwnode{\alpha_2}\ar@{-}[r]&\cdots\ar@{-}[r]
&\belowwnode{\alpha_{n-2}}\ar@{-}[r]
&\belowcnode{\alpha_{n-1}}\ar@2{->}[r]&\belowwnode{\alpha_n}}
\end{gather*}
with connected subgraphs
\begin{gather*}
\xymatrix{\belowwnode{\alpha_1}\ar@{-}[r]&
\belowwnode{\alpha_2}\ar@{-}[r]&
\belowwnode{\alpha_3}\ar@{-}[r]&\cdots\ar@{-}[r]
&\belowwnode{\alpha_{n-2}}
&
&
\belowwnode{\alpha_n}}
\end{gather*}
Since $\alpha_2$ is the unique simple root that is not orthogonal to the highest root $\gamma$,
it follows from the subgraphs that $\frak{l}_{\gamma} \cong \mathfrak{sl}(n-1,\mathbb{C})$ and
$\frak{l}_{n\gamma} \cong \mathfrak{sl}(2,\mathbb{C})$.
Recall from Table~\ref{table-1} that  Type~1b special constituent of ${\mathfrak l} \otimes {\mathfrak
z}({\mathfrak n})$ is the irreducible constituent $V(\mu+\epsilon_{n\gamma})=\frak{l}_{n\gamma}\otimes
{\mathfrak z}({\mathfrak n})$.

\subsection[The $\tilde{\tau}_2|_{V(\mu+\epsilon_{n\gamma})^*}$ map]{The
$\boldsymbol{\tilde{\tau}_2|_{V(\mu+\epsilon_{n\gamma})^*}}$ map}
\label{SS:Type1bA}

We start with observing the $L$-intertwining operator $\tilde{\tau}_2|_{V(\mu+\epsilon_{n\gamma})^*}:
V(\mu+\epsilon_{n\gamma})^* \to \mathcal{P}^2({\mathfrak g}(1))$.
To do so we f\/irst f\/ix convenient root vectors for ${\mathfrak g}$ so that certain computations will be
carried out easily.
Observe that, as ${\mathfrak q}$ is of type~$B_n(n-1)$, the Lie algebra ${\mathfrak g}$ under consideration
is ${\mathfrak g} = \mathfrak{so}(2n+1,\mathbb{C})$.
We take ${\mathfrak h}$ to be the set of block diagonal matrices
\begin{gather*}
H(h_1,\ldots,h_n)=\diag\bigg(
\begin{pmatrix}
0&\mathbf{i}h_1
\\
-\mathbf{i}h_1&0
\end{pmatrix}
,
\begin{pmatrix}
0&\mathbf{i}h_2
\\
-\mathbf{i}h_2&0
\end{pmatrix}
,\ldots,
\begin{pmatrix}
0&\mathbf{i}h_n
\\
-\mathbf{i}h_n&0
\end{pmatrix}
,0\bigg)
\end{gather*}
with $h_j \in \mathbb{C}$ and $\mathbf{i} = \sqrt{-1}$.
The positive roots are $\Delta^+ = \{ \varepsilon_j \pm \varepsilon_k \, | \, 1 \leq j < k \leq n\} \cup \{
\varepsilon_j \, | \, 1 \leq j \leq n\}$ with $\varepsilon_j(H(h_1, \ldots, h_n)) = h_j$.
We take the root vectors $X_\alpha$ as follows:
\begin{enumerate}\itemsep=0pt
\item $\alpha = \pm (\varepsilon_j \pm \varepsilon_k)$ $(j <k)$:
\begin{gather*}
\phantom{X_\alpha=}
\qquad
j
\qquad
k
\\
X_\alpha=
\begin{pmatrix}
0&E_\alpha
\\
-E_\alpha^{t}&0
\end{pmatrix}
\quad
\begin{matrix}
j
\\
k
\end{matrix}
\end{gather*}
with
\begin{gather*}
E_{\varepsilon_j-\varepsilon_k}=\frac{1}{2}
\begin{pmatrix}
1&\mathbf{i}
\\
-\mathbf{i}&1
\end{pmatrix}
,
\qquad
E_{\varepsilon_j+\varepsilon_k}=\frac{1}{2}
\begin{pmatrix}
1&-\mathbf{i}
\\
-\mathbf{i}&-1
\end{pmatrix},
\\
E_{-\varepsilon_j+\varepsilon_k}=\frac{1}{2}
\begin{pmatrix}
-1&\mathbf{i}
\\
-\mathbf{i}&-1
\end{pmatrix},
\qquad
E_{-\varepsilon_j-\varepsilon_k}=\frac{1}{2}
\begin{pmatrix}
-1&-\mathbf{i}
\\
-\mathbf{i}&1
\end{pmatrix}
,
\end{gather*}
where $X_\alpha$ denotes the matrix whose entries are all zero except the $j$th and
$k$th pairs of indices.
\item
$\alpha = \pm \varepsilon_j$:
\begin{gather*}
\phantom{X_\alpha=(}
\quad
j
\quad
2n+1
\\
X_\alpha=
\begin{pmatrix}
0&v_\alpha
\\
-v_\alpha^{t}&0
\end{pmatrix}
\quad
\begin{matrix}
j
\\
2n+1
\end{matrix}
\end{gather*}
with
\begin{gather*}
v_{\varepsilon_j}=\frac{1}{\sqrt{2}}
\begin{pmatrix}
1
\\
-\mathbf{i}
\end{pmatrix}
,
\qquad
v_{-\varepsilon_j}=\frac{1}{\sqrt{2}}
\begin{pmatrix}
-1
\\
-\mathbf{i}
\end{pmatrix}
,
\end{gather*}
where $X_\alpha$ denotes the matrix def\/ined similarly to the previous case.
\end{enumerate}

For $\alpha, \beta \in \Delta$ with $\alpha + \beta \neq 0$, let $N_{\alpha, \beta}$ denote the constant so
that $[X_{\alpha}, X_{\beta}] = N_{\alpha,\beta} X_{\alpha+\beta}$.
Table~\ref{T:Constants1} summarizes the values of $N_{\alpha, \beta}$ for $\alpha + \beta$ a~positive root.
The constant $N_{\alpha, \beta}$ satisf\/ies the property that $N_{-\alpha, -\beta} = -N_{\alpha, \beta}$
(see, for instance, \cite[Theorem~6.6]{Knapp02}).
Thus the values of $N_{\alpha, \beta}$ for $\alpha +
\beta$ a~negative root can also be obtained from Table~\ref{T:Constants1}.
Here, we would like to acknowledge that, for the cases that $\alpha$ and $\beta$ are long roots (formulas
(1)--(12)), we simply adapt the beautiful multiplication table by Knapp in~\cite[Section 10]{Knapp04}.
\begin{table}[t] \caption{The values of $N_{\alpha,\beta}$ for roots $\alpha$ and $\beta$ for
$\mathfrak{so}(2n+1,\mathbb{C})$ with indices $i<j<k$ when $\alpha+\beta$ is a~positive root.}\label{T:Constants1}\vspace{1mm}
\centering
\begin{tabular}{cccccccccc}
\hline
Formula &&&$\alpha$ &&&$\beta$&&&$N_{\alpha,\beta}$
\\
\hline
(1) &&& $\varepsilon_i + \varepsilon_k$ &&& $\varepsilon_j - \varepsilon_k$ &&& $-1$
\\
(2) &&& $\varepsilon_i - \varepsilon_k$ &&& $\varepsilon_j + \varepsilon_k$ &&& $-1$
\\
(3) &&& $\varepsilon_i + \varepsilon_k$ &&& $-\varepsilon_j - \varepsilon_k$ &&& $+1$
\\
(4) &&& $\varepsilon_i - \varepsilon_k$ &&& $-\varepsilon_j + \varepsilon_k$ &&& $+1$
\\
(5) &&& $\varepsilon_i + \varepsilon_j$ &&& $-\varepsilon_j + \varepsilon_k$ &&& $-1$
\\
(6) &&& $\varepsilon_i - \varepsilon_j$ &&& $\varepsilon_j + \varepsilon_k$ &&& $+1$
\\
(7) &&& $\varepsilon_i + \varepsilon_j$ &&& $-\varepsilon_j - \varepsilon_k$ &&& $-1$
\\
(8) &&& $\varepsilon_i - \varepsilon_j$ &&& $\varepsilon_j - \varepsilon_k$ &&& $+1$
\\
(9) &&& $\varepsilon_i + \varepsilon_j$ &&& $-\varepsilon_i + \varepsilon_k$ &&& $+1$
\\
(10) &&& $-\varepsilon_i + \varepsilon_j$ &&& $\varepsilon_i + \varepsilon_k$ &&& $+1$
\\
(11) &&& $\varepsilon_i + \varepsilon_j$ &&& $-\varepsilon_i - \varepsilon_k$ &&& $+1$
\\
(12) &&& $-\varepsilon_i + \varepsilon_j$ &&& $\varepsilon_i - \varepsilon_k$ &&& $+1$
\\
(13) &&& $\varepsilon_j$ &&& $-\varepsilon_j + \varepsilon_k$ &&& $-1$
\\
(14) &&& $-\varepsilon_j$ &&& $\varepsilon_j+\varepsilon_k$ &&& $-1$
\\
(15) &&& $\varepsilon_k$ &&& $\varepsilon_j -\varepsilon_k$ &&& $-1$
\\
(16) &&& $-\varepsilon_k$ &&& $\varepsilon_j + \varepsilon_k$ &&& $+1$
\\
(17) &&& $\varepsilon_j$ &&& $\varepsilon_k$ &&& $-1$
\\
(18) &&& $\varepsilon_j$ &&& $-\varepsilon_k$ &&& $+1$
\\
\hline
\end{tabular}
\end{table}

For $\alpha \in \Delta^+$ we set $H_\alpha = [X_\alpha,X_{-\alpha}]$; namely, we have
\begin{gather*}
H_{\varepsilon_j\pm\varepsilon_k}=\diag\bigg(0,\ldots,0,
\stackrel{j}{\begin{pmatrix}
0&\mathbf{i}
\\
-\mathbf{i}&0
\end{pmatrix}}
,0,\ldots,0,\pm
\stackrel{k}{\begin{pmatrix}
0&\mathbf{i}
\\
-\mathbf{i}&0
\end{pmatrix}}
,0,\ldots,0\bigg)
\end{gather*}
and
\begin{gather*}
H_{\varepsilon_j}=\diag\bigg(0,\ldots,0,
\stackrel{j}{\begin{pmatrix}
0&\mathbf{i}
\\
-\mathbf{i}&0
\end{pmatrix}}
,0,\ldots,0,\bigg).
\end{gather*}

Now observe if $T(X, Y):= \Tr(XY)$ then $T(X_\alpha, X_{-\alpha}) =2 $ for all $\alpha \in \Delta^+$.
As the restriction $T(\cdot, \cdot)|_{{\mathfrak h}_0 \times {\mathfrak h}_0}$ to a~real form ${\mathfrak
h}_0$ of ${\mathfrak h}$ is an inner product on ${\mathfrak h}_0$, the trace form $T(\cdot, \cdot)$ is
a~positive constant multiple of the Killing form $\kappa(\cdot, \cdot)$.
If $b_0$ is the non-zero constant so that $\kappa(X,Y) = b_0 T(X,Y)$ then $\kappa(X_\alpha, X_{-\alpha}) =
2b_0$ for all $\alpha \in \Delta^+$.
Thus the dual vectors $X^*_\alpha$ for $X_\alpha$ with respect to the Killing form are $X^*_\alpha =
(1/(2b_0))X_{-\alpha}$.

Now the element $\omega$ in~\eqref{Eqn:omega} is given by
\begin{gather*}
\omega=\sum_{\gamma_j\in\Delta({\mathfrak z}({\mathfrak n}))}X^*_{\gamma_j}\otimes X_{\gamma_j}=\frac{1}
{2b_0}\sum_{\gamma_j\in\Delta({\mathfrak z}({\mathfrak n}))}X_{-\gamma_j}\otimes X_{\gamma_j}.
\end{gather*}
Thus the map $\tau_2(X)$ in~\eqref{Eqn:Tau2} may be expressed as
\begin{gather}
\label{Eqn:Tau2_2}
\tau_2(X)=\frac{1}{2}\big(\ad(X)^2\otimes\Id\big)\omega=\frac{1}{4b_0}
\sum_{\gamma_j\in\Delta({\mathfrak z}({\mathfrak n}))}\ad(X)^2X_{-\gamma_j}\otimes X_{\gamma_j}.
\end{gather}
Proposition 7.3 of~\cite{KuboThesis1} showed that when ${\mathfrak q}$ is of type~$B_n(n-1)$, the $\tau_2$
map is not identically zero.

To construct dif\/ferential operators from $V(\mu+\epsilon_{n\gamma})^*$, it is necessary to show that the
linear map $\tilde{\tau}_2|_{V(\mu+\epsilon_{n\gamma})^*}:V(\mu+\epsilon_{n\gamma})^* \to
\mathcal{P}^2({\mathfrak g}(1))$ is not identically zero
(see \emph{Step 2} and \emph{Step 3} in Section~\ref{SS:Omega}).
We shall prove this by showing that
$\tilde{\tau}_2(Y^*)(X)$ is a~non-zero polynomial on ${\mathfrak g}(1)$ for some $Y^*$ in
$V(\mu+\epsilon_{n\gamma})^*$.
To this end observe that, as discussed in Section~\ref{SS:Omega_2}, we have
$V(\mu+\epsilon_{n\gamma})=\frak{l}_{n\gamma} \otimes {\mathfrak z}({\mathfrak n})$.
Therefore $V(\mu+\epsilon_{n\gamma})^*=\frak{l}_{n\gamma}^*\otimes{\mathfrak z}({\mathfrak n})^* =
\frak{l}_{n\gamma}\otimes {\mathfrak z}(\bar{\mathfrak n})$.
Since $\gamma \in \Delta({\mathfrak z}({\mathfrak n}))$ and $\frak{l}_{n\gamma} =
\Span_\mathbb{C}\{X_{\alpha_n}, H_{\alpha_n}, X_{-\alpha_n}\}$, we have $X^*_{\alpha_n} \in
\frak{l}_{n\gamma}$ and $X^*_\gamma \in {\mathfrak z}(\bar{\mathfrak n})$; therefore $X^*_{\alpha_n}
\otimes X^*_\gamma \in \frak{l}_{n\gamma} \otimes {\mathfrak z}(\bar{\mathfrak n})
=V(\mu+\epsilon_{n\gamma})^*$.
\begin{Prop}
\label{Prop:Tau2Type1b}
If ${\mathfrak q} = {\mathfrak l} \oplus {\mathfrak n}$ is the parabolic subalgebra of type~$B_n(n-1)$ then
the $L$-intertwining operator $\tilde{\tau}_2|_{V(\mu+\epsilon_{n\gamma})^*}$ is not identically zero.
\end{Prop}
\begin{proof}
We show that $\tilde{\tau}_2(X^*_{\alpha_n} \otimes X^*_\gamma)(X)$ is a~non-zero polynomial on ${\mathfrak
g}(1)$.
As $X^*_{\alpha}$ denotes the dual vector for $X_\alpha$ with respect to the Killing form $\kappa$,
by~\eqref{Eqn:Tau2_2}, the operator $\tilde{\tau}_2(X^*_{\alpha_n} \otimes X^*_{\gamma})(X)
=(X^*_{\alpha_n} \otimes X^*_{\gamma})(\tau_2(X))$ is given by
\begin{gather}
\tilde{\tau}_2\big(X^*_{\alpha_n}\otimes X^*_{\gamma}\big)(X)=(X^*_{\alpha_n}\otimes X^*_{\gamma})(\tau_2(X))
\nonumber
\\
\hphantom{\tilde{\tau}_2\big(X^*_{\alpha_n}\otimes X^*_{\gamma}\big)(X)}{} =\frac{1}{4b_0}\sum_{\gamma_j\in\Delta({\mathfrak z}({\mathfrak n}))}\kappa\big(X^*_{\alpha_n},\ad
(X)^2X_{-\gamma_j}\big)\kappa(X^*_{\gamma},X_{\gamma_j})
\nonumber
\\
\hphantom{\tilde{\tau}_2\big(X^*_{\alpha_n}\otimes X^*_{\gamma}\big)(X)}{} =\frac{1}{4b_0}\kappa\big(X^*_{\alpha_n},\ad(X)^2X_{-\gamma}\big)
=\frac{1}{8b_0^2}\kappa\big(X_{-\alpha_n},\ad(X)^2X_{-\gamma}\big).
\label{Eqn:Kappa_Bn0}
\end{gather}
Write $X = \sum\limits_{\alpha \in \Delta({\mathfrak g}(1))} \eta_\alpha X_{\alpha}$, where $\eta_\alpha
\in {\mathfrak n}^*$ is the coordinate dual to $X_\alpha$ with respect to the Killing form.
Recall from~\eqref{Eqn:Delta} that if $W \subset {\mathfrak g}$ is an $\ad({\mathfrak h})$-invariant
subspace then, for any weight $\nu \in {\mathfrak h}^*$, we write $\Delta_{\nu}(W) = \{ \alpha \in
\Delta(W) \, | \, \nu - \alpha \in \Delta\}$.
Then,
\begin{gather}
\eqref{Eqn:Kappa_Bn0}=\frac{1}{8b_0^2}\kappa\big(X_{-\alpha_n},\ad(X)^2X_{-\gamma}\big)
\nonumber
\\
\phantom{\eqref{Eqn:Kappa_Bn0}}
=\frac{1}{8b_0^2}\sum_{\beta,\delta\in\Delta({\mathfrak g}(1))}
\eta_\beta\eta_\delta\kappa(X_{-\alpha_n},[X_\delta,[X_\beta,X_{-\gamma}]])
\nonumber
\\
\phantom{\eqref{Eqn:Kappa_Bn0}}
=\frac{1}{8b_0^2}\sum_{\alpha,\beta\in\Delta({\mathfrak g}(1))}
\eta_\beta\eta_\delta\kappa([X_{-\alpha_n},X_\delta],[X_\beta,X_{-\gamma}]])
\nonumber
\\
\phantom{\eqref{Eqn:Kappa_Bn0}}
=\frac{1}{8b_0^2}\sum_{\substack{\beta\in\Delta_\gamma({\mathfrak g}(1))
\\
\delta\in\Delta_{\alpha_n}({\mathfrak g}(1))}}\eta_\beta\eta_\delta N_{\beta,-\gamma}N_{-\alpha_n,\delta}
\kappa(X_{\delta-\alpha_n},X_{\beta-\gamma}).
\label{Eqn:Kappa_Bn1}
\end{gather}
One may observe that $\kappa(X_{\delta-\alpha_n},X_{\beta-\gamma}) \neq 0$ if and
only if $\delta = \alpha_n +\gamma - \beta \in \Delta({\mathfrak g}(1))$.
If we write
\begin{gather*}%\label{Eqn:ThetaType1b}
\theta(\beta)=\alpha_n+\gamma-\beta
\end{gather*}
then
\begin{gather}
\eqref{Eqn:Kappa_Bn1}=
\frac{1}{8b_0^2}\sum_{\substack{\beta\in\Delta_\gamma({\mathfrak g}(1))\\\delta\in\Delta_{\alpha_n}({\mathfrak g}(1))}}
\eta_\beta\eta_\delta N_{\beta,-\gamma}N_{-\alpha_n,\delta}
\kappa(X_{\delta-\alpha_n},X_{\beta-\gamma})
\nonumber
\\
\phantom{\eqref{Eqn:Kappa_Bn1}}
=\frac{1}{8b_0^2}\sum_{\beta\in\Delta_\gamma({\mathfrak g}(1))\cap\Delta_{\alpha_n+\gamma}
({\mathfrak g}(1))}\eta_\beta\eta_{\theta(\beta)}N_{\beta,-\gamma}N_{-\alpha_n,\theta(\beta)}
\kappa(X_{\theta(\beta)-\alpha_n},X_{\beta-\gamma})
\nonumber
\\
\phantom{\eqref{Eqn:Kappa_Bn1}}
=\frac{1}{8b_0^2}\sum_{\beta\in\Delta_\gamma({\mathfrak g}(1))\cap\Delta_{\alpha_n+\gamma}
({\mathfrak g}(1))}\eta_\beta\eta_{\theta(\beta)}N_{\beta,-\gamma}N_{-\alpha_n,\theta(\beta)}
\kappa(X_{\gamma-\beta},X_{\beta-\gamma})
\nonumber
\\
\phantom{\eqref{Eqn:Kappa_Bn1}}
=\frac{1}{4b_0}\sum_{\beta\in\Delta_\gamma({\mathfrak g}(1))\cap\Delta_{\alpha_n+\gamma}
({\mathfrak g}(1))}\eta_\beta\eta_{\theta(\beta)}N_{\beta,-\gamma}N_{-\alpha_n,\theta(\beta)}
\kappa(X^*_{\beta-\gamma},X_{\beta-\gamma})
\nonumber
\\
\phantom{\eqref{Eqn:Kappa_Bn1}}
=\frac{1}{4b_0}\sum_{\beta\in\Delta_\gamma({\mathfrak g}(1))\cap\Delta_{\alpha_n+\gamma}
({\mathfrak g}(1))}\eta_\beta\eta_{\theta(\beta)}N_{\beta,-\gamma}N_{-\alpha_n,\theta(\beta)}
\nonumber
\\
\phantom{\eqref{Eqn:Kappa_Bn1}}
=\frac{1}{4b_0}\sum_{\beta\in\Delta_\gamma({\mathfrak g}(1))\cap\Delta_{\alpha_n+\gamma}
({\mathfrak g}(1))}N_{\beta,-\gamma}N_{-\alpha_n,\theta(\beta)}
\kappa(X,X^*_\beta)\kappa(X,X^*_{\theta(\beta)}).
\label{Eqn:Kappa_Bn2}
\end{gather}
As $\alpha_n = \varepsilon_n$ and $\gamma = \varepsilon_1 + \varepsilon_2$, by
inspection, we have
\begin{gather*}
\Delta_\gamma({\mathfrak g}(1))=\{\varepsilon_1\pm\varepsilon_n,\varepsilon_2\pm\varepsilon_n\}
\cup\{\varepsilon_1,\varepsilon_2\}
\end{gather*}
and
\begin{gather*}
\Delta_{\alpha_n+\gamma}({\mathfrak g}(1))=\{\varepsilon_1+\varepsilon_n,\varepsilon_2+\varepsilon_n\}
\cup\{\varepsilon_1,\varepsilon_2\}.
\end{gather*}
In particular, $\Delta_{\alpha_n+\gamma}({\mathfrak g}(1)) \subset
\Delta_\gamma({\mathfrak g}(1))$.
Therefore,
\begin{gather*}
\eqref{Eqn:Kappa_Bn2}=\frac{1}{4b_0}\sum_{\beta\in\Delta_\gamma({\mathfrak g}
(1))\cap\Delta_{\alpha_n+\gamma}({\mathfrak g}(1))}N_{\beta,-\gamma}N_{-\alpha_n,\theta(\beta)}
\kappa(X,X^*_\beta)\kappa(X,X^*_{\theta(\beta)})
\\
\phantom{\eqref{Eqn:Kappa_Bn2}}
=\frac{1}{4b_0}\sum_{\beta\in\Delta_{\alpha_n+\gamma}({\mathfrak g}(1))}N_{\beta,-\gamma}
N_{-\alpha_n,\theta(\beta)}\kappa(X,X^*_\beta)\kappa(X,X^*_{\theta(\beta)}).
\end{gather*}
Therefore we obtain
\begin{gather*}
\tilde{\tau}_2(X^*_{\alpha_n}\otimes X^*_{\gamma})(X)=\frac{1}{4b_0}\sum_{\beta\in\Delta_{\alpha_n+\gamma}
({\mathfrak g}(1))}N_{\beta,-\gamma}N_{-\alpha_n,\theta(\beta)}
\kappa(X,X^*_\beta)\kappa(X,X^*_{\theta(\beta)}).
\end{gather*}
Since $\Delta_{\alpha_n+\gamma}({\mathfrak g}(1)) =\{\varepsilon_1 + \varepsilon_n,
\varepsilon_2 + \varepsilon_n\} \cup \{\varepsilon_1, \varepsilon_2\}$, this reads
\begin{gather}
\tilde{\tau}_2\big(X^*_{\alpha_n}\otimes X^*_{\gamma}\big)(X)
=\frac{1}{4b_0}\sum_{\beta\in\Delta_{\alpha_n+\gamma}({\mathfrak g}(1))}N_{\beta,-\gamma}
N_{-\alpha_n,\theta(\beta)}\kappa(X,X^*_\beta)\kappa(X,X^*_{\theta(\beta)})\nonumber
\\
\hphantom{\tilde{\tau}_2\big(X^*_{\alpha_n}\otimes X^*_{\gamma}\big)(X)}{}
 =\frac{1}{4b_0}\big(N_{\varepsilon_1+\varepsilon_n,-(\varepsilon_1+\varepsilon_2)}
N_{-\varepsilon_n,\varepsilon_2}\kappa(X,X^*_{\varepsilon_1+\varepsilon_n})\kappa(X,X^*_{\varepsilon_2})
\nonumber
\\
\hphantom{\tilde{\tau}_2\big(X^*_{\alpha_n}\otimes X^*_{\gamma}\big)(X)=}{}
+N_{\varepsilon_2,-(\varepsilon_1+\varepsilon_2)}N_{-\varepsilon_n,\varepsilon_1+\varepsilon_n}
\kappa(X,X^*_{\varepsilon_2})\kappa(X,X^*_{\varepsilon_1+\varepsilon_n})
\nonumber
\\
\hphantom{\tilde{\tau}_2\big(X^*_{\alpha_n}\otimes X^*_{\gamma}\big)(X)=}{}
+N_{\varepsilon_2+\varepsilon_n,-(\varepsilon_1+\varepsilon_2)}N_{-\varepsilon_n,\varepsilon_1}
\kappa(X,X^*_{\varepsilon_2+\varepsilon_n})\kappa(X,X^*_{\varepsilon_1})
\nonumber
\\
\hphantom{\tilde{\tau}_2\big(X^*_{\alpha_n}\otimes X^*_{\gamma}\big)(X)=}{}
+N_{\varepsilon_1,-(\varepsilon_1+\varepsilon_2)}N_{-\varepsilon_n,\varepsilon_2+\varepsilon_n}
\kappa(X,X^*_{\varepsilon_1})\kappa(X,X^*_{\varepsilon_2+\varepsilon_n})\big)
\nonumber
\\
\hphantom{\tilde{\tau}_2\big(X^*_{\alpha_n}\otimes X^*_{\gamma}\big)(X)}{} =\frac{1}{4b_0}\big((1)(-1)\kappa(X,X^*_{\varepsilon_1+\varepsilon_n})\kappa(X,X^*_{\varepsilon_2}
)+(-1)(1)\kappa(X,X^*_{\varepsilon_2})\kappa(X,X^*_{\varepsilon_1+\varepsilon_n})
\nonumber
\\
\hphantom{\tilde{\tau}_2\big(X^*_{\alpha_n}\otimes X^*_{\gamma}\big)(X)=}{}
+(-1)(-1)\kappa(X,X^*_{\varepsilon_2+\varepsilon_n})\kappa(X,X^*_{\varepsilon_1}
)(1)(1)\kappa(X,X^*_{\varepsilon_1})\kappa(X,X^*_{\varepsilon_2+\varepsilon_n})\big)
\nonumber
\\
\hphantom{\tilde{\tau}_2\big(X^*_{\alpha_n}\otimes X^*_{\gamma}\big)(X)}{} =\frac{1}{2b_0}\big(\kappa(X,X^*_{\varepsilon_2+\varepsilon_n})\kappa(X,X^*_{\varepsilon_1}
)-\kappa(X,X^*_{\varepsilon_1+\varepsilon_n})\kappa(X,X^*_{\varepsilon_2})\big).
\label{Eqn:LowestType1b}
\end{gather}
Therefore $\tilde{\tau}_2(X^*_{\alpha_n} \otimes X^*_{\gamma})(X)$ is a~non-zero
polynomial on ${\mathfrak g}(1)$.
\end{proof}

\subsection[The special value of the $\Omega_2|_{V(\mu+\epsilon_{n\gamma})^*}$ system]{The special value of
the $\boldsymbol{\Omega_2|_{V(\mu+\epsilon_{n\gamma})^*}}$ system}
%\label{SS:Type1bB}

Now we f\/ind the special value of the $\Omega_2|_{V(\mu+\epsilon_{n\gamma})^*}$ system.
To do so we use Proposition~\ref{Prop:HL}.
To this end recall from Section~\ref{SS:GVM-I} the linear map $\omega_2|_{V(\mu+\epsilon_{n\gamma})^*}:
V(\mu+\epsilon_{n\gamma})^* \to \mathcal{U}(\bar {\mathfrak n})$ def\/ined by
$\omega_2|_{V(\mu+\epsilon_{n\gamma})^*} = \sigma \circ \tilde{\tau}_2|_{V(\mu+\epsilon_{n\gamma})^*}$,
where $\sigma: \Sym^2({\mathfrak g}(-1)) \to \mathcal{U}(\bar {\mathfrak n})$ the symmetrization
operator
(here we identify $\mathcal{P}^2({\mathfrak g}(1)) \cong \Sym^2({\mathfrak g}(-1))$).
If
\begin{gather*}
Y_l^*:=8b_0^3(X^*_{\alpha_n}\otimes X^*_\gamma)
\end{gather*}
then, by~\eqref{Eqn:LowestType1b}, $\omega_2(Y_l^*):=
\omega_2|_{V(\mu+\epsilon_{n\gamma})^*}(Y_l^*)$ is given by
\begin{gather*}
\omega_2(Y_l^*)=4b_0^2\big(\sigma(X^*_{\varepsilon_2+\varepsilon_n}X^*_{\varepsilon_1}
)-\sigma(X^*_{\varepsilon_1+\varepsilon_n}X^*_{\varepsilon_2})\big).
\end{gather*}
As the dual vector $X^*_\alpha$ for $X_\alpha$ with respect to the Killing form is
$X^*_\alpha = (1/2b_0)X_{-\alpha}$, this amounts~to
\begin{gather*}
\omega_2\big(Y_l^*\big)=\sigma(X_{-(\varepsilon_2+\varepsilon_n)}X_{-\varepsilon_1}
)-\sigma(X_{-(\varepsilon_1+\varepsilon_n)}X_{-\varepsilon_2}).
\end{gather*}
Moreover, since $\varepsilon_1+\varepsilon_2 +\varepsilon_n \not \in \Delta$, the
symmetrization is unnecessary.
Therefore we obtain
\begin{gather}
\label{Eqn:Omega2Type1b}
\omega_2\big(Y_l^*\big)=X_{-(\varepsilon_2+\varepsilon_n)}X_{-\varepsilon_1}
-X_{-(\varepsilon_1+\varepsilon_n)}X_{-\varepsilon_2}.
\end{gather}

Now we are going to determine the special value of the $\Omega_2|_{V(\mu+\epsilon_{n\gamma})^*}$ system.
\begin{Thm}
\label{Thm:SpValType1b}
Let ${\mathfrak q}$ be the maximal parabolic subalgebra of type~$B_n(n-1)$.
The $\Omega_2|_{V(\mu+\epsilon_{n\gamma})^*}$ system is conformally invariant on $\mathcal{L}_s$ if and
only if $s = 1$.
\end{Thm}
\begin{proof}
Observe that, as $X_{\alpha_n} \otimes X_\gamma$ is a~highest weight vector for $\frak{l}_{n\gamma}\otimes
{\mathfrak z}({\mathfrak n}) = V(\mu+\epsilon_{n\gamma})$, $X^*_{\alpha_n} \otimes X^*_\gamma$ is a~lowest
weight vector for $V(\mu+\epsilon_{n\gamma})^*$; consequently, $Y_l^*$ is a~lowest weight vector for
$V(\mu+\epsilon_{n\gamma})^*$.
Therefore, by Proposition~\ref{Prop:HL}, to f\/ind the special value for the
$\Omega_2|_{V(\mu+\epsilon_{n\gamma})^*}$ system, it suf\/f\/ices to determine $s \in \mathbb{C}$ so that
$X_{\mu} \cdot (\omega_2(Y_l^*)\otimes 1_{-s}) = 0$ for $\mu$ the highest weight for ${\mathfrak g}(1)$.
By inspection, we have $\mu = \varepsilon_1 + \varepsilon_n$
(see Appendix~\ref{SS:Data}).
Thus we compute $X_{\varepsilon_1 + \varepsilon_n} \cdot
(\omega_2(Y_l^*)\otimes 1_{-s})$.

It follows from~\eqref{Eqn:Omega2Type1b} that
\begin{gather*}
X_{\varepsilon_1+\varepsilon_n}\cdot(\omega_2\big(Y_l^*\big)\otimes1_{-s})=X_{\varepsilon_1+\varepsilon_n}
X_{-(\varepsilon_2+\varepsilon_n)}X_{-\varepsilon_1}\otimes1_{-s}-X_{\varepsilon_1+\varepsilon_n}
X_{-(\varepsilon_1+\varepsilon_n)}X_{-\varepsilon_2}\otimes1_{-s}.
\end{gather*}
We observe from the second term.
By the standard computation, we have
\begin{gather*}
T_2=-X_{\varepsilon_1+\varepsilon_n}X_{-(\varepsilon_1+\varepsilon_n)}X_{-\varepsilon_2}\otimes1_{-s}
=-H_{\varepsilon_1+\varepsilon_n}X_{-\varepsilon_2}\otimes1_{-s}=s\lambda_{n-1}
(H_{\varepsilon_1+\varepsilon_n})X_{-\varepsilon_2}\otimes1_{-s}.
\end{gather*}
Observe that since $\lambda_{n-1} = \sum\limits_{j=1}^{n-1}\varepsilon_j$, we have
$\lambda_{n-1}(H_{\varepsilon_1 + \varepsilon_n}) = 1$.
Thus,
\begin{gather*}
T_2=s\lambda_{n-1}(H_{\varepsilon_1+\varepsilon_n})X_{-\varepsilon_2}\otimes1_{-s}=sX_{-\varepsilon_2}
\otimes1_{-s}.
\end{gather*}
Similarly the standard computation shows that the f\/irst term amounts to
\begin{gather*}
T_1=X_{\varepsilon_1+\varepsilon_n}X_{-(\varepsilon_2+\varepsilon_n)}X_{-\varepsilon_1}\otimes1_{-s}
=-X_{-\varepsilon_2}\otimes1_{-s}.
\end{gather*}
Therefore,
\begin{gather*}
X_{\varepsilon_1+\varepsilon_n}\cdot(\omega_2\big(Y_l^*\big)\otimes1_{-s})=T_1+T_2=(s-1)X_{-\varepsilon_2}
\otimes1_{-s}.
\end{gather*}
Now the assertion follows from Proposition~\ref{Prop:HL}.
\end{proof}
\begin{Rem}
\label{Rem:SpValType1b}
Theorem 7.16 in~\cite{KuboThesis1} shows that the special values $s_2$
for the $\Omega_2$ systems associated to Type~1a special constituents $V(\mu+\ge)$
are given by $s_2 = \big(|\Delta_{\mu+\ge}({\mathfrak g}(1))|/2 \big)-1$,
where $|\Delta_{\mu+\ge}({\mathfrak g}(1))|$ is the number
of elements in $\Delta_{\mu+\ge}({\mathfrak g}(1))$.
Now since $\Delta_{\mu+\epsilon_{n\gamma}}({\mathfrak g}(1))=\Delta_{\alpha_n+\gamma}({\mathfrak g}(1))
=\{\varepsilon_1 + \varepsilon_n,$ $\varepsilon_2 + \varepsilon_n\} \cup \{\varepsilon_1, \varepsilon_2\}$,
the special value $s_2$ may be expressed as
\begin{gather*}
s_2=1=\frac{|\Delta_{\mu+\epsilon_{n\gamma}}({\mathfrak g}(1))|}{2}-1.
\end{gather*}
Thus the special value $s_2$ for the $\Omega_2$ systems associated to Type~1a and
Type~1b special constituents can be given in the same formula.
\end{Rem}

\subsection[The standardness of the map $\varphi_{\Omega_2}$]{The standardness of the map
$\boldsymbol{\varphi_{\Omega_2}}$}
%\label{SS:Type1bC}
In the rest of this section we determine whether or not the map $\varphi_{\Omega_2}$ coming from the
$\Omega_2|_{V(\mu+\epsilon_{n\gamma})^*}$ system is standard.

Observe that, as $V(\mu+\epsilon_{n\gamma}) = \frak{l}_{n\gamma} \otimes {\mathfrak z}({\mathfrak n})$ and
as the highest weights for $\frak{l}_{n\gamma}$ and ${\mathfrak z}({\mathfrak n})$ are $\varepsilon_n$ and
$-\varepsilon_{n-2} - \varepsilon_{n-1}$, respectively, we have
\begin{gather*}
V(\mu+\epsilon_{n\gamma})^*=\frak{l}_{n\gamma}\otimes{\mathfrak z}(\bar{\mathfrak n})=V(-\varepsilon_{n-2}
-\varepsilon_{n-1}+\varepsilon_n).
\end{gather*}
By Theorem~\ref{Thm:SpValType1b}, the special value $s_2$ for the
$\Omega_2|_{V(\mu+\epsilon_{n\gamma})^*}$ system is $s_2 = 1$.
Therefore, by~\eqref{Eqn:GVM2}, the $\Omega_2|_{V(\mu+\epsilon_{n\gamma})^*}$ system yields a~non-zero
$\mathcal{U}({\mathfrak g})$-homomorphism
\begin{gather}
\label{Eqn:Type1bGVM}
\varphi_{\Omega_2}:M_{\mathfrak q}((-\varepsilon_{n-2}-\varepsilon_{n-1}+\varepsilon_n)-\lambda_{n-1}
+\rho)\to M_{\mathfrak q}(-\lambda_{n-1}+\rho).
\end{gather}
\begin{Thm}
\label{Thm:Map_Type1b}
If ${\mathfrak q}$ is the maximal parabolic subalgebra of type~$B_n(n-1)$ then the standard map
$\varphi_{std}$ between the generalized Verma modules in~\eqref{Eqn:Type1bGVM} is zero.
Consequently, the map $\varphi_{\Omega_2}$ is non-standard.
\end{Thm}

\begin{proof}
The idea of the proof is the same as for Theorem~\ref{Thm:Map_D(n-2)}.
Namely, f\/irst show that there exists $\alpha_{\nu} \in \Pi({\mathfrak l})$ so that $-\alpha_\nu -
\lambda_{n-1} +\rho$ is linked to $(-\varepsilon_{n-2} - \varepsilon_{n-1}+\varepsilon_n)-\lambda_{n-1} +
\rho$, then apply Proposition~\ref{Prop:Std}.
Observe that
\begin{gather*}
-\varepsilon_{n-2}-\varepsilon_{n-1}+\varepsilon_n=-2(\varepsilon_{n-1}-\varepsilon_n)-(\varepsilon_{n-2}
-\varepsilon_{n-1})-\varepsilon_n
\end{gather*}
with $\varepsilon_{n-1}-\varepsilon_{n} \in \Delta({\mathfrak g}(1))$ and $\varepsilon_{n-2} -
\varepsilon_{n-1}, \varepsilon_n \in \Pi({\mathfrak l})$
(see Appendix~\ref{SS:Data}).
By a~direct computation one can show that
$(\varepsilon_{n-2}-\varepsilon_{n-1}, \varepsilon_{n-1}-\varepsilon_{n})$ links $-\varepsilon_{n} -
\lambda_{n-1} + \rho$ to $(-\varepsilon_{n-2} - \varepsilon_{n-1}+\varepsilon_n)-\lambda_{n-1} + \rho$ Now
the theorem follows from Proposition~\ref{Prop:Std}.
\end{proof}

\section{Type 3 special constituent}
\label{SS:Type3}

In this section we study the $\Omega_2$ system associated to   Type~3 special constituent.
It follows from Table~\ref{table-1} in the introduction that   Type~3 special constituent occurs only when
the parabolic subalgebra ${\mathfrak q}$ is of type~$C_n(i)$ for $2 \leq i \leq n-1$.
The simple root $\alpha_{\mathfrak q}$ that determines the parabolic subalgebra ${\mathfrak q}$ is then
$\alpha_{\mathfrak q} = \alpha_i$.
We write $\lambda_{{\mathfrak q}} = \lambda_i$ for the fundamental weight $\lambda_{{\mathfrak q}}$ for
$\alpha_{{\mathfrak q}}$.
The deleted Dynkin diagram for ${\mathfrak q}$ is
\begin{gather*}
\xymatrix{\belowwnode{\alpha_1}\ar@{-}[r]&\cdots\ar@{-}[r]
&\belowwnode{\alpha_{i-1}}\ar@{-}[r]&\belowcnode{\alpha_i}\ar@{-}[r]&
\belowwnode{\alpha_{i+1}}\ar@{-}[r]&\cdots \ar@{-}[r]
&\belowwnode{\alpha_{n-1}}\ar@2{<-}[r]&\belowwnode{\alpha_n}}
\end{gather*}
with connected subgraphs
\begin{gather*}
\xymatrix{\belowwnode{\alpha_1}\ar@{-}[r]&
\belowwnode{\alpha_2}\ar@{-}[r]&
\belowwnode{\alpha_3}\ar@{-}[r]&\cdots\ar@{-}[r]
&\belowwnode{\alpha_{i-1}}
&
&
\belowwnode{\alpha_{i+1}}\ar@{-}[r]&\cdots \ar@{-}[r]
&\belowwnode{\alpha_{n-1}}\ar@2{<-}[r]&\belowwnode{\alpha_n}}
\end{gather*}
Since $\alpha_1$ is the unique simple root that is not orthogonal to the highest root $\gamma$,
it follows from the subgraphs that $\frak{l}_{\gamma} \cong \mathfrak{sl}(i,\mathbb{C})$ and
$\frak{l}_{n\gamma} \cong \mathfrak{sp}(n-i,\mathbb{C})$.
Recall from Table~\ref{table-1} that   Type~3 special constituent of ${\mathfrak l} \otimes {\mathfrak
z}({\mathfrak n})$ is the irreducible constituent $V(\mu+\epsilon_\gamma) \subset \frak{l}_{\gamma} \otimes
{\mathfrak z}({\mathfrak n})$.

As in Section~\ref{SS:Type1b}, our f\/irst goal is to show the $L$-intertwining operator
$\tilde{\tau}_2|_{V(\mu+\epsilon_\gamma)^*}$ is not identically zero.
To do so we again f\/ix convenient root vectors for ${\mathfrak g}$.
Observe that, as ${\mathfrak q}$ is of type~$C_n(i)$, the Lie algebra ${\mathfrak g}$ under consideration
is ${\mathfrak g} = \mathfrak{sp}(n,\mathbb{C})$.
For $1\leq j \leq n$, we write $\hat{j} = j+n$.
If $E_{ab}$ denotes the matrix with 1 in the $(a, b)$ entry and 0 elsewhere then we take ${\mathfrak h}$ to
be the set of diagonal matrices
\begin{gather*}
H(h_1,\ldots,h_n)=h_1(E_{11}-E_{\hat{1}\hat{1}})+\cdots+h_n(E_{nn}-E_{\hat{n}\hat{n}})
\end{gather*}
with $h_j \in \mathbb{C}$.
The positive roots are $\Delta^+ = \{\varepsilon_j \pm \varepsilon_k \, | \, 1\leq j < k \leq n\} \cup \{
2\varepsilon_j \, | \, 1\leq j \leq n\}$ with $\varepsilon_j(h_1, \ldots, h_n) = h_j$.
We take the root vectors $X_\alpha$ as follows:
\begin{gather*}
X_{\varepsilon_j-\varepsilon_k}=E_{jk}-E_{\hat{k}\hat{j}},
\qquad
X_{\varepsilon_j+\varepsilon_k}=E_{j\hat{k}}+E_{k\hat{j}},
\\
X_{-(\varepsilon_j+\varepsilon_k)}=E_{\hat{j}k}+E_{\hat{k}j},
\qquad
X_{2\varepsilon_j}=E_{j\hat{j}},
\qquad
X_{-2\varepsilon_j}=E_{\hat{j}j}.
\end{gather*}

For $\alpha, \beta \in \Delta$ with $\alpha + \beta \neq 0$, we again denote by $N_{\alpha, \beta}$ the
constant so that $[X_{\alpha}, X_{\beta}] = N_{\alpha,\beta} X_{\alpha+\beta}$.
Table~\ref{T:Constants2} summarizes the values of $N_{\alpha, \beta}$ for $\alpha + \beta$ a~positive root.
\begin{table}[t] \caption{The values of $N_{\alpha,\beta}$ for roots $\alpha$ and $\beta$ for
$\mathfrak{sp}(2n,\mathbb{C})$ with indices $i<j<k$ when $\alpha+\beta$ is a~positive root.}\label{T:Constants2}\vspace{1mm}
\centering
\begin{tabular}{cccccccccc}
\hline
Formula &&&$\alpha$ &&&$\beta$&&&$N_{\alpha,\beta}$
\\
\hline
(1) &&& $\varepsilon_i + \varepsilon_k$ &&& $\varepsilon_j - \varepsilon_k$ &&& $-1$
\\
(2) &&& $\varepsilon_i - \varepsilon_k$ &&& $\varepsilon_j + \varepsilon_k$ &&& $+1$
\\
(3) &&& $\varepsilon_i + \varepsilon_k$ &&& $-\varepsilon_j - \varepsilon_k$ &&& $+1$
\\
(4) &&& $\varepsilon_i - \varepsilon_k$ &&& $-\varepsilon_j + \varepsilon_k$ &&& $+1$
\\
(5) &&& $\varepsilon_i + \varepsilon_j$ &&& $-\varepsilon_j + \varepsilon_k$ &&& $-1$
\\
(6) &&& $\varepsilon_i - \varepsilon_j$ &&& $\varepsilon_j + \varepsilon_k$ &&& $+1$
\\
(7) &&& $\varepsilon_i + \varepsilon_j$ &&& $-\varepsilon_j - \varepsilon_k$ &&& $+1$
\\
(8) &&& $\varepsilon_i - \varepsilon_j$ &&& $\varepsilon_j - \varepsilon_k$ &&& $+1$
\\
(9) &&& $\varepsilon_i + \varepsilon_j$ &&& $-\varepsilon_i + \varepsilon_k$ &&& $-1$
\\
(10) &&& $-\varepsilon_i + \varepsilon_j$ &&& $\varepsilon_i + \varepsilon_k$ &&& $+1$
\\
(11) &&& $\varepsilon_i + \varepsilon_j$ &&& $-\varepsilon_i - \varepsilon_k$ &&& $+1$
\\
(12) &&& $-\varepsilon_i + \varepsilon_j$ &&& $\varepsilon_i - \varepsilon_k$ &&& $+1$
\\
(13) &&& $\varepsilon_i + \varepsilon_j$ &&& $\varepsilon_i - \varepsilon_j$ &&& $-2$
\\
(14) &&& $\varepsilon_i + \varepsilon_j$ &&& $-\varepsilon_i + \varepsilon_j$ &&& $-2$
\\
(15) &&& $2\varepsilon_j$ &&& $-\varepsilon_j + \varepsilon_k$ &&& $-1$
\\
(16) &&& $2\varepsilon_j$ &&& $-\varepsilon_j - \varepsilon_k$ &&& $+1$
\\
(17) &&& $2\varepsilon_k$ &&& $\varepsilon_j -\varepsilon_k$ &&& $-1$
\\
(18) &&& $-2\varepsilon_k$ &&& $\varepsilon_j + \varepsilon_k$ &&& $+1$
\\
\hline
\end{tabular}
\end{table}

For $\alpha \in \Delta^+$ we set $H_\alpha = [X_\alpha,X_{-\alpha}]$.
Namely,
\begin{gather*}
H_{\varepsilon_j\pm\varepsilon_k}=(E_{jj}-E_{\hat{j}\hat{j}})\pm(E_{kk}-E_{\hat{k}\hat{k}})
\qquad
\text{and}
\qquad
H_{2\varepsilon_j}=E_{jj}-E_{\hat{j}\hat{j}}.
\end{gather*}

Now observe that if $T(X, Y):= \Tr(XY)$ then, as $T(\cdot, \cdot)|_{{\mathfrak h}_0 \times {\mathfrak
h}_0}$ is an inner product on a~real form ${\mathfrak h}_0$ of ${\mathfrak h}$, there exists a~positive
constant $c_0$ so that $\kappa(X,Y) = c_0 T(X,Y)$.
Since $T(X_\alpha, X_{-\alpha})$ takes the value of one for $\alpha$ long and two for $\alpha$ short, we
have
\begin{gather}
\label{Eqn:Kappa}
\kappa(X_\alpha,X_{-\alpha})=
\begin{cases}
c_0&\text{if}~\alpha~\text{is long},
\\
2c_0&\text{if}~\alpha~\text{is short}.
\end{cases}
\end{gather}
Thus the dual vector $X^*_\alpha$ for $X_\alpha$ with respect to the Killing form is given by
\begin{gather}
\label{Eqn:DVect}
X^*_\alpha=
\begin{cases}
(1/c_0)X_{-\alpha}&\text{if}~\alpha~\text{is long},
\\
(1/(2c_0))X_{-\alpha}&\text{if}~\alpha~\text{is short}.
\end{cases}
\end{gather}
Now if $\Delta({\mathfrak z}({\mathfrak n}))_{\text{long}}$ (reps.\
$\Delta({\mathfrak z}({\mathfrak n}))_{\text{short}}$) is the set of long roots (reps.\
short roots) in $\Delta({\mathfrak z}({\mathfrak n}))$ then the element $\omega$ in~\eqref{Eqn:omega} is
given by
\begin{gather}
\omega=\sum_{\gamma_j\in\Delta({\mathfrak z}({\mathfrak n}))}X^*_{\gamma_j}\otimes X_{\gamma_j}
\nonumber
\\
\phantom{\omega}
=\frac{1}{c_0}\sum_{\gamma_j\in\Delta({\mathfrak z}({\mathfrak n}))_\text{long}}X_{-\gamma_j}
\otimes X_{\gamma_j}+\frac{1}{2c_0}\sum_{\gamma_j\in\Delta({\mathfrak z}({\mathfrak n}))_\text{short}}
X_{-\gamma_j}\otimes X_{\gamma_j}.
\label{Eqn:OmegaType3}
\end{gather}
Observe that $\Delta({\mathfrak z}({\mathfrak n})) = \{ \varepsilon_j + \varepsilon_k \, | \, 1\leq j<k\leq
i\} \cup \{2\varepsilon_j \, | \, 1 \leq j \leq i\}$.
Thus,
\begin{gather*}
\Delta({\mathfrak z}({\mathfrak n}))_\text{long}=\{2\varepsilon_j\,|\,1\leq j\leq i\}
\qquad
\text{and}
\qquad
\Delta({\mathfrak z}({\mathfrak n}))_\text{short}=\{\varepsilon_j+\varepsilon_k\,|\,1\leq j<k\leq i\}.
\end{gather*}
Therefore,~\eqref{Eqn:OmegaType3} reads
\begin{gather*}
\omega=\frac{1}{c_0}\sum_{j=1}^i X_{-2\varepsilon_j}\otimes X_{2\varepsilon_j}+\frac{1}{2c_0}
\sum_{1\leq j<k\leq i}X_{-(\varepsilon_j+\varepsilon_k)}\otimes X_{\varepsilon_j+\varepsilon_k}.
\end{gather*}
Thus the $\tau_2$ map in~\eqref{Eqn:Tau2} may be expressed as
\begin{gather*}%\label{Eqn:Tau2_3}
\tau_2(X)=\frac{1}{2}\big(\ad(X)^2\otimes\Id\big)\omega
\\
\phantom{\tau_2(X)}
=\frac{1}{2c_0}\sum_{j=1}^i\ad(X)^2X_{-2\varepsilon_j}\otimes X_{2\varepsilon_j}+\frac{1}{4c_0}
\sum_{1\leq j<k\leq i}\ad(X)^2X_{-(\varepsilon_j+\varepsilon_k)}
\otimes X_{\varepsilon_j+\varepsilon_k}.
\end{gather*}

As for the case that ${\mathfrak q}$ is of type~$B_n(n-1)$, Proposition~7.3 of~\cite{KuboThesis1} showed
that the $\tau_2$ map is not identically zero.

\subsection[Lowest weight vector for $V(\mu+\epsilon_\gamma)^*$]{Lowest weight vector for
$\boldsymbol{V(\mu+\epsilon_\gamma)^*}$}
%\label{SS:Type3A}

As for the case of Type~1b special constituent it is necessary to show that the linear map
$\tilde{\tau}_2|_{V(\mu+\epsilon_\gamma)^*}:V(\mu+\epsilon_\gamma)^* \to \mathcal{P}^2({\mathfrak g}(1))$
is not identically zero.
We will again achieve it by \mbox{showing} that $\tilde{\tau}_2(Y^*_l)(X)$ is a~non-zero polynomial on ${\mathfrak
g}(1)$, where $Y^*_l$ is a~lowest weight vector for \mbox{$V(\mu+\epsilon_\gamma)^*$}.

When a~special constituent is of Type~1b, as $V(\mu+\epsilon_{n\gamma})=\frak{l}_{n\gamma}\otimes
{\mathfrak z}({\mathfrak n})$, it was easy to f\/ind a~lowest weight vector for
$V(\mu+\epsilon_{n\gamma})^*$.
In contrast, in the present case, since $V(\mu+\epsilon_\gamma) \subsetneq \frak{l}_{\gamma} \otimes
{\mathfrak z}({\mathfrak n})$, we cannot use the same idea.
So our f\/irst goal is to f\/ind an explicit form of a~lowest weight vector for $V(\mu+\epsilon_\gamma)^*$.
To do so we now observe a~highest weight vector for $V(\mu+\epsilon_\gamma)$.

If $\pr_{\frak{l}_{\gamma} \otimes {\mathfrak z}({\mathfrak n})}: {\mathfrak l} \otimes {\mathfrak
z}({\mathfrak n}) \to \frak{l}_{\gamma} \otimes {\mathfrak z}({\mathfrak n})$ is the projection map from
${\mathfrak l}\otimes {\mathfrak z}({\mathfrak n})$ onto $\frak{l}_{\gamma} \otimes {\mathfrak
z}({\mathfrak n})$ then we claim that $\pr_{\frak{l}_{\gamma}\otimes{\mathfrak z}({\mathfrak
n})}(\tau_2(X_{\mu} + X_{\epsilon_\gamma}))$ is a~highest weight vector for $V(\mu+\epsilon_\gamma)$.
The following two technical lemmas will simplify the expositions of the proof.
\begin{Lem}
\label{Lem:4.4}
For $Z \in {\mathfrak l}$ and $X_1, X_2 \in {\mathfrak g}(1)$, we have
\begin{gather}
\label{Eqn:Technical}
Z\cdot\big(\ad(X_1)\ad(X_2)\otimes\Id
\big)\omega=\big((\ad([Z,X_1])\ad(X_2)+\ad(X_1)\ad
([Z,X_2]))\otimes\Id\big)\omega,
\end{gather}
where dot $(\cdot)$ denotes the usual Lie algebra action on tensor products.
\end{Lem}
\begin{proof}
As this lemma simply follows from the arguments used in the proof for Proposition 7.5
of~\cite{KuboThesis1}, we omit the proof.
\end{proof}
\begin{Lem}
\label{Lem:4.6}
We have
\begin{gather*}
\pr_{\frak{l}_{\gamma}\otimes{\mathfrak z}({\mathfrak n})}(\tau_2(X_{\mu}+X_{\epsilon_\gamma}))
\\
\qquad
=\frac{1}{2}\Big(\pr_{\frak{l}_{\gamma}\otimes{\mathfrak z}({\mathfrak n})}
(\ad(X_{\mu})\ad(X_{\epsilon_\gamma})\otimes\Id)\omega
+\pr_{\frak{l}_{\gamma}\otimes{\mathfrak z}({\mathfrak n})}
(\ad(X_{\epsilon_\gamma})\ad(X_{\mu})\otimes\Id)
\omega\Big).
\end{gather*}
\end{Lem}

\begin{proof}
Since
\begin{gather}
\tau_2(X_{\mu}+X_{\epsilon_\gamma})=\frac{1}{2}\big(\ad(X_{\mu}+X_{\epsilon_\gamma})^2\otimes\Id\big)\omega
\nonumber
\\
\phantom{\tau_2(X_{\mu}+X_{\epsilon_\gamma})}{}
 =\frac{1}{2}\Big((\ad(X_{\mu})^2\otimes\Id)\omega+(\ad(X_{\mu})\ad
(X_{\epsilon_\gamma})\otimes\Id)\omega
\nonumber
\\
\phantom{\tau_2(X_{\mu}+X_{\epsilon_\gamma})=}{}
+(\ad(X_{\epsilon_\gamma})\ad(X_{\mu})\otimes\Id)\omega+(\ad(X_{\epsilon_\gamma}
)^2\otimes\Id)\omega\Big),
\label{Eqn:4.71}
\end{gather}
we have
\begin{gather}
\pr_{\frak{l}_{\gamma}\otimes{\mathfrak z}({\mathfrak n})}
(\tau_2(X_{\mu}+X_{\epsilon_\gamma}))
\nonumber
\\
\qquad
=\frac{1}{2}
\Big(\pr_{\frak{l}_{\gamma}\otimes{\mathfrak z}({\mathfrak n})}(\ad(X_{\mu})^2\otimes\Id)\omega
+\pr_{\frak{l}_{\gamma}\otimes{\mathfrak z}({\mathfrak n})}
(\ad(X_{\mu})\ad(X_{\epsilon_\gamma})\otimes\Id)\omega
\nonumber
\\
\qquad\phantom{=}
+\pr_{\frak{l}_{\gamma}\otimes{\mathfrak z}({\mathfrak n})}
(\ad(X_{\epsilon_\gamma})\ad
(X_{\mu})\otimes\Id)
\omega
+\pr_{\frak{l}_{\gamma}\otimes{\mathfrak z}({\mathfrak n})}
(\ad(X_{\epsilon_\gamma})^2\otimes\Id)\omega \Big).
\label{Eqn:4.7}
\end{gather}
Observe that $\pr_{\frak{l}_{\gamma} \otimes {\mathfrak z}({\mathfrak n})}(\ad(X_{\mu})^2
\otimes \Id) \omega) = 0$.
Indeed, as $\mu = \varepsilon_1 + \varepsilon_{i+1}$ and $\omega =(1/c_0)\sum\limits_{j=1}^i
X_{-2\varepsilon_j} \otimes X_{2\varepsilon_j} +(1/2c_0)\sum\limits_{1\leq j < k \leq i} X_{-(\varepsilon_j
+ \varepsilon_k)} \otimes X_{\varepsilon_j+\varepsilon_k}$ (see~\eqref{Eqn:OmegaType3}), we have
\begin{gather*}
\big(\ad(X_{\mu})^2\otimes\Id\big)\omega
\\
\qquad =\frac{1}{c_0}\sum_{j=1}^i\ad(X_{\varepsilon_1+\varepsilon_{i+1}})^2X_{-2\varepsilon_j}
\otimes X_{2\varepsilon_j}+\frac{1}{2c_0}\sum_{1\leq j<k\leq i}\ad(X_{\varepsilon_1+\varepsilon_{i+1}}
)^2X_{-(\varepsilon_j+\varepsilon_k)}\otimes X_{\varepsilon_j+\varepsilon_k}
\\
\qquad =\frac{1}{c_0}\ad(X_{\varepsilon_1+\varepsilon_{i+1}})^2X_{-2\varepsilon_1}
\otimes X_{2\varepsilon_1}
\\
\qquad =\frac{1}{c_0}N_{\varepsilon_1+\varepsilon_{i+1},-2\varepsilon_1}N_{\varepsilon_1+\varepsilon_{i+1}
,-\varepsilon_1+\varepsilon_{i+1}}X_{2\varepsilon_{i+1}}\otimes X_{2\varepsilon_1}.
\end{gather*}
Clearly, $X_{2\varepsilon_{i+1}} \otimes X_{2\varepsilon_1} \in \frak{l}_{n\gamma} \otimes {\mathfrak
z}({\mathfrak n})$.
Thus, $\pr_{\frak{l}_{\gamma} \otimes {\mathfrak z}({\mathfrak n})}(\ad(X_{\mu})^2 \otimes
\Id) \omega) = 0$.
It can be shown similarly that $\pr_{\frak{l}_{\gamma} \otimes {\mathfrak z}({\mathfrak
n})}(\ad(X_{\epsilon_\gamma})^2 \otimes \Id) \omega) = 0$.
Now the proposed equality follows from~\eqref{Eqn:4.7}.
\end{proof}
\begin{Prop}
\label{Prop:HWeight}
The vector $\pr_{\frak{l}_{\gamma}\otimes{\mathfrak z}({\mathfrak n})}(\tau_2(X_{\mu} +
X_{\epsilon_\gamma}))$ is a~highest weight vector for $V(\mu+\epsilon_\gamma)$.
\end{Prop}
\begin{proof}
We start with showing that $\pr_{\frak{l}_{\gamma}\otimes{\mathfrak z}({\mathfrak n})}(\tau_2(X_{\mu}
+ X_{\epsilon_\gamma}))$ has weight $\mu+\epsilon_\gamma$.
Since the ${\mathfrak l}$-action commutes with the projection map $\pr_{\frak{l}_{\gamma} \otimes
{\mathfrak z}({\mathfrak n})}$, by Lemma~\ref{Lem:4.6}, for $H \in {\mathfrak h} \subset {\mathfrak l}$, we
have
\begin{gather*}
H\cdot\pr_{\frak{l}_{\gamma}\otimes{\mathfrak z}({\mathfrak n})}
(\tau_2(X_{\mu}+X_{\epsilon_\gamma}))
\\
\qquad
=\frac{1}{2}\Big(\pr_{\frak{l}_{\gamma}\otimes{\mathfrak z}({\mathfrak n})}(H\cdot\ad
(X_{\mu})\ad(X_{\epsilon_\gamma})\otimes\Id)\omega
+\pr_{\frak{l}_{\gamma}\otimes{\mathfrak z}({\mathfrak n})}
(H\cdot\ad(X_{\epsilon_\gamma})\ad(X_{\mu})\otimes\Id)
\omega\Big)
\\
\qquad
=\frac{1}{2}\Big(
\pr_{\frak{l}_{\gamma}\otimes{\mathfrak z}({\mathfrak n})}
\big(\big(\ad([H,X_{\mu}])\ad(X_{\epsilon_\gamma})+\ad(X_{\mu})\ad([H,X_{\epsilon_\gamma}])\big)\otimes\Id\big)
\omega
\\
\qquad\phantom{=}
+\pr_{\frak{l}_{\gamma}\otimes{\mathfrak z}({\mathfrak n})}
\big((\ad([H,X_{\epsilon_\gamma}])\ad(X_{\mu})+\ad(X_{\epsilon_\gamma})\ad([H,X_{\mu}]))\otimes\Id\big)
\omega\Big)
\\
\qquad
=\frac{(\mu+\epsilon_\gamma)(H)}{2}
\Big(\pr_{\frak{l}_{\gamma}\otimes{\mathfrak z}({\mathfrak n})}
(\ad(X_{\mu})\ad(X_{\epsilon_\gamma})\otimes\Id)\omega
+\pr_{\frak{l}_{\gamma}\otimes{\mathfrak z}({\mathfrak n})}
(\ad(X_{\epsilon_\gamma})\ad(X_{\mu})\otimes\Id)
\omega\Big)
\\
\qquad
=(\mu+\epsilon_\gamma)(H)\pr_{\frak{l}_{\gamma}\otimes{\mathfrak z}({\mathfrak n})}
(\tau_2(X_{\mu}+X_{\epsilon_\gamma})).
\end{gather*}
Note that Lemma~\ref{Lem:4.4} is applied from line two to line three.

Next we show that $\pr_{\frak{l}_{\gamma}\otimes{\mathfrak z}({\mathfrak n})}(\tau_2(X_{\mu} +
X_{\epsilon_\gamma}))$ is a~highest weight vector.
Let $\alpha \in \Pi({\mathfrak l})$.
As the ${\mathfrak l}$-action commutes with $\pr_{\frak{l}_{\gamma} \otimes {\mathfrak z}({\mathfrak
n})}$, we f\/irst observe $X_\alpha \cdot \tau_2(X_{\mu} + X_{\epsilon_\gamma})$.
It follows from~\eqref{Eqn:4.71} that $X_{\alpha} \cdot \tau_2(X_{\mu}+X_{\epsilon_\gamma})$ is given by
\begin{gather*}
X_{\alpha}\cdot\tau_2(X_{\mu}+X_{\epsilon_\gamma})
=\frac{1}{2}\Big(X_\alpha\cdot\big(\ad(X_{\mu})^2\otimes\Id\big)\omega+X_\alpha\cdot(\ad
(X_{\mu})\ad(X_{\epsilon_\gamma})\otimes\Id)\omega
\\
\qquad\phantom{=} +X_\alpha\cdot(\ad(X_{\epsilon_\gamma})\ad(X_{\mu})\otimes\Id
)\omega+X_\alpha\cdot\big(\ad(X_{\epsilon_\gamma})^2\otimes\Id\big)\omega\Big).
\end{gather*}
If $Z = X_{\alpha}$ in~\eqref{Eqn:Technical} then, as $[X_{\alpha}, X_{\mu}] = 0$,
we obtain
\begin{gather}
X_{\alpha}\cdot\tau_2(X_{\mu}+X_{\epsilon_\gamma})
\nonumber
\\
\qquad =\frac{1}{2}\big(\ad(X_{\mu})\ad([X_{\alpha},X_{\epsilon_\gamma}])\otimes\Id
\big)\omega+\big(\ad([X_{\alpha},X_{\epsilon_\gamma}])\ad(X_{\mu})\otimes\Id
\big)\omega
\nonumber
\\
\qquad\phantom{=} +\big(\ad([X_{\alpha},X_{\epsilon_\gamma}])\ad(X_{\epsilon_\gamma})\otimes\Id
\big)\omega+\big(\ad(X_{\epsilon_\gamma})\ad([X_{\alpha},X_{\epsilon_\gamma}])\otimes\Id
\big)\omega.
\label{Eqn:LVect}
\end{gather}
Recall from Tables~2 and~4 in~\cite[Section 6]{KuboThesis1} that we have
$\mu+\epsilon_\gamma = \varepsilon_1 + \varepsilon_2$ with $\mu=\varepsilon_1 + \varepsilon_{i+1}$ and
$\epsilon_\gamma = \varepsilon_2 - \varepsilon_{i+1}$.
Since $\Pi({\mathfrak l}) = \{\varepsilon_j - \varepsilon_{j+1} \;:\; 1 \leq j \leq n-1$ with $j
\neq i \} \cup \{2\varepsilon_n\}$, it follows that $\alpha + \epsilon_\gamma \in \Delta$ if and only if
$\alpha = \varepsilon_1 - \varepsilon_2$.
So it suf\/f\/ices to consider the case that $\alpha = \varepsilon_1 - \varepsilon_2$.
As $(\varepsilon_1 - \varepsilon_2) + 2(\varepsilon_2 - \varepsilon_{i+1}) \notin \Delta$, we have
$\ad(X_{\varepsilon_1-\varepsilon_2})\ad(X_{\varepsilon_2 - \varepsilon_{i+1}})
=\ad(X_{\varepsilon_2 - \varepsilon_{i+1}})\ad(X_{\varepsilon_1-\varepsilon_2})$.
Therefore, if $\alpha = \varepsilon_1 - \varepsilon_2$ then, as $\mu = \varepsilon_1 + \varepsilon_{i+1}$
and $\epsilon_\gamma + \alpha = \varepsilon_1 - \varepsilon_{i+1}$,~\eqref{Eqn:LVect} becomes
\begin{gather*}
X_{\varepsilon_1-\varepsilon_2}\cdot\tau_2(X_{\varepsilon_1+\varepsilon_{i+1}}
+X_{\varepsilon_2-\varepsilon_{i+1}})
\\
\qquad
=\frac{N_{\varepsilon_1-\varepsilon_2,\varepsilon_2-\varepsilon_{i+1}}}{2}\bigg(\big(\ad
(X_{\varepsilon_1+\varepsilon_{i+1}})\ad(X_{\varepsilon_1-\varepsilon_{i+1}})\otimes\Id
\big)\omega
\\
\qquad\phantom{=} +\big(\ad(X_{\varepsilon_1-\varepsilon_{i+1}})\ad(X_{\varepsilon_1+\varepsilon_{i+1}}
)\otimes\Id\big)\omega
+\big(2\ad(X_{\varepsilon_2-\varepsilon_{i+1}})\ad
(X_{\varepsilon_1-\varepsilon_{i+1}})\otimes\Id\big)\omega\bigg).
\end{gather*}
Table~\ref{T:Constants2} shows that $N_{\varepsilon_1 - \varepsilon_2, \varepsilon_2-\varepsilon_{i+1}} =
1$.
Therefore,
\begin{gather*}
X_{\varepsilon_1-\varepsilon_2}\cdot\tau_2(X_{\varepsilon_1+\varepsilon_{i+1}}
+X_{\varepsilon_2-\varepsilon_{i+1}})
\\
\qquad =\frac{1}{2}\bigg(\big(\ad(X_{\varepsilon_1+\varepsilon_{i+1}})\ad
(X_{\varepsilon_1-\varepsilon_{i+1}})\otimes\Id\big)\omega+\big(\ad
(X_{\varepsilon_1-\varepsilon_{i+1}})\ad(X_{\varepsilon_1+\varepsilon_{i+1}})\otimes\Id
\big)\omega
\\
\qquad\phantom{=} +\big(2\ad(X_{\varepsilon_2-\varepsilon_{i+1}})\ad
(X_{\varepsilon_1-\varepsilon_{i+1}})\otimes\Id\big)\omega\bigg).
\end{gather*}

Now we consider the contribution from each term separately.
A direct computation shows that the f\/irst term is given by
\begin{gather*}
T_1=\big(\ad(X_{\varepsilon_1+\varepsilon_{i+1}})\ad(X_{\varepsilon_1-\varepsilon_{i+1}}
)\otimes\Id\big)\omega
\\
\hphantom{T_1}{} =\frac{1}{c_0}\sum_{j=1}^i\ad(X_{\varepsilon_1+\varepsilon_{i+1}})\ad
(X_{\varepsilon_1-\varepsilon_{i+1}})X_{-2\varepsilon_j}\otimes X_{2\varepsilon_j}
\\
\hphantom{T_1=}{} +\frac{1}{2c_0}\sum_{1\leq j<k\leq i}\ad(X_{\varepsilon_1+\varepsilon_{i+1}})\ad
(X_{\varepsilon_1-\varepsilon_{i+1}})X_{-(\varepsilon_j+\varepsilon_k)}
\otimes X_{\varepsilon_j+\varepsilon_k}
\\
\hphantom{T_1}{} =\frac{1}{c_0}N_{\varepsilon_1-\varepsilon_{i+1},-2\varepsilon_1}H_{\varepsilon_1+\varepsilon_{i+1}}
\otimes X_{2\varepsilon_1}
\\
\hphantom{T_1=}{} +\frac{1}{2c_0}\sum_{k=2}^i N_{\varepsilon_1-\varepsilon_{i+1}
,-(\varepsilon_1+\varepsilon_k)}N_{\varepsilon_1+\varepsilon_{i+1},-(\varepsilon_k+\varepsilon_{i+1})}
X_{\varepsilon_1-\varepsilon_k}\otimes X_{\varepsilon_1+\varepsilon_k}
\\
\hphantom{T_1}{} =\frac{-1}{c_0}H_{\varepsilon_1+\varepsilon_{i+1}}\otimes X_{2\varepsilon_1}-\frac{1}{2c_0}\sum_{k=2}
^i X_{\varepsilon_1-\varepsilon_k}\otimes X_{\varepsilon_1+\varepsilon_k}.
\end{gather*}
Similarly we have
\begin{gather*}
T_2=\big(\ad(X_{\varepsilon_1-\varepsilon_{i+1}})\ad(X_{\varepsilon_1+\varepsilon_{i+1}}
)\otimes\Id\big)\omega=\frac{1}{c_0}H_{\varepsilon_1-\varepsilon_{i+1}}\otimes X_{2\varepsilon_1}
+\frac{1}{2c_0}\sum_{k=2}^i X_{\varepsilon_1-\varepsilon_k}\otimes X_{\varepsilon_1+\varepsilon_k}
\end{gather*}
and
\begin{gather*}
T_3=\big(2\ad(X_{\varepsilon_2-\varepsilon_{i+1}})\ad(X_{\varepsilon_1-\varepsilon_{i+1}}
)\otimes\Id\big)\omega=\frac{4}{c_0}X_{2\varepsilon_{i+1}}\otimes X_{-(\varepsilon_1+\varepsilon_2)}.
\end{gather*}
Therefore,
\begin{gather*}%\label{Eqn:4.9}
X_{\varepsilon_1-\varepsilon_2}\cdot\tau_2(X_{\varepsilon_1+\varepsilon_{i+1}}
+X_{\varepsilon_2-\varepsilon_{i+1}})
=T_1+T_2+T_3
\\
\qquad =\frac{-1}{c_0}H_{\varepsilon_1+\varepsilon_{i+1}}\otimes X_{2\varepsilon_1}-\frac{1}{2c_0}\sum_{k=2}
^i X_{\varepsilon_1-\varepsilon_k}\otimes X_{\varepsilon_1+\varepsilon_k}
\\
\qquad\phantom{=} +\frac{1}{c_0}H_{\varepsilon_1-\varepsilon_{i+1}}\otimes X_{2\varepsilon_1}+\frac{1}{2c_0}\sum_{k=2}
^i X_{\varepsilon_1-\varepsilon_k}\otimes X_{\varepsilon_1+\varepsilon_k}+\frac{4}{c_0}X_{2\varepsilon_{i+1}
}\otimes X_{-(\varepsilon_1+\varepsilon_2)}
\\
\qquad =\frac{1}{c_0}(H_{\varepsilon_1-\varepsilon_{i+1}}-H_{\varepsilon_1+\varepsilon_{i+1}}
)\otimes X_{2\varepsilon_1}+\frac{4}{c_0}X_{2\varepsilon_{i+1}}\otimes X_{-(\varepsilon_1+\varepsilon_2)}.
\end{gather*}

Observe that ${\mathfrak h} \cap \frak{l}_{\gamma}$ is spanned by the elements
$H_{\varepsilon_j - \varepsilon_{j+1}} =(E_{j j}-E_{\hat{j},\hat{j}}) - (E_{j+1,
j+1}-E_{\widehat{j+1},\widehat{j+1}})$ for $1\leq j \leq i-1$.
Since $H_{\varepsilon_1-\varepsilon_{i+1}} - H_{\varepsilon_1 + \varepsilon_{i+1}} = -2(E_{i+1,
i+1}-E_{\widehat{i+1},\widehat{i+1}})$, it follows that $H_{\varepsilon_1-\varepsilon_{i+1}} -
H_{\varepsilon_1 + \varepsilon_{i+1}} \notin {\mathfrak h} \cap \frak{l}_{\gamma}$.
As $\pr_{\frak{l}_{\gamma} \otimes {\mathfrak z}({\mathfrak n})}(X_{2\varepsilon_{i+1}} \otimes
X_{-(\varepsilon_1 + \varepsilon_2)}) = 0$, we then obtain
\begin{gather*}
X_{\varepsilon_1-\varepsilon_2}\cdot\pr_{\frak{l}_{\gamma}\otimes{\mathfrak z}({\mathfrak n})}
(\tau_2(X_{\varepsilon_1+\varepsilon_{i+1}}+X_{\varepsilon_2-\varepsilon_{i+1}}))
=\pr_{\frak{l}_{\gamma}\otimes{\mathfrak z}({\mathfrak n})}(X_{\varepsilon_1-\varepsilon_2}
\cdot\tau_2(X_{\varepsilon_1+\varepsilon_{i+1}}+X_{\varepsilon_2-\varepsilon_{i+1}}))
\\
 =\frac{1}{c_0}\pr_{\frak{l}_{\gamma}\otimes{\mathfrak z}({\mathfrak n})}
\big((H_{\varepsilon_1-\varepsilon_{i+1}}-H_{\varepsilon_1+\varepsilon_{i+1}})\otimes X_{2\varepsilon_1}
\big)+\frac{4}{c_0}\pr_{\frak{l}_{\gamma}\otimes{\mathfrak z}({\mathfrak n})}
\big(X_{2\varepsilon_{i+1}}\otimes X_{-(\varepsilon_1+\varepsilon_2)}\big)
 =0. \!\!\!\!\tag*{\qed}
\end{gather*}
\renewcommand{\qed}{}
\end{proof}

Now we def\/ine the ``opposite'' $\tau_2$ map $\bar{\tau}_2$ by
\begin{gather*}
\bar{\tau}_2 : \  {\mathfrak g}(-1)\to{\mathfrak g}(0)\otimes{\mathfrak g}(-2),
\qquad
X^* \mapsto\frac{1}{2}\big(\ad(X^*)^2\otimes\Id\big)\bar{\omega}
\end{gather*}
with
\begin{gather*}
\bar{\omega}=\sum_{\gamma_j\in{\mathfrak z}({\mathfrak n})}X_{\gamma}\otimes X^*_{\gamma}
=\frac{1}{c_0}\sum_{j=1}^i X_{2\varepsilon_j}\otimes X_{-2\varepsilon_j}+\frac{1}{2c_0}
\sum_{1\leq j<k\leq i}X_{\varepsilon_j+\varepsilon_k}\otimes X_{-(\varepsilon_j+\varepsilon_k)}.
\end{gather*}
It follows from the same arguments in the proof for Lemma 3.3 and Proposition 7.3 in~\cite{KuboThesis1}
that the $\bar{\tau}_2$ map is not identically zero and $L$-equivariant.
Let $\pr_{\frak{l}_{\gamma} \otimes {\mathfrak z}(\bar{\mathfrak n})}: {\mathfrak l} \otimes
{\mathfrak z}(\bar{\mathfrak n}) \to \frak{l}_{\gamma} \otimes {\mathfrak z}(\bar{\mathfrak n})$ be the
projection map from ${\mathfrak l} \otimes {\mathfrak z}(\bar{\mathfrak n})$ onto $\frak{l}_{\gamma}
\otimes {\mathfrak z}(\bar{\mathfrak n})$.
\begin{Prop}
\label{Prop:LWeight}
The vector $\pr_{\frak{l}_{\gamma}\otimes{\mathfrak z}(\bar{\mathfrak
n})}(\bar{\tau}_2(X_{-\mu} + X_{-\epsilon_\gamma}))$ is a~lowest weight vector for
$V(\mu+\epsilon_\gamma)^*$.
\end{Prop}
\begin{proof}
This proposition immediately follows from the arguments used in the proof for
Proposition~\ref{Prop:HWeight}, by replacing positive (resp.\
negative) roots with negative (resp.\
positive) roots.
\end{proof}

We set
\begin{gather}
\label{Eqn:LV}
Y^*_l:=\frac{8c_0^2}{i+1}\pr_{\frak{l}_{\gamma}\otimes{\mathfrak z}(\bar{\mathfrak n})}(\bar{\tau}
_2(X_{-\mu}+X_{-\epsilon_\gamma})).
\end{gather}
It follows from Proposition~\ref{Prop:LWeight} that $Y^*_l$ is a~lowest weight vector for $V(\mu
+\epsilon_\gamma)^*$.
In the next subsection we compute $\tilde{\tau}_2(Y^*_l)(X)$.
To the end we give an explicit form for $Y^*_l$.
\begin{Lem}
\label{Lem:LV}
We have
\begin{gather*}
Y^*_l=\frac{4c_0}{i+1}X_{\varepsilon_1-\varepsilon_2}\otimes X_{-2\varepsilon_1}-\frac{4c_0}{i+1}
X_{-(\varepsilon_1-\varepsilon_2)}\otimes X_{-2\varepsilon_2}+\frac{2c_0}{i+1}\sum_{k=3}^i X_{-(\varepsilon_2-\varepsilon_k)}\otimes X_{-(\varepsilon_1+\varepsilon_k)}
\\
\hphantom{Y^*_l=}{}
-\frac{2c_0}{i+1}\sum_{k=3}^i X_{-(\varepsilon_1-\varepsilon_k)}\otimes X_{-(\varepsilon_2+\varepsilon_k)}
 -\frac{2c_0}{i+1}H_{\varepsilon_1-\varepsilon_2}\otimes X_{-(\varepsilon_1+\varepsilon_2)}.
\end{gather*}
\end{Lem}
\begin{proof}
By using the same arguments for Lemma~\ref{Lem:4.6}, one can obtain
\begin{gather}
\pr_{\frak{l}_{\gamma}\otimes{\mathfrak z}(\bar{\mathfrak n})}
(\tau_2(X_{-\mu}+X_{-\epsilon_\gamma}))
\nonumber
\\
\qquad
=\frac{1}{2}\Big(\pr_{\frak{l}_{\gamma}\otimes{\mathfrak z}(\bar{\mathfrak n})}
(\ad(X_{-\mu})\ad(X_{-\epsilon_\gamma})\otimes\Id)\bar{\omega}
+\pr_{\frak{l}_{\gamma}\otimes{\mathfrak z}(\bar{\mathfrak n})}
(\ad(X_{-\epsilon_\gamma})\ad(X_{-\mu})\otimes\Id)
\bar{\omega}\Big).  \!\!\!\!
\label{Eqn:4.14}
\end{gather}
A direct computation shows that
\begin{gather*}
\pr_{\frak{l}_{\gamma}\otimes{\mathfrak z}(\bar{\mathfrak n})}(\ad(X_{-\mu})\ad
(X_{-\epsilon_\gamma})\otimes\Id)\bar{\omega})
\\
\qquad =\frac{1}{c_0}\sum_{j=1}^i\pr_{\frak{l}_{\gamma}\otimes{\mathfrak z}(\bar{\mathfrak n})}
\big(\ad(X_{-(\varepsilon_1+\varepsilon_{i+1})})\ad(X_{-(\varepsilon_2-\varepsilon_{i+1})}
)X_{2\varepsilon_j}\otimes X_{-2\varepsilon_j}\big)
\\
\qquad\phantom{=} +\frac{1}{2c_0}\sum_{1\leq j<k\leq i}\pr_{\frak{l}_{\gamma}\otimes{\mathfrak z}(\bar{\mathfrak n})}
\big(\ad(X_{-(\varepsilon_1+\varepsilon_{i+1})})\ad(X_{-(\varepsilon_2-\varepsilon_{i+1})}
)X_{\varepsilon_j+\varepsilon_k}\otimes X_{-(\varepsilon_j+\varepsilon_k)}\big)
\\
\qquad =-\frac{1}{c_0}X_{-(\varepsilon_1-\varepsilon_2)}\otimes X_{-2\varepsilon_2}-\frac{1}{2c_0}
\pr_{\frak{l}_{\gamma}\otimes{\mathfrak z}(\bar{\mathfrak n})}(H_{\varepsilon_1+\varepsilon_{i+1}}
\otimes X_{-(\varepsilon_1+\varepsilon_2)})
\\
\qquad\phantom{=}
-\frac{1}{2c_0}\sum_{k=3}^i X_{-(\varepsilon_1-\varepsilon_k)}
\otimes X_{-(\varepsilon_2+\varepsilon_k)}.
\end{gather*}
Observe that we have
\begin{gather*}
H_{\varepsilon_1+\varepsilon_{i+1}}=\frac{1}{i}\sum_{j=1}^i(E_{jj}-E_{\hat{j}\hat{j}})+\frac{1}{i}\sum_{k=2}
^iH_{\varepsilon_1-\varepsilon_k}+H_{2\varepsilon_{i+1}}.
\end{gather*}
Since $(1/i)\sum\limits_{j=1}^i(E_{jj} - E_{\hat{j}\hat{j}}) \in {\mathfrak
z}({\mathfrak l})$ and $H_{2\varepsilon_{i+1}} \in {\mathfrak h} \cap \frak{l}_{n\gamma}$, it follows that
\begin{gather*}
\pr_{\frak{l}_{\gamma}\otimes{\mathfrak z}(\bar{\mathfrak n})}(H_{\varepsilon_1+\varepsilon_{i+1}}
\otimes X_{-(\varepsilon_1+\varepsilon_2)})=\frac{1}{i}\sum_{k=2}^iH_{\varepsilon_1-\varepsilon_k}
\otimes X_{-(\varepsilon_1+\varepsilon_2)}.
\end{gather*}
Therefore, $\pr_{\frak{l}_{\gamma} \otimes {\mathfrak z}(\bar{\mathfrak
n})}(\ad(X_{-\mu})\ad(X_{-\epsilon_\gamma}) \otimes \Id) \bar{\omega})$ is given by
\begin{gather}
\pr_{\frak{l}_{\gamma}\otimes{\mathfrak z}(\bar{\mathfrak n})}(\ad(X_{-\mu})\ad
(X_{-\epsilon_\gamma})\otimes\Id)\bar{\omega})
\nonumber
\\
\qquad =-\frac{1}{c_0}X_{-(\varepsilon_1-\varepsilon_2)}\otimes X_{-2\varepsilon_2}-\frac{1}{2ic_0}
\sum_{k=2}^iH_{\varepsilon_1-\varepsilon_k}\otimes X_{-(\varepsilon_1+\varepsilon_2)}
\nonumber
\\
\qquad\phantom{=}
-\frac{1}{2c_0}
\sum_{k=3}^i X_{-(\varepsilon_1-\varepsilon_k)}\otimes X_{-(\varepsilon_2+\varepsilon_k)}.
\label{Eqn:4.15}
\end{gather}
Similarly we have
\begin{gather}
\pr_{\frak{l}_{\gamma}\otimes{\mathfrak z}(\bar{\mathfrak n})}(\ad(X_{-\epsilon_\gamma}
)\ad(X_{-\mu})\otimes\Id)\bar{\omega})
\nonumber
\\
\qquad
=\frac{1}{c_0}X_{\varepsilon_1-\varepsilon_2}\otimes X_{-2\varepsilon_1}-\frac{1}{2ic_0}
(H_{\varepsilon_1-\varepsilon_2}-\sum_{k=3}^iH_{\varepsilon_2-\varepsilon_k})
\otimes X_{-(\varepsilon_1+\varepsilon_2)}
\nonumber
\\
\qquad\phantom{=}
+\frac{1}{2c_0}\sum_{k=3}^i X_{-(\varepsilon_2-\varepsilon_k)}
\otimes X_{-(\varepsilon_1+\varepsilon_k)}.
\label{Eqn:4.16}
\end{gather}
By substituting~\eqref{Eqn:4.15} and~\eqref{Eqn:4.16} into~\eqref{Eqn:4.14} and multiplying the resulting
equation by $8c_0^2/(i+1)$, one obtains
\begin{gather*}
Y^*_l=\frac{4c_0}{i+1}X_{\varepsilon_1-\varepsilon_2}\otimes X_{-2\varepsilon_1}-\frac{4c_0}{i+1}
X_{-(\varepsilon_1-\varepsilon_2)}\otimes X_{-2\varepsilon_2}
\\
\hphantom{Y^*_l=}{} +\frac{2c_0}{i+1}\sum_{k=3}^i X_{-(\varepsilon_2-\varepsilon_k)}\otimes X_{-(\varepsilon_1+\varepsilon_k)}
-\frac{2c_0}{i+1}\sum_{k=3}^i X_{-(\varepsilon_1-\varepsilon_k)}\otimes X_{-(\varepsilon_2+\varepsilon_k)}
\\
\hphantom{Y^*_l=}{} -\frac{2c_0}{i(i+1)}\big(\sum_{k=2}^i H_{\varepsilon_1-\varepsilon_k}+H_{\varepsilon_1-\varepsilon_2}
-\sum_{k=3}^iH_{\varepsilon_2-\varepsilon_k})\otimes X_{-(\varepsilon_1+\varepsilon_2)}.
\end{gather*}
Now the proposed equality follows from manipulating the elements in the Cartan subalgebra.
\end{proof}

\subsection[The $\tilde{\tau}_2|_{V(\mu+\epsilon_\gamma)^*}$ map]{The
$\boldsymbol{\tilde{\tau}_2|_{V(\mu+\epsilon_\gamma)^*}}$ map}
%\label{SS:Type3B}

Now we show that the map $\tilde{\tau}_2|_{V(\mu+\epsilon_\gamma)^*}$ is not identically zero.
To do so, we recall several essential ingredients.
First observe that, as the duality is with respect to the Killing form $\kappa$, if $Y^* = X_\alpha \otimes
X_{-\gamma_j} \in \frak{l}_{\gamma} \otimes {\mathfrak z}(\bar{\mathfrak n})$ and $X_\beta \otimes
X_{\gamma_k} \in \frak{l}_{\gamma} \otimes {\mathfrak z}({\mathfrak n})$ then $Y^*(X_\beta \otimes
X_{\gamma_k})$ is given by $Y^*(X_\beta \otimes X_{\gamma_k})= \kappa(X_\alpha,
X_\beta)\kappa(X_{-\gamma_k}, X_{\gamma_j})$.
As observed in~\eqref{Eqn:Kappa}, we have
\begin{gather*}
\kappa(X_\alpha,X_{-\alpha})=
\begin{cases}
c_0 & \text{if}~\alpha~\text{is long},
\\
2c_0 & \text{if}~\alpha~\text{is short}.
\end{cases}
\end{gather*}
Finally we recall from~\eqref{Eqn:Delta} that if $W \subset {\mathfrak g}$ is an
$\ad({\mathfrak h})$-invariant subspace then, for any weight $\nu \in {\mathfrak h}^*$, we write
$\Delta_{\nu}(W) = \{ \alpha \in \Delta(W) \, | \, \nu - \alpha \in \Delta\}$.
\begin{Prop}
\label{Prop:Ttau}
The $L$-intertwining operator $\tilde{\tau}_2|_{V(\mu+\epsilon_\gamma)^*}$ is not identically zero.
\end{Prop}
\begin{proof}
Take lowest weight vector $Y_l^*$ as in~\eqref{Eqn:LV}.
We show that $\tilde{\tau}_2(Y^*_l)(X)$ is a~non-zero polynomial on ${\mathfrak g}(1)$.
As $\tau_2(X) = (1/2)\sum\limits_{\gamma_j \in \Delta({\mathfrak z}({\mathfrak n}))}
\ad(X)^2X^*_{\gamma_j}\otimes X_{\gamma_j}$, by Lemma~\ref{Lem:LV}, the polynomial
$\tilde{\tau}_2(Y^*_l)(X)$ may express as a~sum of f\/ive terms.
We consider the contribution from each term separately, and start with observing the contribution from the f\/irst term   % wewewe
%We have
\begin{gather}
T_1=\frac{2c_0}{i+1}\sum_{\gamma_j\in\Delta({\mathfrak z}({\mathfrak n}))}
\kappa\big(X_{\varepsilon_1-\varepsilon_2},\ad(X)^2X^*_{\gamma_j}\big)\kappa(X_{-2\varepsilon_1},X_{\gamma_j})
\nonumber
\\
\phantom{T_1}
=\frac{2c_0}{i+1}\kappa\big(X_{\varepsilon_1-\varepsilon_2},\ad(X)^2X_{-2\varepsilon_1}\big).
\label{Eqn:T1_1}
\end{gather}
Write $X = \sum\limits_{\alpha \in \Delta({\mathfrak g}(1))} \eta_\alpha X_{\alpha}$, where $\eta_\alpha
\in {\mathfrak n}^*$ is the coordinate dual to $X_\alpha$ with respect to the Killing form.
Then,
\begin{gather}
\eqref{Eqn:T1_1}=\frac{2c_0}{i+1}\kappa\big(X_{\varepsilon_1-\varepsilon_2},\ad
(X)^2X_{-2\varepsilon_1}\big)
\nonumber
\\
\phantom{\eqref{Eqn:T1_1}}
=\frac{2c_0}{i+1}\sum_{\alpha,\beta\in\Delta({\mathfrak g}(1))}
\eta_\alpha\eta_\beta\kappa(X_{\varepsilon_1-\varepsilon_2},[X_\beta,[X_\alpha,X_{-2\varepsilon_1}]])
\nonumber
\\
\phantom{\eqref{Eqn:T1_1}}
=\frac{2c_0}{i+1}\sum_{\alpha,\beta\in\Delta({\mathfrak g}(1))}
\eta_\alpha\eta_\beta\kappa([X_{\varepsilon_1-\varepsilon_2},X_\beta],[X_\alpha,X_{-2\varepsilon_1}])
\nonumber
\\
\phantom{\eqref{Eqn:T1_1}}
=\frac{2c_0}{i+1}\sum_{\substack{\alpha\in\Delta_{2\varepsilon_1}({\mathfrak g}(1))\\
\beta\in\Delta^{\varepsilon_1-\varepsilon_2}({\mathfrak g}(1))}}
N_{\varepsilon_1-\varepsilon_2,\beta}
N_{\alpha,-2\varepsilon_1}\eta_\alpha\eta_\beta\kappa(X_{\beta+(\varepsilon_1-\varepsilon_2)}
,X_{\alpha-2\varepsilon_1}),
\label{Eqn:T1_2}
\end{gather}
where $\Delta^{\varepsilon_1-\varepsilon_2}({\mathfrak g}(1))
= \{ \alpha \in \Delta({\mathfrak g}(1)) \, | \, (\varepsilon_1 - \varepsilon_2) + \alpha \in \Delta\}$.
Observe that $\kappa(X_{\beta+(\varepsilon_1-\varepsilon_2)}, X_{\alpha-2\varepsilon_1}) \neq 0$ if and
only if $\alpha \in \Delta_{2\varepsilon_1}({\mathfrak g}(1))$.
Indeed, f\/irst one may see from~\eqref{Eqn:T1_2} that $\kappa(X_{\beta+(\varepsilon_1-\varepsilon_2)},
X_{\alpha-2\varepsilon_1}) \neq 0$ if and only if $\beta = (\varepsilon_1 + \varepsilon_2) -\alpha$;
equivalently, the value of the Killing form is non-zero if and only if $\alpha \in
\Delta_{2\varepsilon_1}({\mathfrak g}(1))\cap \Delta_{\varepsilon_1+\varepsilon_2}({\mathfrak g}(1))$.
By inspection, we have $\Delta_{2\varepsilon_1}({\mathfrak g}(1)) = \{ \varepsilon_1 \pm \varepsilon_j \; |
\; i+1 \leq j \leq n\}$.
Thus, for any $\alpha \in \Delta_{2\varepsilon_1}({\mathfrak g}(1))$, it follows that
$(\varepsilon_1+\varepsilon_2)-\alpha \in \Delta$.
Therefore $\Delta_{2\varepsilon_1}({\mathfrak g}(1))\cap \Delta_{\varepsilon_1+\varepsilon_2}({\mathfrak
g}(1)) = \Delta_{2\varepsilon_1}({\mathfrak g}(1))$.
Now, if $(\alpha, \beta)$ denotes a~pair so that $\kappa(X_{\beta+(\varepsilon_1-\varepsilon_2)},
X_{\alpha-2\varepsilon_1}) \neq 0$ then $(\alpha, \beta) = (\varepsilon_1\pm \varepsilon_j, \varepsilon_2
\mp \varepsilon_j)$ for $i+1 \leq j \leq n$ with respect to the signs.
Therefore,
\begin{gather*}
\eqref{Eqn:T1_2}=\frac{2c_0}{i+1}\sum_{\substack{\alpha\in\Delta_{2\varepsilon_1}({\mathfrak g}(1))\\
\beta\in\Delta^{\varepsilon_1-\varepsilon_2}({\mathfrak g}(1))}}N_{\varepsilon_1-\varepsilon_2,\beta}
N_{\alpha,-2\varepsilon_1}\eta_\alpha\eta_\beta\kappa(X_{\beta+(\varepsilon_1-\varepsilon_2)}
,X_{\alpha-2\varepsilon_1})
\\
\phantom{\eqref{Eqn:T1_2}}
=\frac{2c_0}{i+1}\sum_{\alpha\in\Delta_{2\varepsilon_1}({\mathfrak g}(1))}
N_{\varepsilon_1-\varepsilon_2,(\varepsilon_1+\varepsilon_2)-\alpha}N_{\alpha,-2\varepsilon_1}
\eta_\alpha\eta_{(\varepsilon_1+\varepsilon_2)-\alpha}\kappa(X_{2\varepsilon_1-\alpha}
,X_{\alpha-2\varepsilon_1})
\\
\phantom{\eqref{Eqn:T1_2}}
=\frac{2c_0}{i+1}\sum_{j=i+1}^n N_{\varepsilon_1-\varepsilon_2,\varepsilon_2-\varepsilon_j}
N_{\varepsilon_1+\varepsilon_j,-2\varepsilon_1}\eta_{\varepsilon_1+\varepsilon_j}
\eta_{\varepsilon_2-\varepsilon_j}\kappa(X_{\varepsilon_1-\varepsilon_j},X_{-(\varepsilon_1-\varepsilon_j)})
\\
\phantom{\eqref{Eqn:T1_2}}
\phantom{=} +\frac{2c_0}{i+1}\sum_{j=i+1}^n N_{\varepsilon_1-\varepsilon_2,\varepsilon_2+\varepsilon_j}
N_{\varepsilon_1-\varepsilon_j,-2\varepsilon_1}\eta_{\varepsilon_1-\varepsilon_j}
\eta_{\varepsilon_2+\varepsilon_j}\kappa(X_{\varepsilon_1+\varepsilon_j},X_{-(\varepsilon_1+\varepsilon_j)})
\\
\phantom{\eqref{Eqn:T1_2}}
=\frac{2c_0}{i+1}\sum_{j=i+1}^n(1)(1)\eta_{\varepsilon_1+\varepsilon_j}
\eta_{\varepsilon_2-\varepsilon_j}\kappa(X_{\varepsilon_1-\varepsilon_j},X_{-(\varepsilon_1-\varepsilon_j)})
\\
\phantom{\eqref{Eqn:T1_2}}
\phantom{=} +\frac{2c_0}{i+1}\sum_{j=i+1}^n(1)(-1)\eta_{\varepsilon_1-\varepsilon_j}\eta_{\varepsilon_2+\varepsilon_j}
\kappa(X_{\varepsilon_1+\varepsilon_j},X_{-(\varepsilon_1+\varepsilon_j)})
\\
\phantom{\eqref{Eqn:T1_2}}
=\frac{4c_0^2}{i+1}\sum_{j=i+1}^n\eta_{\varepsilon_1+\varepsilon_j}\eta_{\varepsilon_2-\varepsilon_j}
-\eta_{\varepsilon_1-\varepsilon_j}\eta_{\varepsilon_2+\varepsilon_j}
\\
\phantom{\eqref{Eqn:T1_2}}
=\frac{4c_0^2}{i+1}\sum_{j=i+1}^n\kappa(X,X^*_{\varepsilon_1+\varepsilon_j}
)\kappa(X,X^*_{\varepsilon_2-\varepsilon_j})-\kappa(X,X^*_{\varepsilon_1-\varepsilon_j}
)\kappa(X,X^*_{\varepsilon_2+\varepsilon_j}).
\end{gather*}

By a~similar computation one obtains
\begin{gather*}
T_2=\frac{-2c_0}{i+1}\sum_{\gamma_j\in\Delta({\mathfrak z}({\mathfrak n}))}
\kappa\big(X_{-(\varepsilon_1-\varepsilon_2)},\ad(X)^2X^*_{\gamma_j}\big)\kappa(X_{-2\varepsilon_2}
,X_{\gamma_j})
\\[-0.5mm]
\phantom{T_2}
=\frac{4c_0^2}{i+1}\sum_{j=i+1}^n\kappa(X,X^*_{\varepsilon_1+\varepsilon_j}
)\kappa(X,X^*_{\varepsilon_2-\varepsilon_j})-\kappa(X,X^*_{\varepsilon_1-\varepsilon_j}
)\kappa(X,X^*_{\varepsilon_2+\varepsilon_j}),
\\[-0.5mm]
T_3=\frac{c_0}{i+1}\sum_{\gamma_j\in\Delta({\mathfrak z}({\mathfrak n}))}\sum_{k=3}
^i\kappa\big(X_{-(\varepsilon_2-\varepsilon_k)},\ad(X)^2X^*_{\gamma_j}
\big)\kappa(X_{-(\varepsilon_1+\varepsilon_k)},X_{\gamma_j})
\\[-0.5mm]
\phantom{T_3}
=\frac{2(i-2)c_0^2}{i+1}\sum_{j=i+1}^n\kappa(X,X^*_{\varepsilon_1+\varepsilon_j}
)\kappa(X,X^*_{\varepsilon_2-\varepsilon_j})-\kappa(X,X^*_{\varepsilon_1-\varepsilon_j}
)\kappa(X,X^*_{\varepsilon_2+\varepsilon_j}),
\\[-0.5mm]
T_4=\frac{-c_0}{i+1}\sum_{\gamma_j\in\Delta({\mathfrak z}({\mathfrak n}))}\sum_{k=3}
^i\kappa\big(X_{-(\varepsilon_1-\varepsilon_k)},\ad(X)^2X^*_{\gamma_j}
\big)\kappa(X_{-(\varepsilon_2+\varepsilon_k)},X_{\gamma_j})
\\[-0.5mm]
\phantom{T_4}
=\frac{2(i-2)c_0^2}{i+1}\sum_{j=i+1}^n\kappa(X,X^*_{\varepsilon_1+\varepsilon_j}
)\kappa(X,X^*_{\varepsilon_2-\varepsilon_j})-\kappa(X,X^*_{\varepsilon_1-\varepsilon_j}
)\kappa(X,X^*_{\varepsilon_2+\varepsilon_j}),
\end{gather*}
and
\begin{gather*}
T_5=\frac{-c_0}{i+1}\sum_{\gamma_j\in\Delta({\mathfrak z}({\mathfrak n}))}
\kappa\big(H_{\varepsilon_1-\varepsilon_2},\ad(X)^2X^*_{\gamma_j}
\big)\kappa(X_{-(\varepsilon_1+\varepsilon_2)},X_{\gamma_j})
\\[-0.5mm]
\phantom{T_5}
=\frac{4c_0^2}{i+1}\sum_{j=i+1}^n\kappa(X,X^*_{\varepsilon_1+\varepsilon_j}
)\kappa(X,X^*_{\varepsilon_2-\varepsilon_j})-\kappa(X,X^*_{\varepsilon_1-\varepsilon_j}
)\kappa(X,X^*_{\varepsilon_2+\varepsilon_j}).
\end{gather*}

Therefore $\tilde{\tau}_2(Y^*_l)(X)$ may be given as
\begin{gather}
\tilde{\tau}_2\big(Y^*_l\big)(X)=T_1+T_2+T_3+T_4+T_5
\nonumber
\\[-0.5mm]
\hphantom{\tilde{\tau}_2\big(Y^*_l\big)(X)}{} =\frac{4c_0^2}{i+1}\sum_{j=i+1}^n\kappa(X,X^*_{\varepsilon_1+\varepsilon_j}
)\kappa(X,X^*_{\varepsilon_2-\varepsilon_j})-\kappa(X,X^*_{\varepsilon_1-\varepsilon_j}
)\kappa(X,X^*_{\varepsilon_2+\varepsilon_j})
\nonumber
\\[-0.5mm]
\hphantom{\tilde{\tau}_2\big(Y^*_l\big)(X)=}{} +\frac{4c_0^2}{i+1}\sum_{j=i+1}^n\kappa(X,X^*_{\varepsilon_1+\varepsilon_j}
)\kappa(X,X^*_{\varepsilon_2-\varepsilon_j})-\kappa(X,X^*_{\varepsilon_1-\varepsilon_j}
)\kappa(X,X^*_{\varepsilon_2+\varepsilon_j})
\nonumber
\\[-0.5mm]
\hphantom{\tilde{\tau}_2\big(Y^*_l\big)(X)=}{} +\frac{2(i-2)c_0^2}{i+1}\sum_{j=i+1}^n\kappa(X,X^*_{\varepsilon_1+\varepsilon_j}
)\kappa(X,X^*_{\varepsilon_2-\varepsilon_j})-\kappa(X,X^*_{\varepsilon_1-\varepsilon_j}
)\kappa(X,X^*_{\varepsilon_2+\varepsilon_j})
\nonumber
\\[-0.5mm]
\hphantom{\tilde{\tau}_2\big(Y^*_l\big)(X)=}{} +\frac{2(i-2)c_0^2}{i+1}\sum_{j=i+1}^n\kappa(X,X^*_{\varepsilon_1+\varepsilon_j}
)\kappa(X,X^*_{\varepsilon_2-\varepsilon_j})-\kappa(X,X^*_{\varepsilon_1-\varepsilon_j}
)\kappa(X,X^*_{\varepsilon_2+\varepsilon_j})
\nonumber
\\[-0.5mm]
\hphantom{\tilde{\tau}_2\big(Y^*_l\big)(X)=}{} +\frac{4c_0^2}{i+1}\sum_{j=i+1}^n\kappa(X,X^*_{\varepsilon_1+\varepsilon_j}
)\kappa(X,X^*_{\varepsilon_2-\varepsilon_j})-\kappa(X,X^*_{\varepsilon_1-\varepsilon_j}
)\kappa(X,X^*_{\varepsilon_2+\varepsilon_j})
\nonumber
\\[-0.5mm]
\hphantom{\tilde{\tau}_2\big(Y^*_l\big)(X)}{} =4c_0^2\sum_{j=i+1}^n\kappa(X,X^*_{\varepsilon_1+\varepsilon_j}
)\kappa(X,X^*_{\varepsilon_2-\varepsilon_j})-\kappa(X,X^*_{\varepsilon_1-\varepsilon_j}
)\kappa(X,X^*_{\varepsilon_2+\varepsilon_j}).
\label{Eqn:TauCn}
\end{gather}
Hence $\tilde{\tau}_2(Y^*_l)(X)$ is a~non-zero polynomial on ${\mathfrak g}(1)$.
\end{proof}

\subsection{The special value}
%\label{SS:Type3C}

Now we are going to f\/ind the special value for the $\Omega_2|_{V(\mu+\epsilon_\gamma)^*}$ system.
As for  Type~1b case, to f\/ind the special value, we use Proposition~\ref{Prop:HL}.
Recall from Section~\ref{SS:GVM-I} the linear map $\omega_2|_{V(\mu+\epsilon_\gamma)^*}:
V(\mu+\epsilon_\gamma)^* \to \mathcal{U}(\bar {\mathfrak n})$ def\/ined by
$\omega_2|_{V(\mu+\epsilon_{n\gamma})^*} = \sigma \circ \tilde{\tau}_2|_{V(\mu+\epsilon_{n\gamma})^*}$,
where $\sigma: \Sym^2({\mathfrak g}(-1)) \to \mathcal{U}(\bar {\mathfrak n})$ is the symmetrization
operator.
If $Y^*_l$ is the lowest weight vector for $V(\mu+\epsilon_\gamma)^*$ def\/ined in~\eqref{Eqn:LV} then it
follows from~\eqref{Eqn:TauCn} that $\omega(Y^*_l):=\omega_{V(\mu+\epsilon_\gamma)^*}$ is given by
\begin{gather*}
\omega\big(Y^*_l\big)=4c_0^2\sum_{j=i+1}^n\sigma(X^*_{\varepsilon_1+\varepsilon_j}
X^*_{\varepsilon_2-\varepsilon_j})-\sigma(X^*_{\varepsilon_1-\varepsilon_j}X^*_{\varepsilon_2+\varepsilon_j}
).
\end{gather*}
By~\eqref{Eqn:DVect}, this amounts to
\begin{gather*}
\label{Eqn:omegaCn}
\omega(Y^*_l)=\sum_{j=i+1}^n\sigma\big(X_{-(\varepsilon_1+\varepsilon_j)}X_{-(\varepsilon_2-\varepsilon_j}
)\big)-\sigma\big(X_{-(\varepsilon_1-\varepsilon_j)}X_{-(\varepsilon_2+\varepsilon_j)}\big).
\end{gather*}

The following lemma will simplify arguments for a~proof for Theorem~\ref{Thm:SpValType3} below.
\begin{Lem}
\label{Lem:Sym}
For $X, Y, Z \in {\mathfrak g}$, in $\mathcal{U}({\mathfrak g})$, we have
\begin{gather*}
X\cdot\sigma(YZ)=\sigma([X,Y]Z)+\sigma(Y[X,Z]).
\end{gather*}
\end{Lem}
\begin{proof}
A direct computation.
\end{proof}

Now we are ready to determine the special value for the $\Omega_2|_{V(\mu+\epsilon_\gamma)^*}$ system.
\begin{Thm}
\label{Thm:SpValType3}
Let ${\mathfrak q}$ be the maximal parabolic subalgebra of type~$C_n(i)$ for $2\leq i \leq n-1$.
The $\Omega_2|_{V(\mu+\epsilon_\gamma)^*}$ system is conformally invariant on $\mathcal{L}_s$ if and only
if $s = n-i+1$.
\end{Thm}
\begin{proof}
By Proposition~\ref{Prop:HL}, to prove this theorem, it suf\/f\/ices to show that $X_{\mu} \cdot
(\omega(Y^*_l)\otimes 1_{-s})=0$ in $M_{\mathfrak q}(\mathbb{C}_{-s})$ if and only if $s = n-i+1$, where
$\mu$ is the highest weight for ${\mathfrak g}(1)$.
It follows from~\eqref{Eqn:omegaCn} that $X_{\mu} \cdot (\omega(Y^*_l)\otimes 1_{-s})$ may be a~sum of two
terms.
As $\mu = \varepsilon_1 + \varepsilon_{i+1}$, the f\/irst term is
\begin{gather*}
T_1=\sum_{j=i+1}^n\big(X_{\varepsilon_1+\varepsilon_{i+1}}\cdot\sigma(X_{-(\varepsilon_1+\varepsilon_j)}
X_{-(\varepsilon_2-\varepsilon_j)})\big)\otimes1_{-s}.
\end{gather*}
By Lemma~\ref{Lem:Sym}, this may be expressed as
\begin{gather}
T_1=\sum_{j=i+1}^n\sigma([X_{\varepsilon_1+\varepsilon_{i+1}},X_{-(\varepsilon_1+\varepsilon_j)}
]X_{-(\varepsilon_2-\varepsilon_j)})\otimes1_{-s}
\nonumber
\\
\phantom{T_1=}{} +\sum_{j=i+1}^n\sigma(X_{-(\varepsilon_1+\varepsilon_j)}[X_{\varepsilon_1+\varepsilon_{i+1}}
,X_{-(\varepsilon_2-\varepsilon_j)}])\otimes1_{-s}
\nonumber
\\
\phantom{T_1}{}
=\sum_{j=i+1}^n\sigma([X_{\varepsilon_1+\varepsilon_{i+1}},X_{-(\varepsilon_1+\varepsilon_j)}
]X_{-(\varepsilon_2-\varepsilon_j)})
\nonumber
\\
\phantom{T_1}{}
=\sigma(H_{\varepsilon_1+\varepsilon_{i+1}}X_{-(\varepsilon_2-\varepsilon_j)})\otimes1_{-s}
\nonumber
\\
\phantom{T_1=}{}
+\sum_{j=i+2}^nN_{\varepsilon_1+\varepsilon_{i+1},-(\varepsilon_1+\varepsilon_j)}
\sigma(X_{-(\varepsilon_j-\varepsilon_{i+1})}X_{-(\varepsilon_2-\varepsilon_j)})\otimes1_{-s}.
\label{Eqn:Sp1}
\end{gather}
Since $\sigma(ab) = (1/2)(ab + ba)$, we have
\begin{gather}
\eqref{Eqn:Sp1}
=\sigma(H_{\varepsilon_1+\varepsilon_{i+1}}X_{-(\varepsilon_2{-}\varepsilon_j)})\!\otimes\! 1_{-s}\!
+\!\!\sum_{j=i+2}^nN_{\varepsilon_1{+}\varepsilon_{i+1},{-}(\varepsilon_1{+}\varepsilon_j)}
\sigma(X_{{-}(\varepsilon_j{-}\varepsilon_{i+1})}X_{-(\varepsilon_2{-}\varepsilon_j)})\!\otimes \!1_{-s}
\nonumber
\\
\phantom{\eqref{Eqn:Sp1}}
=\frac{1}{2}\big(H_{\varepsilon_1+\varepsilon_{i+1}}X_{-(\varepsilon_2-\varepsilon_j)}
+X_{-(\varepsilon_2-\varepsilon_j)}H_{\varepsilon_1+\varepsilon_{i+1}}\big)\otimes1_{-s}
\label{Eqn:Sp2}
\\
\phantom{\eqref{Eqn:Sp1}=} +\frac{1}{2}\sum_{j=i+2}^nN_{\varepsilon_1+\varepsilon_{i+1},-(\varepsilon_1+\varepsilon_j)}
\big(X_{-(\varepsilon_j-\varepsilon_{i+1})}X_{-(\varepsilon_2-\varepsilon_j)}
+X_{-(\varepsilon_2-\varepsilon_j)}X_{-(\varepsilon_j-\varepsilon_{i+1})}\big)\otimes1_{-s}.
\nonumber
\end{gather}
If $\lambda_i=\sum\limits_{j=1}^i\varepsilon_j$ is the fundamental weight for
$\alpha_i$ then, as $H \cdot 1_{-s} = \lambda_i(H) 1_{-s}$ for $H \in {\mathfrak h}$, a~direct computation
shows that
\begin{gather*}
\eqref{Eqn:Sp2}
=\frac{1}{2}\big(H_{\varepsilon_1+\varepsilon_{i+1}}X_{-(\varepsilon_2-\varepsilon_j)}
+X_{-(\varepsilon_2-\varepsilon_j)}H_{\varepsilon_1+\varepsilon_{i+1}}\big)\otimes1_{-s}
\\
\phantom{\eqref{Eqn:Sp2}=}
+\frac{1}{2}\sum_{j=i+2}^nN_{\varepsilon_1+\varepsilon_{i+1},-(\varepsilon_1+\varepsilon_j)}
\big(X_{-(\varepsilon_j-\varepsilon_{i+1})}X_{-(\varepsilon_2-\varepsilon_j)}
+X_{-(\varepsilon_2-\varepsilon_j)}X_{-(\varepsilon_j-\varepsilon_{i+1})}\big)\otimes1_{-s}
\\
\phantom{\eqref{Eqn:Sp2}}
=-\frac{1}{2}(\varepsilon_2-\varepsilon_{i+1})(H_{\varepsilon_1+\varepsilon_{i+1}}
)X_{-(\varepsilon_2-\varepsilon_{i+1})}\otimes1_{-s}-s\lambda_i(H_{\varepsilon_1+\varepsilon_{i+1}}
)X_{-(\varepsilon_2-\varepsilon_{i+1})}\otimes1_{-s}
\\
\phantom{\eqref{Eqn:Sp2}=}
+\frac{1}{2}\sum_{j=i+2}^n N_{\varepsilon_1+\varepsilon_{i+1},-(\varepsilon_1+\varepsilon_j)}
N_{\varepsilon_{i+1}-\varepsilon_j,-(\varepsilon_2-\varepsilon_j)}X_{-(\varepsilon_2-\varepsilon_{i+1})}
\otimes1_{-s}
\\
\phantom{\eqref{Eqn:Sp2}}
=-\frac{1}{2}(1)X_{-(\varepsilon_2-\varepsilon_{i+1})}\otimes1_{-s}
-s(1)X_{-(\varepsilon_2-\varepsilon_{i+1})}\otimes1_{-s}
\\
\phantom{\eqref{Eqn:Sp2}=}
+\frac{1}{2}\sum_{j=i+2}^n(1)(1)X_{-(\varepsilon_2-\varepsilon_{i+1})}\otimes1_{-s}
\\
\phantom{\eqref{Eqn:Sp2}}
=-\frac{1}{2}(2s-n+i)X_{-(\varepsilon_2-\varepsilon_{i+1})}\otimes1_{-s}.
\end{gather*}
Similarly, by a~direct computation, the second term amounts to
\begin{gather*}
T_2=-\sum_{j=i+1}^n\big(X_\mu\cdot\sigma(X_{-(\varepsilon_1-\varepsilon_j)}
X_{-(\varepsilon_2+\varepsilon_j)})\big)\otimes1_{-s}
\\
\phantom{T_2}
=-\frac{1}{2}N_{\varepsilon_1+\varepsilon_{i+1},-(\varepsilon_2+\varepsilon_{i+1})}
N_{\varepsilon_1-\varepsilon_2,-(\varepsilon_1-\varepsilon_{i+1})}X_{-(\varepsilon_2-\varepsilon_{i+1})}
\otimes1_{-s}
\\
\phantom{T_2=}
-\frac{1}{2}\sum_{j=i+1}^n N_{\varepsilon_1+\varepsilon_{i+1},-(\varepsilon_i-\varepsilon_j)}
N_{\varepsilon_{i+1}+\varepsilon_j,-(\varepsilon_2+\varepsilon_j)}X_{-(\varepsilon_2-\varepsilon_{i+1})}
\otimes1_{-s}
\\
\phantom{T_2}
=-\frac{1}{2}(1)(-1)X_{-(\varepsilon_2-\varepsilon_{i+1})}\otimes1_{-s}
\\
\phantom{T_2=}
-\frac{1}{2}(-2)(1)X_{-(\varepsilon_2-\varepsilon_{i+1})}\otimes1_{-s}-\frac{1}{2}\sum_{j=i+2}
^n(-1)(1)X_{-(\varepsilon_2-\varepsilon_{i+1})}\otimes1_{-s}
\\
\phantom{T_2}
=\frac{1}{2}(n-i+2)X_{-(\varepsilon_2-\varepsilon_{i+1})}\otimes1_{-s}.
\end{gather*}
Therefore, $X_{\mu} \cdot (\omega(Y^*_l)\otimes 1_{-s})$ is given by
\begin{gather*}
X_{\mu}\cdot(\omega\big(Y^*_l\big)\otimes1_{-s})=T_1+T_2
\\
\qquad =-\frac{1}{2}(2s-n+i)X_{-(\varepsilon_2-\varepsilon_{i+1})}\otimes1_{-s}+\frac{1}{2}
(n-i+2)X_{-(\varepsilon_2-\varepsilon_{i+1})}\otimes1_{-s}
\\
\qquad =-\frac{1}{2}(2s-n+i-(n-i+2))X_{-(\varepsilon_2-\varepsilon_{i+1})}\otimes1_{-s}
\\
\qquad =-(s-(n-i+1))X_{-(\varepsilon_2-\varepsilon_{i+1})}\otimes1_{-s}.
\end{gather*}
Hence $X_{\mu} \cdot (\omega(Y^*_l)\otimes 1_{-s}) = 0$ if and only if $s = n-i+1$.
\end{proof}

\subsection[The standardness of the map $\varphi_{\Omega_2}$]{The standardness of the map
$\boldsymbol{\varphi_{\Omega_2}}$}
%\label{SS:Type3D}

In the remainder of this section we determine the standardness of the map $\varphi_{\Omega_2}$ coming from
the conformally invariant $\Omega_2|_{V(\mu+\epsilon_\gamma)^*}$ system.

Observe that if $w_0$ is the longest Weyl group element for $\frak{l}_{\gamma}$ then the highest weight
$\nu$ for $V(\mu+\epsilon_\gamma)^*=V(\varepsilon_1+\varepsilon_2)^*$ is $\nu=-w_0(\varepsilon_1
+\varepsilon_2)=-\varepsilon_{i-1}-\varepsilon_i$.
By Theorem~\ref{Thm:SpValType3}, the special value~$s_2$ for the
$\Omega_2|_{V(\mu+\epsilon_{n\gamma})^*}$ system is $s_2 = n-i+1$.
Therefore, by~\eqref{Eqn:GVM2}, the $\Omega_2|_{V(\mu+\epsilon_{n\gamma})^*}$ system yields a~non-zero
$\mathcal{U}({\mathfrak g})$-homomorphism
\begin{gather}
\label{Eqn:Type3GVM}
\varphi_{\Omega_2}:M_{\mathfrak q}((-\varepsilon_{i-1}-\varepsilon_{i}
)-(n-i+1)\lambda_i+\rho)\to M_{\mathfrak q}(-(n-i+1)\lambda_i+\rho).
\end{gather}
\begin{Thm}
\label{Thm:Map_Type3}
If ${\mathfrak q}$ is the maximal parabolic subalgebra of type~$C_n(i)$ for $2 \leq i \leq n-1$ then the
standard map $\varphi_{std}$ between the generalized Verma modules in~\eqref{Eqn:Type3GVM} is zero.
Consequently, the map $\varphi_{\Omega_2}$ is non-standard.
\end{Thm}
\begin{proof}
By using the argument similar to one given in the proof for Theorem~\ref{Thm:Map_D(n-2)}, one can easily
show that $(\varepsilon_{i-1}-\varepsilon_i, \varepsilon_i-\varepsilon_n)$ links $-2\varepsilon_n -
(n-i+1)\lambda_i +\rho$ to $(-\varepsilon_{i=1}-\varepsilon_i)- (n-i+1)\lambda_i +\rho$.
Now the theorem follows from Proposition~\ref{Prop:Std}.
\end{proof}

\appendix

\section{Miscellaneous data}
\label{SS:Data}

This appendix summarizes the miscellaneous data for the maximal parabolic subalgebras ${\mathfrak
q}={\mathfrak l} \oplus {\mathfrak g}(1) \oplus {\mathfrak z}({\mathfrak n})$ of types $B_n(n-1)$, $C_n(i)$
($2\leq i \leq n-1$), and $D_n(n-2)$.
For the data for other maximal parabolic subalgebras of quasi-Heisenberg type see, for example, Appendix~A
of~\cite{KuboThesis1}.

\subsection[${B}_n(n-1)$]{$\boldsymbol{{B}_n(n-1)}$}

\begin{enumerate}\itemsep=0pt
\item The deleted Dynkin diagram:
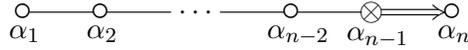
\begin{gather*}
\xymatrix{\belowwnode{\alpha_1}\ar@{-}[r]&
\belowwnode{\alpha_2}\ar@{-}[r]&\cdots\ar@{-}[r]
&\belowwnode{\alpha_{n-2}}\ar@{-}[r]
&\belowcnode{\alpha_{n-1}}\ar@2{->}[r]&\belowwnode{\alpha_n}}
\end{gather*}

\item The subgraph for $\frak{l}_{\gamma}$:
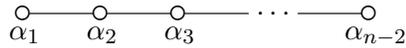
\begin{gather*}
\xymatrix{\belowwnode{\alpha_1}\ar@{-}[r]&
\belowwnode{\alpha_2}\ar@{-}[r]&
\belowwnode{\alpha_3}\ar@{-}[r]&\cdots\ar@{-}[r]
&\belowwnode{\alpha_{n-2}}}
\end{gather*}

\item The subgraph for $\frak{l}_{n\gamma}$:
\begin{gather*}
\xymatrix{\belowwnode{\alpha_n}}
\end{gather*}
\end{enumerate}

We have $\alpha_\gamma = \alpha_2$.
The highest weight $\mu$ and the set of weights $\Delta({\mathfrak g}(1))$ for ${\mathfrak g}(1)$ are $\mu
= \varepsilon_1+\varepsilon_{n}$ and $\Delta({\mathfrak g}(1)) = \{\varepsilon_j \pm \varepsilon_n \, | \,
1 \leq j \leq n-1 \} \cup \{\varepsilon_j \, | \, 1\leq j \leq n-1 \}$.
The highest weight $\gamma$ and the set of weights ${\mathfrak g}({\mathfrak z}({\mathfrak n}))$ for
${\mathfrak z}({\mathfrak n}))$ are $\gamma = \varepsilon_1 + \varepsilon_2$ and $\Delta({\mathfrak
z}({\mathfrak n})) = \{ \varepsilon_j + \varepsilon_k \, | \, 1\leq j<k\leq n-1\}$.
The highest root $\xi_\gamma$ and the set of positive roots $\Delta^+(\frak{l}_{\gamma})$ for
$\frak{l}_{\gamma}$ are $\xi_\gamma = \varepsilon_1 - \varepsilon_{n-1}$ and $\Delta^+(\frak{l}_{\gamma}) =
\{ \varepsilon_j - \varepsilon_k \, | \, 1\leq j < k\leq n-1\}$.
The highest root $\xi_{n\gamma}$ and the set of positive roots $\Delta^+(\frak{l}_{n\gamma})$ for
$\frak{l}_{n\gamma}$ are $\xi_{n\gamma} = \varepsilon_n$ and $\Delta^+(\frak{l}_{n\gamma}) = \{
\varepsilon_n \}$.

\subsection[${C}_n(i)$, $2\leq i \leq n-1$]{$\boldsymbol{{C}_n(i)}$, $\boldsymbol{2\leq i \leq n-1}$}

\begin{enumerate}\itemsep=0pt

\item The deleted Dynkin diagram:
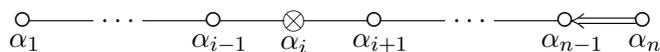
\begin{gather*}
\xymatrix{\belowwnode{\alpha_1}\ar@{-}[r]& \cdots\ar@{-}[r]& \belowwnode{\alpha_{i-1}}\ar@{-}
[r]& \belowcnode{\alpha_i}\ar@{-}[r]& \belowwnode{\alpha_{i+1}}\ar@{-}[r]& \cdots\ar@{-}
[r]& \belowwnode{\alpha_{n-1}}\ar@2{<-}[r]& \belowwnode{\alpha_n}}
\end{gather*}

\item The subgraph for $\frak{l}_{\gamma}$:
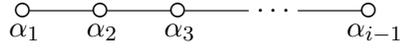
\begin{gather*}
\xymatrix{\belowwnode{\alpha_1}\ar@{-}[r]& \belowwnode{\alpha_2}\ar@{-}[r]& \belowwnode{\alpha_3}
\ar@{-}[r]& \cdots\ar@{-}[r]& \belowwnode{\alpha_{i-1}}}
\end{gather*}

\item The subgraph for $\frak{l}_{n\gamma}$:
\begin{gather*}
\xymatrix{\belowwnode{\alpha_{i+1}}\ar@{-}[r]& \cdots\ar@{-}[r]& \belowwnode{\alpha_{n-1}}
\ar@2{<-}[r]& \belowwnode{\alpha_n}}
\end{gather*}
\end{enumerate}

We have $\alpha_\gamma = \alpha_1$.
The highest weight $\mu$ and the set of weights $\Delta({\mathfrak g}(1))$ for ${\mathfrak g}(1)$ are $\mu
= \varepsilon_1+\varepsilon_{i+1}$ and $\Delta({\mathfrak g}(1)) = \{\varepsilon_j \pm \varepsilon_k \; |
\; 1 \leq j \leq i$ and $i+1 \leq k \leq n\}$.
The highest weight $\gamma$ and the set of weights $\Delta({\mathfrak z}({\mathfrak n}))$ for ${\mathfrak
z}({\mathfrak n})$ are $\gamma = 2\varepsilon_1$ $\Delta({\mathfrak z}({\mathfrak n})) = \{ \varepsilon_j +
\varepsilon_k \, | \, 1\leq j<k\leq i\} \cup \{2\varepsilon_j \, | \, 1 \leq j \leq i\}$.
The highest root $\xi_\gamma$ and the set of positive roots $\Delta^+(\frak{l}_{\gamma})$ for
$\frak{l}_{\gamma}$ are $\xi_\gamma = \varepsilon_1 - \varepsilon_i$ and $\Delta^+(\frak{l}_{\gamma}) = \{
\varepsilon_j - \varepsilon_k \, | \, 1\leq j < k\leq i\}$ The highest root $\xi_{n\gamma}$ and the set of
positive roots $\Delta(\frak{l}_{n\gamma})$ for $\frak{l}_{n\gamma}$ are $\xi_{n\gamma} =
2\varepsilon_{i+1}$ and $\Delta^+(\frak{l}_{n\gamma}) = \{ \varepsilon_j \pm \varepsilon_k \, | \, i+1 \leq
j < k \leq n \} \cup \{ 2\varepsilon_j \, | \, i+1 \leq j \leq n\}$.

\subsection[${D}_n(n-2)$]{$\boldsymbol{{D}_n(n-2)}$}

\begin{enumerate}\itemsep=0pt

\item The deleted Dynkin diagram:
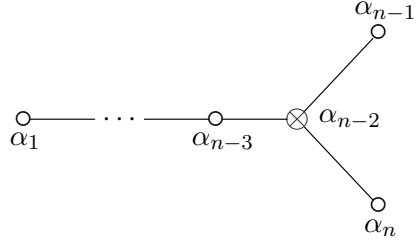
\begin{gather*}
\xymatrix{& & & & \abovewnode{\alpha_{n-1}}
\\
\belowwnode{\alpha_1}\ar@{-}[r]& \cdots\ar@{-}[r]& \belowwnode{\alpha_{n-3}}\ar@{-}
[r]& *-<\nodesize>{\crossandcirclesymbol}\save[]+<20pt,0pt>*\txt{$\alpha_{n-2}$}\restore\ar@{-}
[ur]\ar@{-}[dr]&
\\
& & & & \belowwnode{\alpha_n}}
\end{gather*}

\item The subgraph for $\frak{l}_{\gamma}$:
\begin{gather*}
\xymatrix{\belowwnode{\alpha_1}\ar@{-}[r]& \belowwnode{\alpha_2}\ar@{-}[r]& \belowwnode{\alpha_3}
\ar@{-}[r]& \cdots\ar@{-}[r]& \belowwnode{\alpha_{n-3}}}
\end{gather*}

\item The subgraph for $\frak{l}_{n\gamma}^-$:
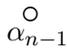
\begin{gather*}
\xymatrix{\belowwnode{\alpha_{n-1}}}
\end{gather*}

\item The subgraph for $\frak{l}_{n\gamma}^+$:
\begin{gather*}
\xymatrix{\belowwnode{\alpha_{n}}}
\end{gather*}

\end{enumerate}

We have $\alpha_\gamma = \alpha_2$.
The highest weight $\mu$ and the set of weights $\Delta({\mathfrak g}(1))$ for ${\mathfrak g}(1)$ are $\mu
= \varepsilon_1+\varepsilon_{n-1}$ and $\Delta({\mathfrak g}(1)) = \{e_j \pm e_k \, | \, 1 \leq j \leq n-2$
 and $k =n-1, n\}$.
The highest weight $\gamma$ and the set of weights $\Delta({\mathfrak z}({\mathfrak n}))$ for ${\mathfrak
z}({\mathfrak n})$ are $\gamma = \varepsilon_1 + \varepsilon_2$ and $\Delta({\mathfrak z}({\mathfrak n})) =
\{ e_j + e_k \, | \, 1\leq j<k\leq n-2\}$.
The highest root $\xi_\gamma$ and the set of positive roots $\Delta^+(\frak{l}_{\gamma})$ for
$\frak{l}_{\gamma}$ are $\xi_\gamma = \varepsilon_1 - \varepsilon_{n-2}$ and $\Delta^+(\frak{l}_{\gamma}) =
\{ e_j - e_k \, | \, 1\leq j < k\leq n-2\}$.
The highest root $\xi_{n\gamma}^-$ and the set of positive roots $\Delta^+(\frak{l}_{n\gamma}^-)$ are
$\xi_{n\gamma}^- = \varepsilon_{n-1} - \varepsilon_{n}$ and $\Delta^+(\frak{l}_{n\gamma}^-) = \{
\varepsilon_{n-1} - \varepsilon_{n} \}$.
The highest root $\xi_{n\gamma}^+$ and the set of positive roots $\Delta^+(\frak{l}_{n\gamma}^+)$ are
$\xi_{n\gamma}^+ = \varepsilon_{n-1} + \varepsilon_{n}$ and $\Delta^+(\frak{l}_{n\gamma}^+) = \{
\varepsilon_{n-1} + \varepsilon_{n} \}$.

\subsection*{Acknowledgements}
 The author was supported by the Global COE program at the Graduate School of
Mathematical Sciences, the University of Tokyo, Japan.
He would like to be thankful for the referees for their careful reading and invaluable comments.

\pdfbookmark[1]{References}{ref}
\LastPageEnding


\begin{thebibliography}{99}
\footnotesize\itemsep=0pt

\bibitem{BKZ08}
Barchini L., Kable A.C., Zierau R., Conformally invariant systems of
  dif\/ferential equations and prehomogeneous vector spaces of {H}eisenberg
  parabolic type, \href{http://dx.doi.org/10.2977/prims/1216238304}{\textit{Publ. Res. Inst. Math. Sci.}} \textbf{44} (2008),
  749--835.

\bibitem{BKZ09}
Barchini L., Kable A.C., Zierau R., Conformally invariant systems of
  dif\/ferential operators, \href{http://dx.doi.org/10.1016/j.aim.2009.01.006}{\textit{Adv. Math.}} \textbf{221} (2009), 788--811.

\bibitem{BE89}
Baston R.J., Eastwood M.G., The {P}enrose transform. Its interaction with
  representation theory, \textit{Oxford Mathematical Monographs, Oxford Science
  Publications}, The Clarendon Press, Oxford University Press, New York, 1989.

\bibitem{Boe85}
Boe B.D., Homomorphisms between generalized {V}erma modules, \href{http://dx.doi.org/10.2307/1999964}{\textit{Trans.
  Amer. Math. Soc.}} \textbf{288} (1985), 791--799.

\bibitem{BC86}
Boe B.D., Collingwood D.H., Intertwining operators between holomorphically
  induced modules, \href{http://dx.doi.org/10.2140/pjm.1986.124.73}{\textit{Pacific~J. Math.}} \textbf{124} (1986), 73--84.

\bibitem{Bourbaki08}
Bourbaki N., \'{E}l\'ements de math\'ematique. Groupes et
  alg{\`e}bres de Lie, Chapitres 4--6, Masson, Paris, 1981.

\bibitem{CS90}
Collingwood D.H., Shelton B., A duality theorem for extensions of induced
  highest weight modules, \href{http://dx.doi.org/10.2140/pjm.1990.146.227}{\textit{Pacific~J. Math.}} \textbf{146} (1990),
  227--237.

\bibitem{DES91}
Davidson M.G., Enright T.J., Stanke R.J., Dif\/ferential operators and highest
  weight representations, \href{http://dx.doi.org/10.1090/memo/0455}{\textit{Mem. Amer. Math. Soc.}} \textbf{94} (1991),
  no.~455, iv+102~pages.



\bibitem{Dobrev88}
Dobrev V., Canonical construction of dif\/ferential operators intertwining
  representations of real semisimple, \href{http://dx.doi.org/10.1016/0034-4877(88)90050-X}{\textit{Rep. Math. Phys.}} \textbf{25}
  (1988), 159--181.

\bibitem{Dobrev-a12}
Dobrev V., Invariant dif\/ferential operators for non-compact {L}ie algebras
  parabolically related to conformal {L}ie algebras, \href{http://dx.doi.org/10.1007/JHEP02(2013)015}{\textit{J.~High Energy
  Phys.}} \textbf{2013} (2013), no.~2, 015, 41~pages, \href{http://arxiv.org/abs/1208.0409}{arXiv:1208.0409}.

\bibitem{Dobrev-a13}
Dobrev V., Invariant dif\/ferential operators for non-compact {L}ie groups: the
  reduced ${\rm SU}(3,3)$ multiplets, \href{http://arxiv.org/abs/1312.5998}{arXiv:1312.5998}.

\bibitem{Ehrenpreis88}
Ehrenpreis L., Hypergeometric functions, in Algebraic Analysis, {V}ol.~{I},
  Academic Press, Boston, MA, 1988, 85--128.

\bibitem{HJ82}
Harris M., Jakobsen H.P., Singular holomorphic representations and singular
  modular forms, \href{http://dx.doi.org/10.1007/BF01457310}{\textit{Math. Ann.}} \textbf{259} (1982), 227--244.

\bibitem{Huang93}
Huang J.S., Intertwining dif\/ferential operators and reducibility of generalized
  {V}erma modules, \href{http://dx.doi.org/10.1007/BF01459504}{\textit{Math. Ann.}} \textbf{297} (1993), 309--324.

\bibitem{Jakobsen85}
Jakobsen H.P., Basic covariant dif\/ferential operators on {H}ermitian symmetric
  spaces, \textit{Ann. Sci. \'Ecole Norm. Sup.} \textbf{18} (1985), 421--436.

\bibitem{KKR03}
Johnson K.D., Kor{\'a}nyi A., Reimann H.M., Equivariant f\/irst order
  dif\/ferential operators for parabolic geometries, \href{http://dx.doi.org/10.1016/S0019-3577(03)90053-X}{\textit{Indag. Math.~(N.S.)}}
  \textbf{14} (2003), 385--393.

\bibitem{Juhl09}
Juhl A., Families of conformally covariant dif\/ferential operators,
  {$Q$}-curvature and holography, \href{http://dx.doi.org/10.1007/978-3-7643-9900-9}{\textit{Progress in Mathematics}}, Vol.~275,
  Birkh\"auser Verlag, Basel, 2009.

\bibitem{Kable11}
Kable A.C., {$K$}-f\/inite solutions to conformally invariant systems of
  dif\/ferential equations, \href{http://dx.doi.org/10.2748/tmj/1325886280}{\textit{Tohoku Math.~J.}} \textbf{63} (2011),
  539--559.

\bibitem{Kable12A}
Kable A.C., Conformally invariant systems of dif\/ferential equations on f\/lag
  manifolds for {$G_2$} and their {$K$}-f\/inite solutions, \textit{J.~Lie
  Theory} \textbf{22} (2012), 93--136.

\bibitem{Kable12C}
Kable A.C., The {H}eisenberg ultrahyperbolic equation: {$K$}-f\/inite and
  polynomial solutions, \textit{Kyoto~J. Math.} \textbf{52} (2012), 839--894.

\bibitem{Kable12B}
Kable A.C., The {H}eisenberg ultrahyperbolic equation: the basic solutions as
  distributions, \href{http://dx.doi.org/10.2140/pjm.2012.258.165}{\textit{Pacific~J. Math.}} \textbf{258} (2012), 165--197.

\bibitem{Kable13}
Kable A.C., On certain conformally invariant systems of dif\/ferential equations,
  \href{http://nyjm.albany.edu:8000/j/2013/19_189.html}{\textit{New York~J. Math.}} \textbf{19} (2013), 189--251.

\bibitem{Klimyk68}
Klimyk A.U., Decomposition of a direct product of irreducible representations
  of a semisimple Lie algebra into a direct sum of irreducible representations,
  in Thirteen Papers on Algebra and Analysis, \textit{Amer. Math. Soc. Translations}, Vol.~76,
Amer. Math. Soc., Providence, R.I., 1968, 63--74.

\bibitem{Knapp02}
Knapp A.W., Lie groups beyond an introduction, \textit{Progress in
  Mathematics}, Vol.~140, 2nd ed., Birkh\"auser Boston Inc., Boston, MA, 2002.

\bibitem{Knapp04}
Knapp A.W., Nilpotent orbits and some small unitary representations of
  indef\/inite orthogonal groups, \href{http://dx.doi.org/10.1016/S0022-1236(03)00254-4}{\textit{J.~Funct. Anal.}} \textbf{209} (2004),
  36--100.

\bibitem{KO03a}
Kobayashi T., {\O}rsted B., Analysis on the minimal representation of {${\rm
  O}(p,q)$}. {I}.~{R}ealization via conformal geometry, \href{http://dx.doi.org/10.1016/S0001-8708(03)00012-4}{\textit{Adv. Math.}}
  \textbf{180} (2003), 486--512, \href{http://arxiv.org/abs/math.RT/0111083}{math.RT/0111083}.

\bibitem{KO03b}
Kobayashi T., {\O}rsted B., Analysis on the minimal representation of {${\rm
  O}(p,q)$}. {II}.~{B}ranching laws, \href{http://dx.doi.org/10.1016/S0001-8708(03)00013-6}{\textit{Adv. Math.}} \textbf{180} (2003),
  513--550, \href{http://arxiv.org/abs/math.RT/0111085}{math.RT/0111085}.

\bibitem{KO03c}
Kobayashi T., {\O}rsted B., Analysis on the minimal representation of {${\rm
  O}(p,q)$}. {III}.~{U}ltrahyperbolic equations on {${\mathbb R}^{p-1,q-1}$},
  \href{http://dx.doi.org/10.1016/S0001-8708(03)00014-8}{\textit{Adv. Math.}} \textbf{180} (2003), 551--595, \href{http://arxiv.org/abs/math.RT/0111086}{math.RT/0111086}.

\bibitem{KOSS-a13}
Kobayashi T., {\O}rsted B., Somberg P., Soucek V., Branching laws for {V}erma
  modules and applications in parabolic geometry.~{I}, \href{http://arxiv.org/abs/1305.6040}{arXiv:1305.6040}.

\bibitem{KP-a13}
Kobayashi T., Pevzner M., {R}ankin--{C}ohen operators for symmetric pairs,
  \href{http://arxiv.org/abs/1301.2111}{arXiv:1301.2111}.

\bibitem{KR00}
Kor{\'a}nyi A., Reimann H.M., Equivariant f\/irst order dif\/ferential operators on
  boundaries of symmetric spaces, \href{http://dx.doi.org/10.1007/s002229900030}{\textit{Invent. Math.}} \textbf{139} (2000),
  371--390.

\bibitem{Kostant75}
Kostant B., Verma modules and the existence of quasi-invariant dif\/ferential
  operators, in Non-Commutative Harmonic Analysis ({A}ctes {C}olloq.,
  {M}arseille-{L}uminy, 1974), \textit{Lecture Notes in Math.}, Vol.~446,
  Springer, Berlin, 1975, 101--128.

\bibitem{Kubo11}
Kubo T., A system of third-order dif\/ferential operators conformally invariant
  under {$\mathfrak{sl}(3,{\mathbb C})$} and {$\mathfrak{so}(8,{\mathbb C})$},
  \href{http://dx.doi.org/10.2140/pjm.2011.253.439}{\textit{Pacific~J. Math.}} \textbf{253} (2011), 439--453, \href{http://arxiv.org/abs/1104.1999}{arXiv:1104.1999}.

\bibitem{KuboThesis2}
Kubo T., On the homomorphisms between the generalized {V}erma modules arising
  from conformally invariant systems, \textit{J.~Lie Theory} \textbf{23}
  (2013), 847--883, \href{http://arxiv.org/abs/1209.5516}{arXiv:1209.5516}.

\bibitem{KuboThesis1}
Kubo T., Special values for conformally invariant systems associated to maximal
  parabolics of quasi-{H}eisenberg type, \textit{Trans. Amer. Math. Soc.}, {t}o
  appear, \href{http://arxiv.org/abs/1209.1861}{arXiv:1209.1861}.

\bibitem{Lepowsky77}
Lepowsky J., A generalization of the {B}ernstein--{G}elfand--{G}elfand
  resolution, \href{http://dx.doi.org/10.1016/0021-8693(77)90254-X}{\textit{J.~Algebra}} \textbf{49} (1977), 496--511.

\bibitem{Matumoto06}
Matumoto H., The homomorphisms between scalar generalized {V}erma modules
  associated to maximal parabolic subalgebras, \href{http://dx.doi.org/10.1215/S0012-7094-05-13113-1}{\textit{Duke Math.~J.}}
  \textbf{131} (2006), 75--119, \href{http://arxiv.org/abs/math.RT/0309454}{math.RT/0309454}.

\bibitem{Matumoto-a12}
Matumoto H., On the homomorphisms between scalar generalized {V}erma modules,
  \href{http://arxiv.org/abs/1205.6748}{arXiv:1205.6748}.

\bibitem{Somberg-a13}
Somberg P., The {F}-method and a branching problem for generalized {V}erma
  modules associated to (Lie $G_2$, $so(7)$), \href{http://arxiv.org/abs/1303.7311}{arXiv:1303.7311}.

\bibitem{Wallach79}
Wallach N.R., The analytic continuation of the discrete series.~{II},
  \href{http://dx.doi.org/10.2307/1998680}{\textit{Trans. Amer. Math. Soc.}} \textbf{251} (1979), 19--37.

\end{thebibliography}
\end{document}